\documentclass{amsart}
\usepackage{geometry}
\usepackage[dvips]{graphicx}
\usepackage{amscd,amsfonts,amssymb,amsmath,amsthm,latexsym}

\usepackage{graphics}
\usepackage[all]{xy}
\xyoption{curve}
\xyoption{import}
\xyoption{arc}
\xyoption{ps}
\theoremstyle{plain}
    \newtheorem{teo}{Theorem}
    \newtheorem{lema}[teo]{Lemma}
    \newtheorem{coro}[teo]{Corollary}
    \newtheorem{prop}[teo]{Proposition}

    \newtheorem{conj}[teo]{Conjecture}
\theoremstyle{definition}
    \newtheorem{defi}[teo]{Definition}
    \newtheorem{ex}[teo]{Example}
\theoremstyle{remark}
    
    \newtheorem{obs}{Remark}
    
    \newtheorem{case}{Case}

\newcommand{\suchthat}{\ | \ }
\newcommand{\field}{K}
\newcommand{\ra}{R\langle\langle A\rangle\rangle}
\newcommand{\rap}{R\langle\langle A'\rangle\rangle}
\newcommand{\usualra}{R\langle A\rangle}
\newcommand{\idealM}{\mathfrak{m}}
\newcommand{\cyc}{\operatorname{cyc}}
\newcommand{\racyc}{R\langle\langle A\rangle\rangle_{\cyc}}
\newcommand{\jacobas}{\mathcal{P}(A,S)}
\newcommand{\jacobapsp}{\mathcal{P}(A',S')}

\newcommand{\astrivial}{(A_{\operatorname{triv}},S_{\operatorname{triv}})}
\newcommand{\asreduced}{(A_{\operatorname{red}},S_{\operatorname{red}})}
\newcommand{\atriv}{A_{\operatorname{triv}}}
\newcommand{\ored}{\operatorname{red}}
\newcommand{\otriv}{\operatorname{triv}}
\newcommand{\surf}{(\Sigma, M)}
\newcommand{\punct}{P}
\newcommand{\marked}{M}
\newcommand{\arc}{i}
\newcommand{\arcsinsurf}{\mathbf{A}^\circ\surf}

\newcommand{\qtau}{Q(\tau)}
\newcommand{\qsigma}{Q(\sigma)}
\newcommand{\atau}{A(\tau)}
\newcommand{\asigma}{A(\sigma)}
\newcommand{\ratau}{R\langle\langle A(\tau)\rangle\rangle}

\newcommand{\stau}{S(\tau)}
\newcommand{\astau}{(A(\tau),S(\tau))}
\newcommand{\qstau}{(Q(\tau),S(\tau))}
\newcommand{\qssigma}{(Q(\sigma),S(\sigma))}
\newcommand{\assigma}{(A(\sigma),S(\sigma))}

\newcommand{\rasigma}{R\langle\langle A(\sigma)\rangle\rangle}
\newcommand{\ssigma}{S(\sigma)}
\newcommand{\sptau}{S'(\tau)}
\newcommand{\spsigma}{S'(\sigma)}
\newcommand{\surfnoM}{\Sigma}
\newcommand{\premuti}{\widetilde{\mu}_i}
\newcommand{\muti}{\mu_i}
\newcommand{\tildeas}{(\widetilde{A},\widetilde{S})}
\newcommand{\tildeastau}{(\widetilde{A(\tau)},\widetilde{S(\tau)})}
\newcommand{\tildeatau}{\widetilde{A(\tau)}}
\newcommand{\tildestau}{\widetilde{S(\tau)}}
\newcommand{\rtildeatau}{R\langle\langle\widetilde{A(\tau)}\rangle\rangle}

\newcommand{\overatau}{\overline{A(\tau)}}
\newcommand{\roveratau}{R\langle\langle\overline{(A(\tau))}\rangle\rangle}
\newcommand{\overqtau}{\overline{Q(\tau)}}
\newcommand{\overstau}{\overline{S(\tau)}}
\newcommand{\overqstau}{(\overline{Q(\tau)},\overline{S(\tau)})}
\newcommand{\surfpn}{(\Sigma, M\cup P_n)}
\newcommand{\surfpnminone}{(\Sigma, M\cup P_{n-1})}
\newcommand{\surfpone}{(\Sigma, M\cup P_1)}
\newcommand{\punctn}{P_n}
\newcommand{\qsigmaone}{Q(\sigma_1)}
\newcommand{\qsigmatwo}{Q(\sigma_2)}
\newcommand{\qsigmanminone}{Q(\sigma_{n-1})}
\newcommand{\qtaun}{Q(\tau_n)}
\newcommand{\qssigmanminone}{(Q(\sigma_{n-1}),S(\sigma_{n-1}))}
\newcommand{\qstaun}{(Q(\tau_n),S(\tau_n))}

\newcommand{\trspace}{\operatorname{Tr}}
\newcommand{\premutwo}{\widetilde{\mu}_2}
\newcommand{\mutwo}{\mu_2}
\newcommand{\unredqtau}{\widehat{Q(\tau)}}
\newcommand{\unredstau}{\widehat{S}(\tau)}
\newcommand{\runredatau}{R\langle\langle\widehat{A}(\tau)\rangle\rangle}
\newcommand{\unredastau}{(\widehat{A}(\tau),\widehat{S}(\tau))}
\newcommand{\unredassigma}{(\widehat{A}(\sigma),\widehat{S}(\sigma))}

\begin{document}

\title{Quivers with potentials associated to triangulated surfaces}
\author{Daniel Labardini-Fragoso}
\address{Department of Mathematics, Northeastern University, Boston, MA 02115}
\email{labardini-fra.d@neu.edu}
\date{\today}
\thanks{This work was partially supported by Prof. Andrei Zelevinsky's NSF grant and Prof. Jos\'e Antonio de la Pe\~na's SNI grant.}
\maketitle

\begin{abstract}
We attempt to relate two recent developments: cluster algebras associated to triangulations of
surfaces by Fomin-Shapiro-Thurston, and quivers with potentials and their mutations introduced by Derksen-Weyman-Zelevinsky. To each ideal triangulation of a bordered surface with marked points we associate a quiver with potential, in such a way that whenever two ideal triangulations are related by a flip of an arc, the respective quivers with potentials are related by a mutation with respect to the flipped arc. We prove that if the surface has non-empty boundary, then the quivers with potentials associated to its triangulations are rigid and hence non-degenerate.
\end{abstract}

\tableofcontents

%main text:
\section{Introduction}
This work is inspired by two beautiful papers: On the one hand, in \cite{FST}, S.
Fomin, M. Shapiro and D. Thurston associate to each bordered surface with marked points a cluster algebra, each of whose exchange matrices is
defined in terms of the (signed) adjacencies between the arcs of a triangulation of the surface. They prove that the seeds of this cluster algebra are
related by a mutation if and only if the triangulations to which the seeds are associated are related by a flip. In particular, if two
triangulations are related by a flip, then the (skew-symmetric) matrices associated to them are related by a mutation.

On the other hand, in \cite{DWZ}, H. Derksen, J. Weyman and A. Zelevinsky introduce the
notion of quivers with potentials (QPs for short), that is, pairs consisting of a quiver and a special element of its (complete) path
algebra, and define the mutations of such objects, ultimately leading to the notion of mutation of representations, thus providing a new
representation-theoretic interpretation for quiver mutations originated in the theory of cluster algebras, interpretation that generalizes the
classical Bernstein-Gelfand-Ponomarev reflection functors.

Here we make a first attempt to relate the two above mentioned papers. We associate to each triangulation $\tau$ of a bordered surface with
marked points a potential $S(\tau)$ on the quiver $Q(\tau)$ defined by its signed adjacency matrix $B(\tau)$. The idea is quite simple: each
interior triangle of $\tau$ gives rise to an oriented triangle in $Q(\tau)$, and each puncture has an oriented cycle of $Q(\tau)$ around it;
what we do is to add such oriented triangles and cycles to get the potential $S(\tau)$. We then extend a result from
\cite{FST}: not only the quivers $\qtau$ and $\qsigma$ are related by a mutation if the ideal triangulations $\tau$ and $\sigma$ are related by a
flip, but also the QPs $\qstau$ and $\qssigma$.

There is a delicate point in the process of QP-mutation: The underlying quiver of the mutated QP
depends on the potential of the original QP. This implies that the underlying quiver of the mutated QP may not coincide with the quiver
obtained under ``ordinary" quiver mutation from the underlying quiver of the original QP. The difficulty relies on the fact that it is the
potential that allows us (or not) to delete the 2-cycles from the quiver. A QP is non-degenerate if, after any sequence of QP-mutations, the
potential always allows us to delete all the 2-cycles of the corresponding quiver.

A family of QPs for which non-degeneracy is guaranteed without the necessity of checking all possible sequences of QP-mutations is the family
of rigid QPs. A QP is rigid if every cycle on the underlying quiver is cyclically equivalent to an element of the ideal generated by the cyclic
derivatives of the underlying potential. But even to decide whether a given QP is rigid or not may be a difficult task, as it is to give an
algorithm to decide it. Here we explicitly define a family of rigid QPs associated to surfaces with non-empty boundary.

We now describe the contents of the paper in more detail. Section \ref{background} is divided into three parts. In the first part we recall the definition and involutive properties of mutations of quivers. In the second part we describe the combinatorial setup for this paper, namely, the notions of surface with marked points, arcs, ideal triangulations and their flips. We also introduce the quivers associated to ideal triangulations, and mention the compatibility between the mutation of such quivers and the flips of the triangulations. For technical reasons, we define also the \emph{unreduced} version of these quivers as the result of adding certain 2-cycles to them.

In the third part of Section \ref{background} we summarize the basics of the theory of QP-mutations: we define the notions of complete path algebra, potential, cyclical and right equivalence, reduced and trivial parts of a QP, mutation of a QP with respect to a vertex, and restriction of a QP to a subset of the vertex set. We also prove a result that takes place in the general theory of QP-mutations (and not only in the surface-related setup), namely, that the operation of restriction commutes with the operation of QP-mutation (a similar result was established in \cite{DWZ} for rigid QPs).

In Section \ref{QPofatriangulation} we associate a QP $(Q(\tau),S(\tau))$ to each ideal triangulation $\tau$ of a surface $\surfnoM$ with marked points $M$. This QP is defined as the reduced part of a QP on the arrow span of the unreduced signed adjacency quiver $\unredqtau$. The quiver $Q(\tau)$ turns out to be the quiver associated to (the signed adjacency matrix $B(\tau)$ of) $\tau$ in \cite{FST}. After defining $(Q(\tau),S(\tau))$, we prove that ideal triangulations related by a flip give raise to QPs related by QP-mutation. Unfortunately, the proof of this fact is done with an analysis case-by-case, the reason being that slight changes in the configuration of the arcs surrounding the arc to be flipped can dramatically affect the associated QP.

In Section \ref{rigidityfornonemptyboundary} we prove that the QPs associated to surfaces with non-empty boundary are rigid and hence non-degenerate. Therefore, since the class of quivers associated to ideal triangulations is closed under mutation, we will have given an explicit construction of a non-degenerate QP for each of the quivers in this class, and moreover, since ideal triangulations related by a flip give rise to QPs related by a mutation, these QPs will represent a well defined QP-mutation class for the surface. We conjecture that the QPs we associate to surfaces with empty boundary are non-degenerate as well, but non-rigid. We close Section \ref{rigidityfornonemptyboundary} showing that if the surface has non-empty boundary, then the Jacobian algebras of the QPs we assign to its triangulations are finite-dimensional.

Thanks to talks of T. Br\"{u}stle and R. Schiffler (\cite{B} and \cite{Sicra}) at the International Conference on Representations of Algebras (ICRA XII) held in Toru\'n, Poland, in August 2007, the author became aware of some independent works (\cite{ABCP}, \cite{CCS04}, \cite{PMaster}, \cite{S}) about cluster-tilted algebras and gentle algebras that arise as quivers with relations associated to surface triangulations. Those works deal with bordered surfaces with marked points in a more restricted set up than ours (\cite{CCS04} and \cite{S} deal  mainly with unpunctured and once-punctured disks, while \cite{ABCP} and \cite{PMaster} considers only unpunctured surfaces), and it turns out that the Jacobian algebras of our QPs specialize to their constructions. Namely, the relations they get are precisely the result of applying the cyclic derivatives to the potentials we give here. However, they focus on other interesting and important properties of the corresponding algebras that we had not considered initially. In the forthcoming continuation \cite{L} of this paper we will study the QP-mutation behavior of a class of representations whose beautiful construction was presented by T. Br\"{u}stle in \cite{B} for surfaces without punctures; we will study also some generalizations of this construction to the case of punctured surfaces and relate them to the cluster algebras associated to the surfaces.

\section{Background on triangulations of surfaces and quivers with potentials}\label{background}

\subsection{Quiver mutations}

Recall that a \emph{quiver} is a finite directed graph, that is, a quadruple $Q=(Q_0,Q_1,h,t)$, where $Q_0$ is the (finite) set of
\emph{vertices} of $Q$, $Q_1$ is the (finite) set of \emph{arrows}, and $h:Q_1\rightarrow Q_0$ and $t:Q_1\rightarrow Q_0$ are the \emph{head}
and \emph{tail} functions. Recall also the common notation $a:i\rightarrow j$ to indicate that $a$ is an arrow of $Q$ with $t(a)=i$, $h(a)=j$.
We will always deal only with loop-free quivers, that is, with quivers that have no arrow $a$ with $t(a)=h(a)$.

A \emph{path of length} $d>0$ in $Q$ is a sequence $a_1a_2\ldots a_d$ of arrows with $t(a_j)=h(a_{j+1})$ for $j=1,\ldots,d-1$. A path
$a_1a_2\ldots a_d$ of length $d>0$ is a $d$\emph{-cycle} if $h(a_1)=t(a_d)$. A quiver is \emph{2-acyclic} if it has no 2-cycles.

Paths are composed as functions, that is, if $a=a_1\ldots a_d$ and $b=b_1\ldots b_{d'}$ are paths with $h(b)=t(a)$, then the product $ab$ is defined as the path $a_1,\ldots a_db_1\ldots b_{d'}$, which starts at $t(b_{d'})$ and ends at $h(a_1)$. See Figure \ref{prodofpaths}.

 \begin{figure}[!h]
                \caption{Paths are composed as functions}\label{prodofpaths}
                \centering
$$
\bullet\overset{b_{d'}}{\longrightarrow}\ldots\overset{b_1}{\longrightarrow}\bullet\overset{a_d}{\longrightarrow}\ldots\overset{a_1} {\longrightarrow}\bullet
$$
 \end{figure}

For $i\in Q_0$, an $i$\emph{-hook} in $Q$ is any path $ab$ of length 2 such that $a,b\in Q_1$ are arrows with $t(a)=i=h(b)$.

\begin{defi}\label{threesteps} Given a quiver $Q$ and a vertex $i\in Q_0$ such that $Q$ has no $2$-cycles incident at $i$, we define the
\emph{mutation} of $Q$ in direction $i$ as the quiver $\muti(Q)$ with vertex set $Q_0$ that results after applying the following three-step
procedure:
\begin{itemize}
\item[(Step 1)] For each $i$-hook $ab$ introduce an arrow $[ab]:t(b)\rightarrow h(a)$.
\item[(Step 2)] Replace each arrow $a:i\rightarrow h(a)$ of $Q$ with an arrow $a^*:h(a)\rightarrow i$, and each arrow $b:t(b)\rightarrow i$
of $Q$ with an arrow $b^*:i\rightarrow t(b)$.
\item[(Step 3)] Choose a maximal collection of disjoint 2-cycles and remove them.
\end{itemize}
We call the quiver obtained after the $1^{st}$ and $2^{nd}$ steps the \emph{premutation} $\premuti(Q)$.
\end{defi}

\begin{obs}\begin{enumerate}\item Note that the mutation $\muti$ is defined for non-necessarily 2-acyclic quivers, but in order to be able to perform mutation at any vertex of a quiver, we need it to be 2-acyclic.
\item The choice of the maximal collection in the $3^{rd}$ step is not given by a canonical procedure. However, up to this choice, $\muti$ is an involution on the class of $2$-acyclic quivers, that is, $\muti^2(Q)\cong Q$ for every 2-acyclic quiver $Q$.
\end{enumerate}
\end{obs}

\subsection{Triangulations of surfaces and their flips}\label{backontriangs}

In this subsection we briefly review the material on triangulations of surfaces and their signed adjacency matrices and flips. The reader will find much
deeper discussions on the subject in \cite{FST}.

\begin{defi}[\cite{FST}, Definition 2.1] A \emph{bordered surface with marked points} is a pair $\surf$, where $\surfnoM$ is a compact connected oriented Riemann surface
with (possibly empty) boundary, and $\marked$ is a finite set of points on $\surfnoM$, called \emph{marked points}, such that $\marked$ is
non-empty and has at least one point from each connected component of the boundary of $\surfnoM$. The marked points that lie in the interior of
$\surfnoM$ will be called \emph{punctures}, and the set of punctures of $\surf$ will be denoted $\punct$. We will always assume that $\surf$ is
none of the following:
\begin{itemize}
\item a sphere with less than five punctures;
\item an unpunctured monogon, digon or triangle;
\item a once-punctured monogon.
\end{itemize}
Here, by a monogon (resp. digon, triangle) we mean a disk with exactly one (resp. two, three) marked point(s) on the boundary.
\end{defi}

\begin{obs} The reason for excluding the surfaces in the second and third items of the above definition is the fact that their triangulations (in the sense of Definition \ref{idealtriangdef} below) are empty or there is only one such. The reason for excluding the spheres with less than five punctures is a bit more subtle (self-folded triangles present some ``unpleasant" properties in these surfaces).
\end{obs}

\begin{defi}[\cite{FST}, Definition 2.2] Let $\surf$ be a bordered surface with marked points. An (\emph{ordinary}) \emph{arc} in $\surf$ is a
curve $\arc$ in $\surfnoM$ such that:
\begin{itemize}
\item the endpoints of $\arc$ are marked points in $\marked$;
\item $\arc$ does not intersect itself, except that its endpoints may coincide;
\item the relative interior of $\arc$ is disjoint from $\marked$ and from the boundary of $\surfnoM$;
\item $\arc$ does not cut out an unpunctured monogon or an unpunctured digon.
\end{itemize}
We consider two arcs $\arc_1$ and $\arc_2$ to be the same whenever they are isotopic in $\surfnoM$ rel $\marked$, that is whenever there exists an isotopy $H:I\times\surfnoM\rightarrow\surfnoM$ such that $H(0,x)=x$ for all $x\in\surfnoM$, $H(1,\arc_1)=\arc_2$, and $H(t,m)=m$ for all $t\in I$ and all $m\in\marked$. An arc whose endpoints coincide will be called
a \emph{loop}. (Do not confuse the notion of loop in $\surf$ with a loop in a quiver). We denote the set of (isotopy classes of) arcs in
$\surf$ by $\arcsinsurf$.
\end{defi}

Two arcs are \emph{compatible} if there are arcs in their respective isotopy classes whose relative interiors do not intersect (cf. \cite{FST}, Definition 2.4).

\begin{prop} Given any collection of pairwise compatible arcs, it is always possible to find representatives in their isotopy classes
whose relative interiors do not intersect each other.
\end{prop}

\begin{defi}\label{idealtriangdef} An \emph{ideal triangulation} of $\surf$ is any maximal collection of pairwise compatible arcs whose relative interiors do not intersect each
other (cf. \cite{FST}, Definition 2.6).
\end{defi}

If $\tau$ is an ideal triangulation of $\surf$ and we take a connected component of the complement in $\surfnoM$ of the union of the
arcs in $\tau$, the closure $\triangle$ of this component will be called an \emph{ideal triangle} of $\tau$. An ideal triangle $\triangle$ is
called \emph{interior} if its intersection with the boundary of $\surfnoM$ consists only of (possibly none) marked points. Otherwise it will be
called \emph{non-interior}. An interior ideal triangle $\triangle$ is \emph{self-folded} if it contains exactly two arcs of $\tau$ (note that
every interior ideal triangle contains at least two and at most three arcs of $\tau$, while each non-interior ideal triangle contains at least
one and at most two arcs).

% SELF-FOLDED TRIANGLE
        \begin{figure}[!h]
                \caption{Self-folded triangle}\label{selffoldedtriang}
                \centering
                \includegraphics[scale=.4]{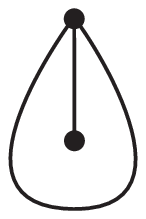}
        \end{figure}

The number $n$ of arcs in an ideal triangulation of $\surf$ is determined by the genus $g$ of $\surfnoM$, the number $b$ of boundary components
of $\surfnoM$, the number $p$ of punctures and the number $c$ of marked points on the boundary of $\surfnoM$, according to the formula
$n=6g+3b+3p+c-6$, which can be proved using the definition and basic properties of the Euler characteristic. Hence $n$ is an invariant of
$\surf$, called the \emph{rank} of $\surf$ (because it coincides with the rank of the cluster algebra associated to $\surf$, see \cite{FST}).

Let $\tau$ be an ideal triangulation of $\surf$ and let $\arc\in\tau$ be an arc. If $\arc$ is not the folded side of a self-folded triangle,
then there exists exactly one arc $\arc'$, different from $\arc$, such that $\sigma=(\tau\setminus\{\arc\})\cup\{\arc'\}$ is an ideal
triangulation of $\surf$. We say that $\sigma$ is obtained by applying a \emph{flip} to $\tau$, or by \emph{flipping} the arc $\arc$ (cf. \cite{FST}, Definition 3.5), and write
$\sigma=f_i(\tau)$.
In order to be able to flip the folded sides of self-folded triangles, one has to enlarge the set of arcs with which
triangulations are formed. This is done by introducing the notion of \emph{tagged arc}. Since we will deal only with ordinary arcs in this paper, we refer the reader to \cite{FST} and \cite{FT} for the definition and properties of tagged arcs and \emph{tagged triangulations}.

\begin{prop}[\cite{FST}, Propositions 3.8 and 7.10.]\label{seqofflips} Any two ideal triangulations are related by a sequence of flips. If $\surf$ is not a surface with empty boundary and exactly one puncture, then any two tagged triangulations are related by a sequence of flips.
\end{prop}

To each ideal triangulation $\tau$ we associate a skew-symmetric $n\times n$ integer matrix $B(\tau)$ whose rows and columns correspond to the
arcs of $\tau$ (cf. \cite{FST}, Definition 4.1). Let $\pi_\tau:\tau\rightarrow\tau$ be the function that is the identity on the set of arcs
that are not folded sides of self-folded triangles of $\tau$, and sends the folded side of a self-folded triangle to the unique loop of $\tau$
enclosing it. For each non-self-folded ideal triangle $\triangle$ of $\tau$, let $B^\triangle=b^\triangle_{ij}$ be the $n\times n$ integer
matrix defined by
\begin{equation}
b^\triangle_{ij}=
\begin{cases} 1\ \ \ \ \ \text{if $\triangle$ has sides $\pi_\tau(i)$ and $\pi_\tau(j)$, with $\pi_\tau(j)$ following
$\pi_\tau(i)$}\\ \ \ \ \ \ \ \ \ \ \text{in the clockwise order defined by the orientation of $\surfnoM$;}\\ -1\ \ \ \text{if the same holds,
but in the counter-clockwise
order;}\\
0\ \ \ \ \ \text{otherwise.}
\end{cases}
\end{equation}
The \emph{signed adjacency matrix} $B(\tau)$ is then defined as
\begin{equation}
B(\tau)=\underset{\triangle}{\sum}B^\triangle,
\end{equation}
where the sum runs over all non-self-folded triangles of $\tau$. Note that $B(\tau)$ is skew-symmetric, and all its entries have absolute value
less than 3.
%
%To define the signed adjacency matrix of a tagged triangulation, we need the notion of signature of a tagged triangulation (cf. \cite{FST},
%Definition 9.1). Let $\tau$ be a tagged triangulation of $\surf$. The \emph{signature} of $\tau$ is the function
%$\delta_\tau:P\rightarrow\{-1,0,1\}$ defined by
%\begin{equation}
%\delta_\tau(p)=
%\begin{cases} 1\ \ \ \ \ \text{if all ends incident to $p$ are tagged plain;}\\ -1\ \ \ \text{if all ends incident to $p$ are tagged notched;}\\
%0\ \ \ \ \ \text{otherwise.}
%\end{cases}
%\end{equation}
%Note that if $\delta_\tau(p)=0$, then there are precisely two arcs in $\tau$ incident to $p$, the untagged versions of these arcs coincide and
%they carry the same tag at the end different from $p$. We define $B(\tau)=B(\tau^\circ)$, where $\tau^\circ$ is the ideal triangulation
%obtained from $\tau$ according to the following rules:
%\begin{itemize}
%\item delete all tags at the punctures $p$ with $\delta_\tau(p)=\pm1$;
%\item for each puncture $p$ with $\delta_\tau(p)=0$, replace the arc $\arc$ notched at $p$ by a loop enclosing $p$ and $\arc$. See Figure \ref{signaturezero}.
%\end{itemize}
%
%% SIGNATURE ZERO
%        \begin{figure}[!h]
%                \caption{}\label{signaturezero}
%                \centering
%                \includegraphics[scale=.3]{signaturezero.EPS}
%        \end{figure}

The matrix $B(\tau)$ gives rise to the \emph{signed adjacency quiver} $\qtau$, whose vertices are the arcs in $\tau$, with $b_{ij}$ arrows from
$i$ to $j$ whenever $b_{ij}>0$. Since $B(\tau)$ is skew-symmetric, $\qtau$ is a 2-acyclic quiver.

\begin{teo}[\cite{FST}, Proposition 4.8] Let $\tau$ and $\sigma$ be ideal triangulations. If $\sigma$ is obtained from $\tau$
by flipping the arc $\arc$ of $\tau$, then $Q(\sigma)=\muti(\qtau)$.
\end{teo}

\begin{obs}\begin{enumerate}\item This Theorem holds in the more general situation where $\tau$ and $\sigma$ are tagged triangulations related by a flip, see \cite{FST}, Lemma 9.7.
\item The assignment [skew-symmetric $n\times n$ integer matrix $B$] $\mapsto$ [2-acyclic quiver $Q$] is general and defines a
bijection between skew-symmetric $n\times n$ integer matrices and 2-acyclic quivers on $n$ vertices. See \cite{FZ2} for far-reaching
discussions on this matter.
\item Signed adjacency matrices for triangulations had already appeared in \cite{FG} and \cite{GSV}, though in a more restricted setup.
\end{enumerate}
\end{obs}

\begin{defi} If a puncture is incident to exactly two arcs $\arc_1$ and $\arc_2$ of the ideal triangulation $\tau$, then $\qtau$ has no arrows between $\arc_1$ and $\arc_2$ (see Figure \ref{noarrows}). For each such pair of arcs we add to $\qtau$ an arrow from $\arc_1$ to $\arc_2$ and an arrow from $\arc_2$ to $\arc_1$, and call the resulting quiver the \emph{unreduced signed adjacency quiver} $\unredqtau$.
% NO ARROWS IN A DIGON
        \begin{figure}[!h]
                \caption{}\label{noarrows}
                \centering
                \includegraphics[scale=.5]{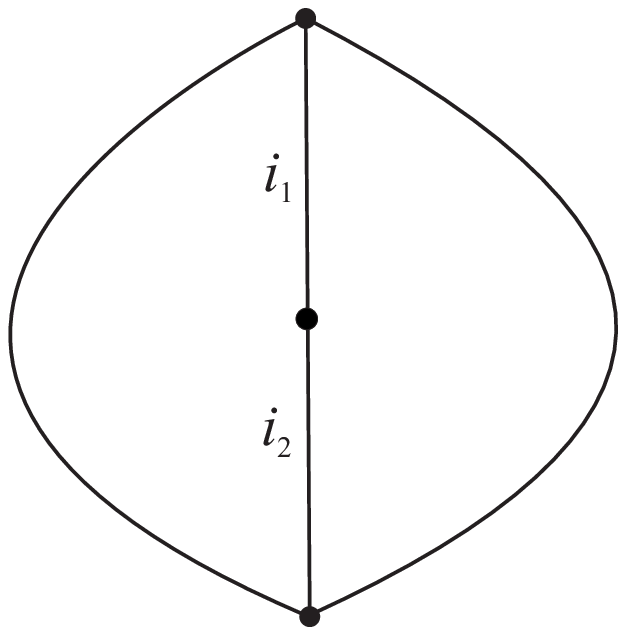}
        \end{figure}
\end{defi}

It is clear that $\qtau$ can be obtained from $\unredqtau$ by deleting all 2-cycles.

\begin{ex} In Figure \ref{squaretriangs} we can see some ideal triangulations of the once-punctured square and their signed adjacency quivers drawn on them.
% TRIANGULATIONS OF THE ONCE-PUNCTURED SQUARE
        \begin{figure}[!h]
                \caption{Some signed adjacency quivers}\label{squaretriangs}
                \centering
                \includegraphics[scale=.4]{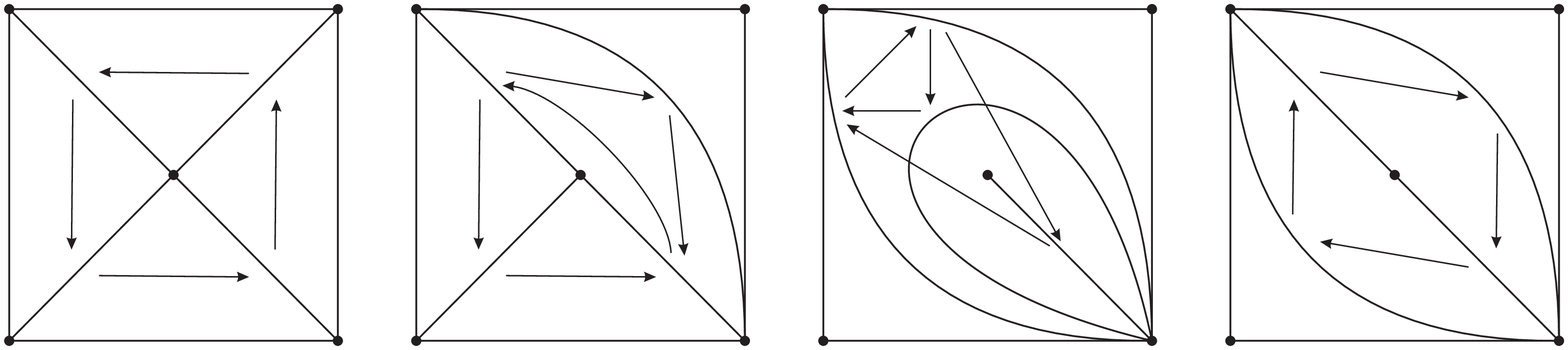}
        \end{figure}
In Figure \ref{squaretriangsunred} we have the same triangulations, but their unreduced signed adjacency quivers instead.
% TRIANGULATIONS OF THE ONCE-PUNCTURED SQUARE (UNREDUCED VERSION)
        \begin{figure}[!h]
                \caption{Some unreduced signed adjacency quivers}\label{squaretriangsunred}
                \centering
                \includegraphics[scale=.4]{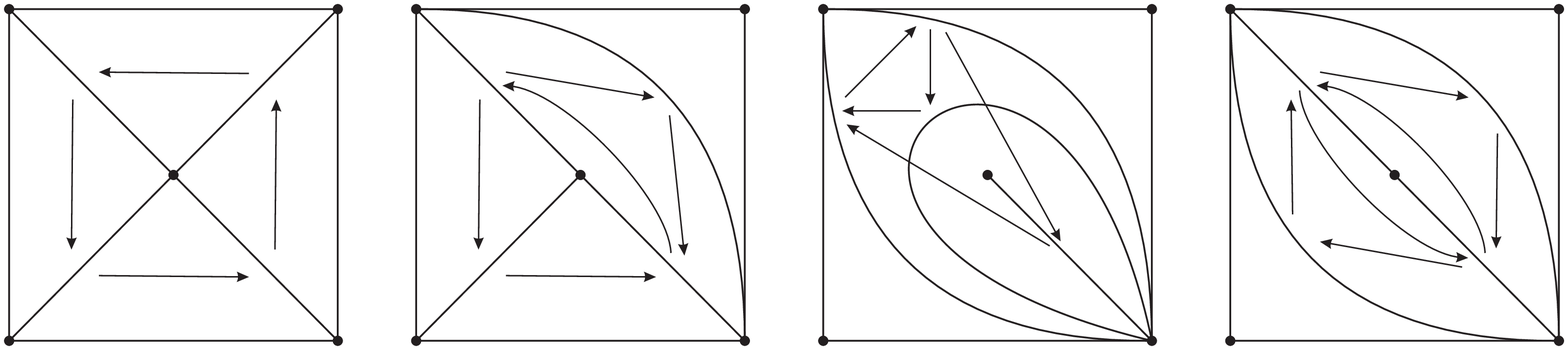}
        \end{figure}
Notice that among these four triangulations, the only one for which the signed adjacency quiver is different from the unreduced one is the triangulation appearing farthest right in both Figures \ref{squaretriangs} and \ref{squaretriangsunred}.
\end{ex}

\subsection{Quivers with potentials and their mutations}

In this subsection we give the background on quivers with potentials and their mutations we shall use in the remaining sections. For a more
detailed and elegant treatment of the subject, we refer the reader to \cite{DWZ}, all of whose notation we will adopt here. In particular,
$\field$ will always denote a field. A survey of the topics treated in \cite{DWZ} can be found in \cite{Z}.

Given a quiver $Q$, we denote by $R$ the $\field$-vector space with basis $\{e_i\suchthat i\in Q_0\}$. If we define $e_ie_j=\delta_{ij}e_i$,
then $R$ becomes naturally a commutative semisimple $\field$-algebra, which we call the \emph{vertex span} of $Q$; each $e_i$ is called the
\emph{path of length zero} at $i$. We define the \emph{arrow span} of $Q$ as the $\field$-vector space $A$ with basis the set of arrows $Q_1$.
Note that $A$ is an $R$-bimodule if we define $e_ia=\delta_{i,h(a)}a$ and $ae_j=a\delta_{t(a),j}$ for $i\in Q_0$ and $a\in Q_1$. For $d\geq 0$
we denote by $A^d$ the $\field$-vector space with basis all the paths of length $d$ in $Q$; this space has a natural $R$-bimodule structure as well. Notice that $A^0=R$ and $A^1=A$.

The \emph{complete path algebra} of $Q$ is the $\field$-vector space consisting of all possibly infinite linear combinations of paths in $Q$,
that is,
\begin{equation}
\ra=\underset{d=0}{\overset{\infty}{\prod}}A^d;
\end{equation}
with multiplication induced by concatenation of paths (cf. \cite{DWZ}, Definition 2.2). Note that $\ra$ is a $\field$-algebra and an $R$-bimodule, and has the usual \emph{path
algebra}
\begin{equation}
\usualra=\underset{d=0}{\overset{\infty}{\bigoplus}}A^d
\end{equation}
as $\field$-subalgebra and sub-$R$-bimodule. Moreover, $\usualra$ is dense in $\ra$ under the $\idealM$-adic topology, whose fundamental system
of open neighborhoods around $0$ is given by the powers of $\idealM=\idealM(A)=\underset{d\geq 1}{\prod}A^d$, the ideal of $\ra$ generated by the arrows. A crucial property of this
topology is the following:
\begin{equation}
\text{a sequence $(x_n)_{n\in\mathbb{N}}$ of elements of $\ra$ converges if and only if for every $d\geq 0$},
\end{equation}
\begin{center}
the sequence $(x_n^{(d)})_{n\in\mathbb{N}}$ stabilizes as $n\rightarrow\infty$, in which case
$\underset{n\rightarrow\infty}{\lim}x_n=\underset{d\geq 0}{\sum}\underset{n\rightarrow\infty}{\lim}x_n^{(d)}$,
\end{center}
where $x_n^{(d)}$ denotes the degree-$d$ component of $x_n$.

Even though the action of $R$ on $\ra$ (and $\usualra$) is not central, it is compatible with the multiplication of $\ra$ in the sense that if
$a$ and $b$ are paths in $Q$, then $e_{h(a)}ab=ae_{t(a)}b=abe_{t(b)}$. Therefore we will say that $\ra$ (and $\usualra$) are
$R$-algebras. Accordingly, any $K$-algebra homomorphism $\varphi$ between (complete) path algebras will be called an $R$-algebra homomorphism
if the underlying quivers have the same set of vertices and $\varphi(r)=r$ for every $r\in R$. It is easy to see that every $R$-algebra
homomorphism between complete path algebras is continuous. The following is an extremely useful criterion to decide if a given linear map
$\varphi:\ra\rightarrow\rap$ between complete path algebras (on the same set of vertices) is an $R$-algebra homomorphism or an $R$-algebra
isomorphism:
\begin{equation}
\text{Every pair $(\varphi^{(1)},\varphi^{(2)})$ of $R$-bimodule homomorphisms $\varphi^{(1)}:A\rightarrow A'$, $\varphi^{(2)}:A\rightarrow\idealM(A')^2$,}
\end{equation}
\begin{center}
extends uniquely to a continuous $R$-algebra homomorphism $\varphi:\ra\rightarrow\rap$ such that $\varphi|_A=(\varphi^{(1)},\varphi^{(2)})$.
Furthermore, $\varphi$ is $R$-algebra isomorphism if and only if $\varphi^{(1)}$ is an $R$-bimodule isomorphism.
\end{center}

A \emph{potential} on $A$ (or $Q$) is any element of $\ra$ all of whose terms are cyclic paths of positive length (cf. \cite{DWZ}, Definition 3.1). The set of all potentials on
$A$ is denoted by $\racyc$, it is a closed vector subspace of $\ra$. Two potentials $S,S'\in\racyc$ are \emph{cyclically equivalent} if $S-S'$
lies in the closure of the vector subspace of $\ra$ spanned by all the elements of the form $a_1\ldots a_d-a_2\ldots a_da_1$ with $a_1\ldots
a_d$ a cyclic path of positive length (cf. \cite{DWZ}, Definition 3.2).

A \emph{quiver with potential} is a pair $(A,S)$ (or $(Q,S)$), where $S$ is a potential on $A$ such that no two different cyclic paths
appearing in the expression of $S$ are cyclically equivalent (cf. \cite{DWZ}, Definition 4.1). (For instance, the pair $(A,xa_1\ldots a_d-ya_{i+1}\ldots a_da_1\ldots a_i)$ is
not a quiver with potential for any choice of different non-zero scalars $x,y\in\field$). We will use the shorthand QP to abbreviate ``quiver
with potential". The \emph{direct sum} of two QPs $(A,S)$ and $(A',S')$ on the same set of vertices is the QP $(A,S)\oplus(A',S')=(A\oplus A',S+S')$.

If $(A,S)$ and $(A',S')$ are QPs on the same set of vertices, we say that $(A,S)$ is \emph{right-equivalent} to $(A',S')$ if there exists a
\emph{right-equivalence} between them, that is, an $R$-algebra isomorphism $\varphi:\ra\rightarrow\rap$ such that $\varphi(S)$ is cyclically
equivalent to $S'$ (cf. \cite{DWZ}, Definition 4.2).

For each arrow $a\in Q_1$ and each cyclic path $a_1\ldots a_d$ in $Q$ we define the \emph{cyclic derivative}
\begin{equation}
\partial_a(a_1\ldots a_d)=\underset{i=1}{\overset{d}{\sum}}\delta_{a,a_i}a_{i+1}\ldots a_da_1\ldots a_{i-1},
\end{equation}
(where $\delta_{a,a_i}$ is the \emph{Kronecker delta}) and extend $\partial_a$ by linearity and continuity to obtain a map $\partial_a:\racyc\rightarrow\ra$ (cf. \cite{DWZ}, Definition 3.1). Note that we have
$\partial_a(S)=\partial_a(S')$ whenever the potentials $S$ and $S'$ are cyclically equivalent.

The \emph{Jacobian ideal} $J(S)$ is the closure of the two-sided ideal of $\ra$ generated by $\{\partial_a(S)\suchthat a\in Q_1\}$, and the
\emph{Jacobian algebra} $\jacobas$ is the quotient algebra $\ra/J(S)$ (cf. \cite{DWZ}, Definition 3.1). Jacobian ideals and Jacobian algebras are invariant under
right-equivalences, in the sense that if $\varphi:\ra\rightarrow\rap$ is a right-equivalence between $(A,S)$ and $(A',S')$, then $\varphi$
sends $J(S)$ onto $J(S')$ and therefore induces an isomorphism $\jacobas\rightarrow\jacobapsp$ (cf. \cite{DWZ}, Proposition 3.7).

%A QP $(Q,S)$ is \emph{trivial} if $Q$ has exactly $2N$ arrows (for some $N>0$) that can be listed in the form $a_1,b_1,\ldots,a_N,b_N$ with
%$a_jb_j$ a $2$-cycle for $j=1,\ldots,N$, and there is an $R$-algebra automorphism $\varphi:\ra\rightarrow\ra$ such that $\varphi(S)$ is
%cyclically equivalent to $a_1b_1+\ldots+a_Nb_N$ (cf. \cite{DWZ}, Definition 4.3 and Proposition 4.4).
A QP $(Q,S)$ is \emph{trivial} if $S\in A^2$ and $\{\partial_a(S)\suchthat a\in Q_1\}$ spans $A$ (cf. \cite{DWZ}, Definition 4.3, see also Proposition 4.4 therein). We say that a QP $(A,S)$ is \emph{reduced} if the degree-$2$ component of $S$ is $0$, that is,
if the expression of $S$ involves no $2$-cycles. Note that the underlying quiver of a reduced QP may have 2-cycles. We say that a quiver $Q$
(or its arrow span, or any QP on it) is \emph{2-acyclic} if it has no $2$-cycles.

\begin{teo}[Splitting Theorem, \cite{DWZ}, Theorem 4.6]\label{splittingthm} For every QP $(A,S)$ there exist a trivial QP
$\astrivial$ and a reduced QP $\asreduced$ such that $(A,S)$ is right-equivalent to the direct sum $\astrivial\oplus\asreduced$. Furthermore,
the right-equivalence class of each of the QPs $\astrivial$ and $\asreduced$ is determined by the right-equivalence class of $(A,S)$.
\end{teo}

In the situation of Theorem \ref{splittingthm}, the QP $\asreduced$ (resp. $\astrivial$) is called the \emph{reduced part} (resp. \emph{trivial
part}) of $(A,S)$ (cf. \cite{DWZ}, Definition 4.13); this terminology is well defined up to right-equivalence.

We now turn to the definition of mutation of a QP. Let $(A,S)$ be a QP on the vertex set $Q_0$ and let $i\in Q_0$. Assume that $Q$ has no
2-cycles incident to $i$. Thus, if necessary, we replace $S$ with a cyclically equivalent potential so that we can assume that every cyclic
path appearing in the expression of $S$ does not begin at $i$. This allows us to define $[S]$ as the potential on $\widetilde{A}$ obtained from
$S$ by replacing each $i$-hook $ab$ with the arrow $[ab]$ (see Definition \ref{threesteps}). Also, we define $\Delta_i(Q)=\sum b^*a^*[ab]$,
where the sum runs over all $i$-hooks $ab$ of $Q$.

\begin{defi}[\cite{DWZ}, equations (5.3) and (5.8) and Definition 5.5] Under the assumptions and notation just stated, we define the \emph{premutation} of $(A,S)$ in direction $i$ as the QP
$\premuti(A,S)=\tildeas$, where $\widetilde{A}$ is the arrow span of $\premuti(Q)$ (see Definition \ref{threesteps}) and
$\widetilde{S}=[S]+\Delta_i(Q)$. The \emph{mutation} $\muti(A,S)$ of $(A,S)$ in direction $i$ is then defined as the reduced part of
$\premuti(A,S)$.
\end{defi}

\begin{teo}[\cite{DWZ}, Theorem 5.2 and Corollary 5.4] Premutations and mutations are well defined up to right-equivalence. That is, if $(A,S)$ and $(A',S')$ are right-equivalent QPs with no 2-cycles incident to the vertex $i$, then the QP $\premuti(A,S)$ is right-equivalent to $\premuti(A',S')$ and the QP $\muti(A,S)$ is right-equivalent to $\muti(A',S')$.
\end{teo}

\begin{teo}[\cite{DWZ}, Theorem 5.7] Mutations are involutive up to right-equivalence. More specifically, if $(A,S)$ is a QP such that $A$ is the arrow span of a quiver that has no 2-cycles incident to the vertex $i$, then $\muti^2(A,S)$ is right-equivalent to $(A,S)$.
\end{teo}

``Unfortunately'', as the following easy example shows, 2-acyclicity is \textbf{not} a QP-mutation invariant (in contrast to the ordinary quiver mutation, where 2-acyclicity is ensured by definition).

\begin{ex} Consider the potentials $S_1=0$ and $S_2=abc$ on the quiver
\begin{displaymath}
\xymatrix{ & 2 \ar[dr]^b & \\
1 \ar[ur]^c & & 3 \ar[ll]^a }
\end{displaymath}
If we perform the premutation $\premutwo$ on $(A,S_1)$ and $(A,S_2)$, we get $(\widetilde{A},\widetilde{S_1})$ and $(\widetilde{A},\widetilde{S_2})$, where $\widetilde{A}$ is the arrow span of the quiver
\begin{displaymath}
\xymatrix{ & 2 \ar[dl]_{c^*} & \\
1 \ar@/^/[rr]^{[bc]} & & 3 \ar[ul]_{b^*} \ar@/^/[ll]^a }
\end{displaymath}
and $\widetilde{S_1}=c^*b^*[bc]$, $\widetilde{S_2}=a[bc]+c^*b^*[bc]$. Note that $(\widetilde{A},\widetilde{S_1})$ is already reduced (hence equals $\mutwo(A,S_1)$), while the reduced part of $(\widetilde{A},\widetilde{S_2})$ is $(\overline{A},0)$ (hence $\mutwo(A,S_2)=(\overline{A},0)$), where $\overline{A}$ is the arrow span of the quiver
\begin{displaymath}
\xymatrix{ & 2 \ar[dl]_{c^*} & \\
1 & & 3 \ar[ul]_{b^*} }
\end{displaymath}
In particular, we cannot apply the mutations in direction 1 or 3 to $(\widetilde{A},\widetilde{S_1})$, but we can apply them to $(\overline{A},0)$.
\end{ex}

\begin{defi}[\cite{DWZ}, Definition 7.2] A QP $(A,S)$ is \emph{non-degenerate} if it is 2-acyclic and the quiver of the QP obtained after any possible sequence of QP-mutations is 2-acyclic.
\end{defi}

\begin{teo}[\cite{DWZ}, Proposition 7.3 and Corollary 7.4] If the base field $\field$ is uncountable, then every 2-acyclic quiver admits a
non-degenerate QP.
\end{teo}

A QP $(A,S)$ is \emph{rigid} if every cycle in $Q$ is cyclically equivalent to an element of the Jacobian ideal $J(S)$ (cf. \cite{DWZ}, Definition 6.10 and equation 8.1). Rigidity is invariant under QP-mutation.

\begin{teo}[\cite{DWZ}, Corollary 6.11, Proposition 8.1 and Corollary 8.2] Every reduced rigid QP is 2-acyclic. The class of reduced rigid QPs is closed under QP-mutation. Consequently, every rigid reduced QP is non-degenerate.
\end{teo}

\begin{prop}[\cite{DWZ}, Corollary 6.6]\label{findimisinvariant} Let $(A,S)$ be a non-degenerate QP and $i\in Q_0$ any vertex, then the Jacobian algebra $\mathcal{P}(A,S)$ is finite-dimensional if and only if so is $\mathcal{P}(\muti(A,S))$. In other words, finite-dimensionality of Jacobian algebras is invariant under QP-mutations.
\end{prop}

We finish this section describing the operation of \emph{restriction} of a QP to a subset of the set of vertices.

\begin{defi}[\cite{DWZ}, Definition 8.8]\label{restriction} Let $(A,S)$ be a QP and $I$ be a subset of the vertex set $Q_0$. The \emph{restriction} of $(A,S)$ to $I$ is the QP $(A|_I,S|_I)$ on the vertex set $Q_0$, with $A|_I=\underset{i,j\in I}{\bigoplus}A_{ij}$ and $S|_I=\psi_I(S)$, where $\psi_I:R\langle\langle A\rangle\rangle\rightarrow R\langle\langle A|_I\rangle\rangle$ is the $R$-algebra homomorphism  such that $\psi_I(a)=a$ for $a\in A|_I$ and $\psi_I(b)=0$ for each arrow $b\notin A|_I$.
\end{defi}

\begin{obs}\label{remarkisolated} Notice that if $I$ is a proper subset of $Q_0$, then the elements of $Q_0\setminus I$ are isolated vertices of the restriction to $I$, that is, there are no arrows of $A|_I$ whose head or tail belongs to $Q_0\setminus I$.
\end{obs}

In Proposition 8.9 of \cite{DWZ}, Derksen-Weyman-Zelevinsky prove that restriction preserves rigidity and finite-dimensionality of Jacobian algebras. Here we prove that it preserves non-degeneracy as well. As a preparation for the proof we recall the construction, given in \cite{DWZ} to prove the Splitting Theorem, of a reduced part and a trivial part of a QP. Let $(A,S)$ be any QP and denote by $S^{(2)}$ the degree-2 component of
$S$. Assume $S^{(2)}\neq 0$ (otherwise $(A,S)$ is already reduced and there is no construction to be done). Up to a right-equivalence that acts as the identity on the arrows of $Q$ not appearing in $S^{(2)}$, we can
assume that
\begin{equation}\label{formofS}
S=\underset{j=1}{\overset{N}{\sum}}(a_jb_j+a_ju_j+v_jb_j)+S',
\end{equation}
where each $a_jb_j$ is a 2-cycle, the $2N$ distinct arrows $a_1,b_1,\ldots,a_N,b_N$ form a basis of $\atriv=\partial S^{(2)}$, each of $u_j$
and $v_j$ belongs to $\idealM^2$, and $S'\in\idealM^3$ is a potential neither of whose terms involves any of the arrows $a_j$ or $b_j$. If
$u_j=v_j=0$ for all $j$, then we already have the decomposition of Theorem \ref{splittingthm}. Otherwise, one defines a unitriangular
automorphism $\varphi_1:\ra\rightarrow\ra$ by setting
\begin{equation}\label{rulesofre}
\varphi_1(a_j)=a_j-v_j,\ \ \varphi_1(b_j)=b_j-u_j,\ \ \varphi_1(c)=c\text{ for }c\in Q_1\setminus\{a_1,b_1,\ldots,a_N,b_N\}.
\end{equation}
This automorphism of $\ra$ is a right-equivalence between $(A,S)$ and $(A,S_1)$, where $S_1$ is a potential with $S_1^{(2)}=S^{(2)}$ and such
that, when written in the form (\ref{formofS}), all the $u_j$-factors and $v_j$-factors belong to $\idealM^3$. If all these factors are $0$, we
already have reached the decomposition of Theorem \ref{splittingthm}. Otherwise, we define a unitriangular automorphism $\varphi_2$ of $\ra$ by
the rules (\ref{rulesofre}) defined in terms of the $u_j$-factors and $v_j$-factors of $S_1$. Then $\varphi_2\varphi_1$ is a right-equivalence
between $(A,S)$ and $(A,S_2)$, where $S_2$ is a potential with $S_2^{(2)}=S_1^{(2)}=S^{(2)}$ and such that, when written in the form
(\ref{formofS}), the $u_j$-factors and $v_j$-factors belong to $\idealM^5$.

If we keep repeating the above procedure \emph{ad infinitum}, we get a sequence $(S_n)_{n\geq1}$, with the corresponding $u_j$-factors and $v_j$-factors belonging to higher and higher powers of $\idealM$. In the limit,
$\varphi=\underset{n\rightarrow\infty}{\lim}\varphi_n\ldots\varphi_1$ will be a right-equivalence between $(A,S)$ and
$(A,\underset{n\rightarrow\infty}{\lim}S_n)$, where $\underset{n\rightarrow\infty}{\lim}S_n$ is a potential whose degree-2 component is
$S^{(2)}$ and such that, when written in the form (\ref{formofS}), all its $u_j$-factors and $v_j$-factors are $0$. This provides the required
right-equivalence of Theorem \ref{splittingthm}.

\begin{obs} In many concrete examples (like the ones given in the present work), there is no need of considering the above limit process.
\end{obs}

\begin{lema}\label{redres=resred} Let $(A,S)$ be a QP, and let $I$ be any subset of the vertex set $Q_0$. There exist a reduced and a trivial QP, $\asreduced$ and $\astrivial$, respectively, such that $(A,S)$ is right-equivalent to $\asreduced\oplus\astrivial$, and with the property that the restriction $(A|_I,S|_I)$ is right-equivalent to $(A_{\operatorname{red}}|_I,S_{\operatorname{red}}|_I)\oplus(A_{\operatorname{triv}}|_I,S_{\operatorname{triv}}|_I)$.
\end{lema}

\begin{proof} The vector subspace of $A$ generated by the cyclic derivatives of the degree-2 component of $S$ is $A_{\otriv}$, and we clearly have $(A|_I)_{\operatorname{triv}}=A_{\operatorname{triv}}|_I$ (see \cite{DWZ}, equations 4.3 and 4.4). Therefore we will also have $(A|_I)_{\operatorname{red}}=A_{\operatorname{red}}|_I$. Denote these spaces by $B_{\otriv}=A_{\operatorname{triv}}|_I$ and $B_{\ored}=A_{\operatorname{red}}|_I$, respectively, and let $T_{\ored}\in R\langle\langle B_{\ored}\rangle\rangle$ and $T_{\otriv}\in R\langle\langle B_{\otriv}\rangle\rangle$ be potentials such that $(B_{\ored},T_{\ored})$ is a reduced QP, $(B_{\otriv},T_{\otriv})$ is a trivial QP, and there exists an $R$-algebra isomorphism $\varphi: R\langle\langle A|_I\rangle\rangle\rightarrow R\langle\langle B_{\ored}\oplus B_{\otriv}\rangle\rangle$ such that $\varphi(S|_I)$ is cyclically equivalent to $T_{\ored}+T_{\otriv}$. Also, let us write $S=S|_I+S'$, where $S'\in R\langle\langle A\rangle\rangle$ is a potential each of whose terms involves at least one arrow that does not belong to $A|_I$.

We can extend $\varphi$ to an $R$-algebra isomorphism $\widehat{\varphi}:R\langle\langle A\rangle\rangle\rightarrow R\langle\langle A_{\ored}\oplus A_{\otriv}\rangle\rangle$ by defining $\widehat{\varphi}(b)=b$ for $b\notin A|_I$. The potential $\widehat{\varphi}(S)$ is cyclically equivalent to $T_{\ored}+T_{\otriv}+\widehat{\varphi}(S')$. Let us denote the degree-2 component of $\widehat{\varphi}(S')$ by $T_{\otriv}'$. Note that every term of the potential $\widehat{\varphi}(S')$ involves at least one arrow that does not belong to $A|_I$; in particular, every arrow appearing in $T_{\otriv}'$ is incident to a vertex outside $I$. Note also that the arrows that appear in $T_{\otriv}$ can appear in $\widehat{\varphi}(S')$, but do not appear in $T_{\ored}$. These remarks make it clear that if we decompose the QP $(A_{\ored}\oplus A_{\otriv},T_{\ored}+T_{\otriv}+\widehat{\varphi}(S'))$ as the direct sum of a reduced and a trivial QP according to the procedure described before this lemma, we will have $S_n|_I=T_{\ored}+T_{\otriv}$ for all $n\geq 1$. Therefore, the restriction of the reduced part of $(A,S)$ to $I$ is (right-equivalent to) the reduced part of the restriction of $(A,S)$ to $I$.
\end{proof}

\begin{prop}\label{resmut=mutres} Let $(A,S)$ be a QP and $I$ a subset of $Q_0$. For $i\in I$, the mutation $\muti(A|_I,S|_I)$ is right-equivalent to the restriction of $\muti(A,S)$ to $I$.
\end{prop}

\begin{proof} An easy check shows that $([S]+\triangle_i(Q))|_I=[S|_I]+\triangle_i(Q|_I)$. Therefore, the premutation $\premuti(A|_I,S|_I)$ is equal to the restriction of $\premuti(A,S)$ to $I$. The proposition then follows from Lemma \ref{redres=resred}.
\end{proof}

\begin{coro} If $(A,S)$ is a non-degenerate QP, then for every subset $I$ of $Q_0$ the restriction $(A|_I,S|_I)$ is non-degenerate as well. In other words, restriction preserves non-degeneracy.
\end{coro}

\section{The QP of a triangulation}\label{QPofatriangulation}

Let $\surf$ be a bordered surface with marked points, with $\punct\subseteq\marked$ the set of punctures of $\surf$. For each $p\in\punct$
choose a non-zero scalar $x_p\in \field$; this choice is going to kept fixed for every triangulation of $\surf$.

\begin{defi}\label{QPfortriangulation} Let $\tau$ be an ideal triangulation of $\surf$. Based on our choice $(x_p)_{p\in\punct}$ we associate to $\tau$ a potential $\stau\in\ratau$ as follows. Let $\widehat{A}(\tau)$ denote the arrow span of $\widehat{Q(\tau)}$.
\begin{itemize} \item For each interior non-self-folded ideal triangle $\triangle$ of $\tau$ that gives rise to an oriented triangle of $\widehat{\qtau}$, let $\widehat{S}^\triangle$ be such oriented triangle up to cyclical equivalence.
\item If the interior non-self-folded ideal triangle $\triangle$ with sides $j$, $k$ and $l$, is adjacent to two self-folded triangles like in the configuration of Figure \ref{adsftriangs},
% TWO ADJACENT SELF-FOLDED TRIANGLES
        \begin{figure}[!h]
                \caption{}\label{adsftriangs}
                \centering
                \includegraphics[scale=.5]{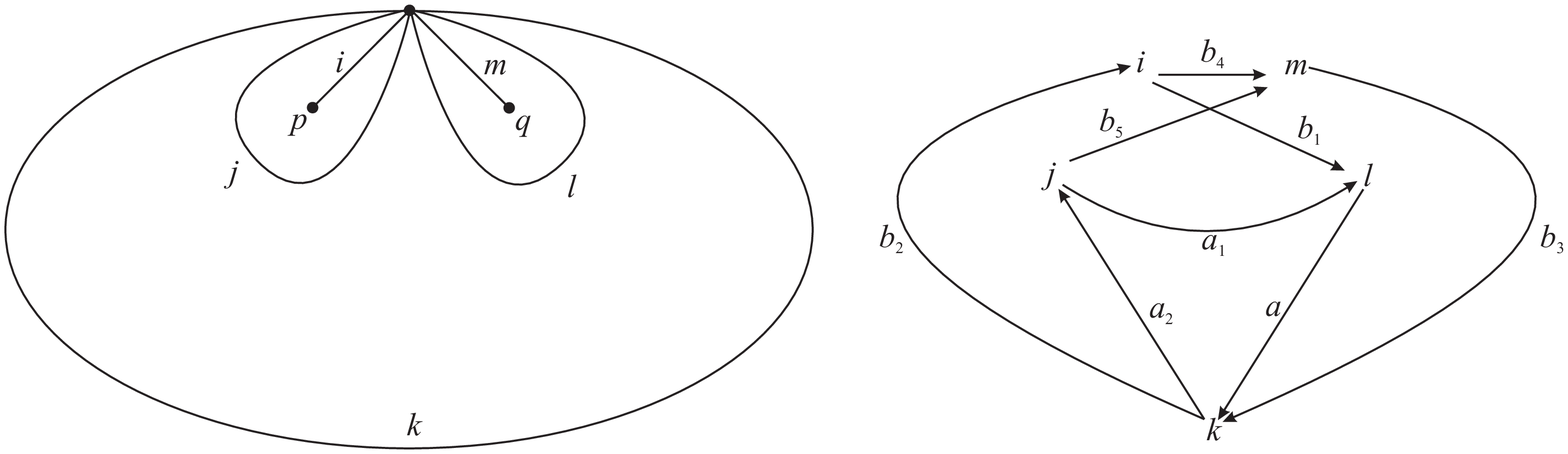}
        \end{figure}\\
define $\widehat{T}^\triangle=\frac{b_2b_3b_4}{x_px_q}$ (up to cyclical equivalence), where $p$ and $q$ are the punctures enclosed in the self-folded triangles adjacent to $\triangle$. Otherwise, if it is adjacent to less than two self-folded triangles, define
$\widehat{T}^\triangle=0$.
\item If a puncture $p$ is adjacent to exactly one arc $i$ of $\tau$, then $i$ is the folded side of a self-folded triangle of $\tau$ and around $i$ we have the configuration shown in Figure \ref{sftriangle}.
% SELF-FOLDED TRIANGLE
        \begin{figure}[!h]
                \caption{}\label{sftriangle}
                \centering
                \includegraphics[scale=.5]{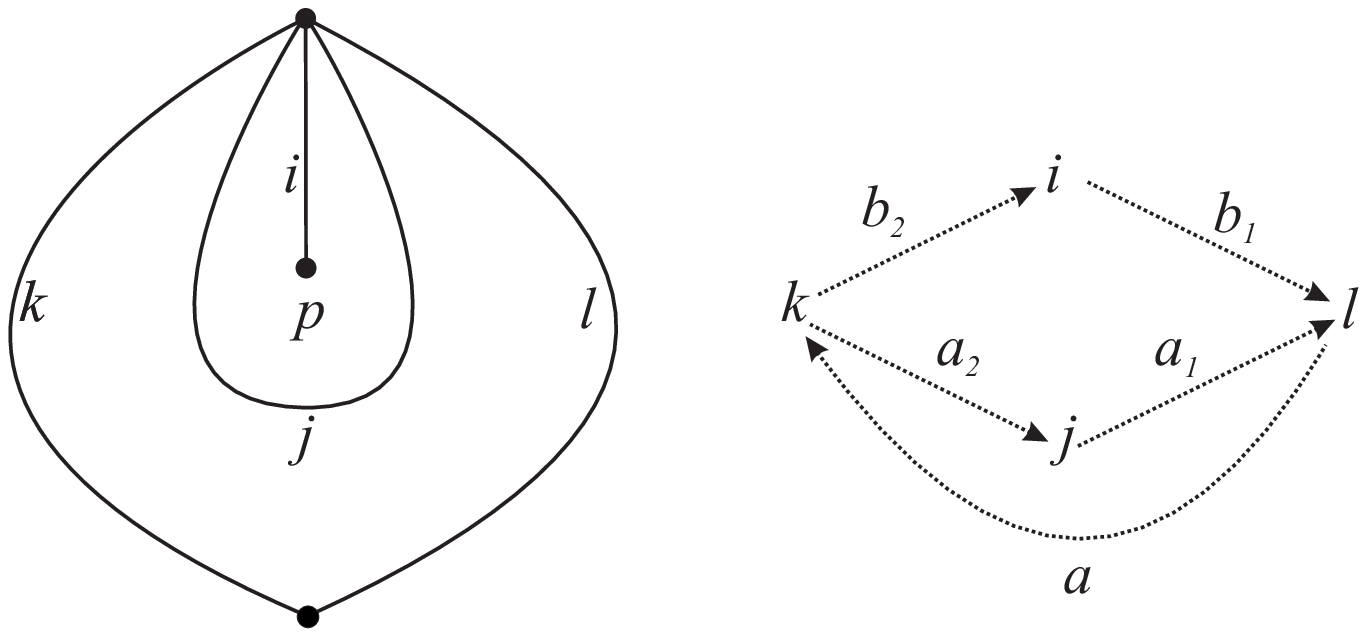}
        \end{figure}\\
In case both $k$ and $l$ are indeed arcs of $\tau$ (and not part of the boundary of $\surfnoM$), then we define $\widehat{S}^p=-\frac{ab_1b_2}{x_p}$ (up to cyclical equivalence).
\item If a puncture $p$ is adjacent to more than one arc, delete all the loops adjacent to $p$ that enclose a self-folded triangle. The arrows
between the remaining arcs adjacent to $p$ form a unique cycle $a^p_1\ldots a^p_d$ that exhausts all such remaining arcs and gives a complete round around $p$ in the counter-clockwise orientation defined by the orientation of $\surfnoM$. We define $\widehat{S}^p=x_pa^p_1\ldots a^p_d$ (up to cyclical equivalence).
\end{itemize}
The \emph{unreduced potential} $\unredstau\in\runredatau$ of $\tau$ is then defined by
\begin{equation}
\unredstau=\underset{\triangle}{\sum}(\widehat{S}^\triangle+\widehat{T}^\triangle)+\underset{p\in\punct}{\sum}\widehat{S}^p,
\end{equation}
where the first sum runs over all interior non-self-folded triangles.% In general, if $\tau$ is a tagged triangulation of $\surf$, we
%define $\unredstau=\sunredtaunot\in\unredrataunot=\runredatau$, where $\tau^\circ$ is the ideal triangulation obtained from $\tau$ using the signature
%$\delta_\tau$ of $\tau$.

Finally, we define $\astau$ to be the (right-equivalence class of the) reduced part of $\unredastau$.
\end{defi}

\begin{obs} Note that, since $\surf$ is not a sphere with less than five punctures, each non-self-folded ideal triangle is adjacent to at most two self-folded triangles.
\end{obs}

To illustrate Definition \ref{QPfortriangulation}, we give some examples.

\begin{ex} If $\surf$ has no punctures (so that the boundary of $\surfnoM$ is non-empty and all marked points lie on the boundary), then for every triangulation $\tau$ of $\surf$ we have $\astau=\unredastau$ and all the terms of $\stau$ are oriented triangles of $\qtau$ arising from interior triangles of $\tau$.
\end{ex}

\begin{ex} If the ideal triangulation $\tau$ of $\surf$ does not have self-folded triangles and is such that each puncture $p$ is incident to at least three arcs of $\tau$, then $\unredastau=(A(\tau),S(\tau))$ and $S(\tau)=\sum_{\triangle}\widehat{S}^\triangle+\sum_{p\in P}x_pa^p_1\ldots a^p_d$ is the sum of the oriented triangles of $Q(\tau)$ arising from interior ideal triangles of $\tau$ and non-zero scalar multiples of cycles around the punctures.
\end{ex}

\begin{ex} Consider the ideal triangulation $\tau$ of the twice-punctured hexagon shown in Figure \ref{uglytrianghexagon}.
% UGLY TRIANGULATION ON TWICE-PUNCTURED HEXAGON
        \begin{figure}[!h]
                \caption{}\label{uglytrianghexagon}
                \centering
                \includegraphics[scale=.3]{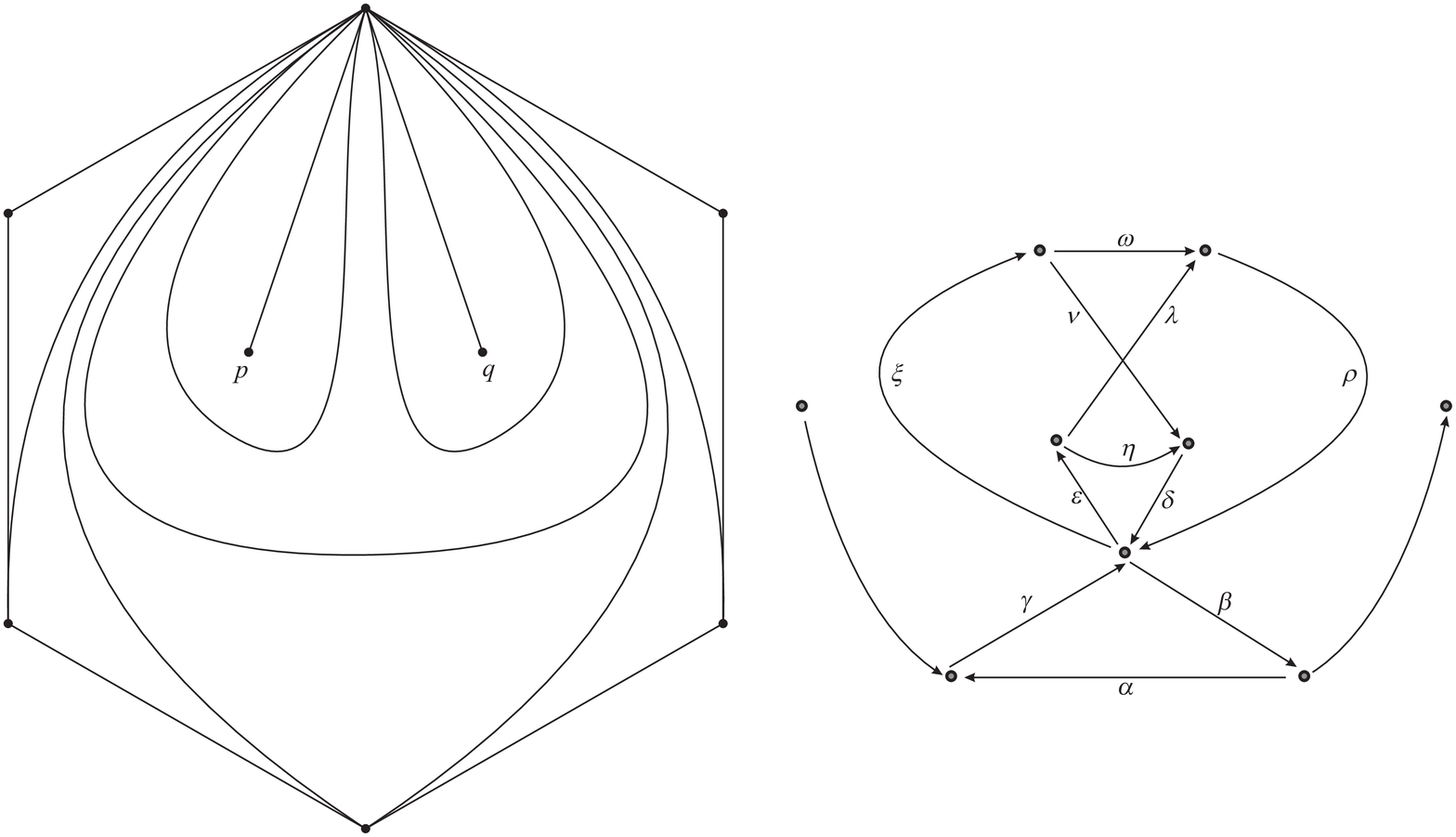}
        \end{figure}
Its (unreduced) signed adjacency quiver is shown on the right. We have $\widehat{S}(\tau)=\alpha\beta\gamma+\delta\eta\varepsilon+
\frac{\rho\omega\xi}{x_px_q}-\frac{\delta\nu\xi}{x_p}-\frac{\varepsilon\rho\lambda}{x_q}$. Since $(\widehat{A}(\tau),\widehat{S}(\tau))$ is clearly 2-acyclic (and hence reduced), we have $\astau=\unredastau$.
\end{ex}

\begin{ex}\label{lessuglyexample} The ideal triangulation $\sigma$ shown in Figure \ref{lessugly} can be obtained from the triangulation $\tau$ of the previous example by a flip.
% LESS UGLY TRIANGULATION ON TWICE-PUNCTURED HEXAGON
        \begin{figure}[!h]
                \caption{}\label{lessugly}
                \centering
                \includegraphics[scale=.3]{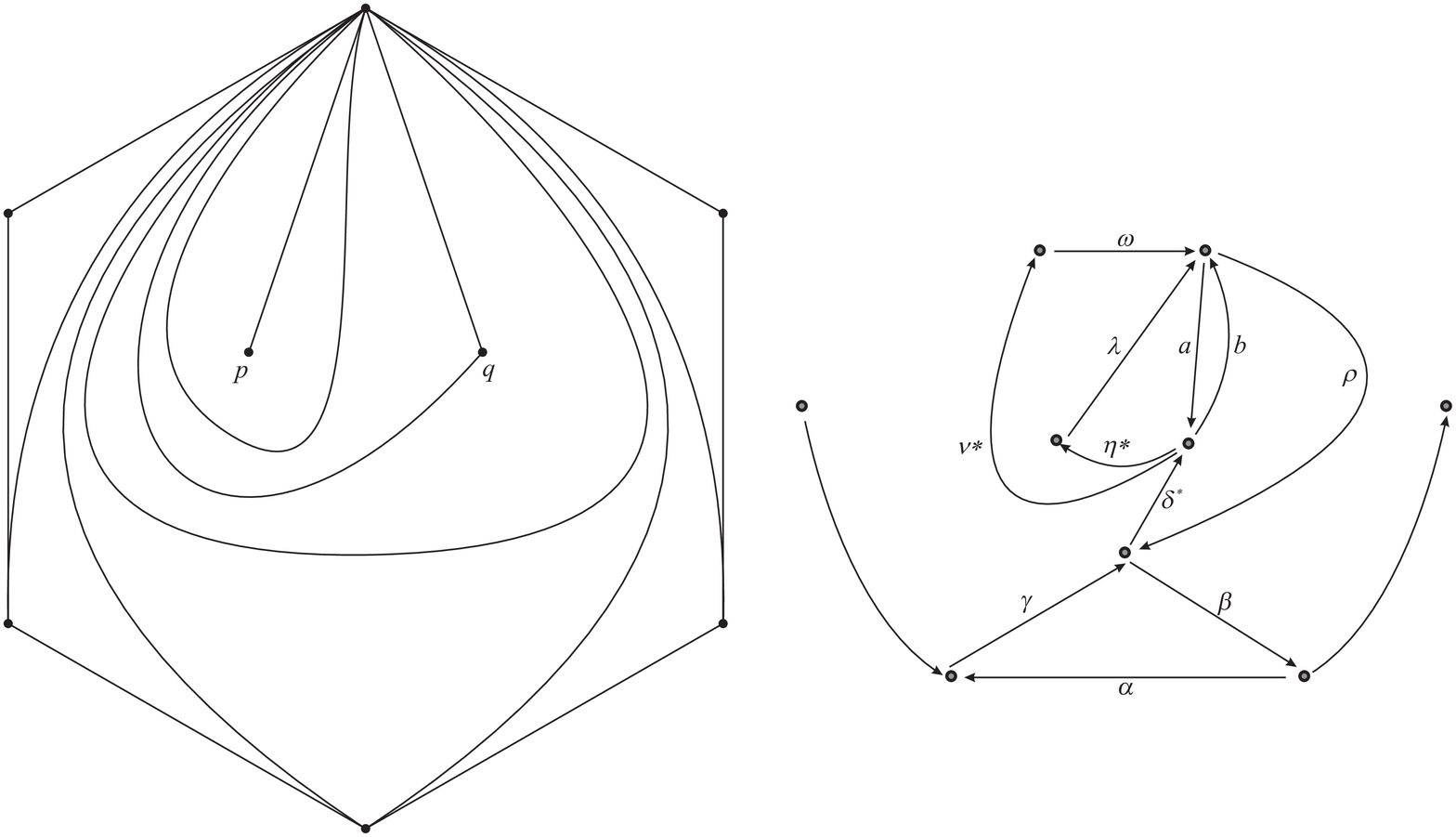}
        \end{figure}
Its unreduced signed adjacency quiver $\widehat{Q(\sigma)}$ is shown on the right, and $\widehat{S}(\sigma)=\alpha\beta\gamma+\delta^*\rho b+a\lambda\eta^*-\frac{a\omega\nu^*}{x_p}+x_qab$. In this example $(\widehat{Q(\sigma)},\widehat{S}(\sigma))$ is not reduced. The $R$-algebra isomorphism $\varphi:R\langle\langle \widehat{A}(\sigma)\rangle\rangle\rightarrow R\langle\langle\widehat{A}(\sigma)\rangle\rangle$ whose action on the arrows is given by $a\mapsto a-\frac{\delta^*\rho}{x_q}$, $b\mapsto b-\frac{\lambda\eta^*}{x_q}+\frac{\omega\nu^*}{x_px_q}$, and the identity on the rest of the arrows, is a right-equivalence between $\unredassigma$ and $(\widehat{A}(\sigma),x_qab+\alpha\beta\gamma-\frac{\delta^*\rho\lambda\eta^*}{x_q}+\frac{\delta^*\rho\omega\nu^*}{x_px_q})$. Therefore, the QP associated to $\sigma$ is, up to right-equivalence, $(Q(\sigma),S(\sigma))$, where $S(\sigma)=\alpha\beta\gamma-\frac{\delta^*\rho\lambda\eta^*}{x_q}-\frac{\delta^*\rho\omega\nu^*}{x_px_q}$.
\end{ex}

\begin{ex}\label{squaredigonexample} Consider the ideal triangulation of the once-punctured square shown in Figure \ref{squarewithdigon}.
% DIGON ON ONCE-PUNCTURED SQUARE
        \begin{figure}[!h]
                \caption{}\label{squarewithdigon}
                \centering
                \includegraphics[scale=.35]{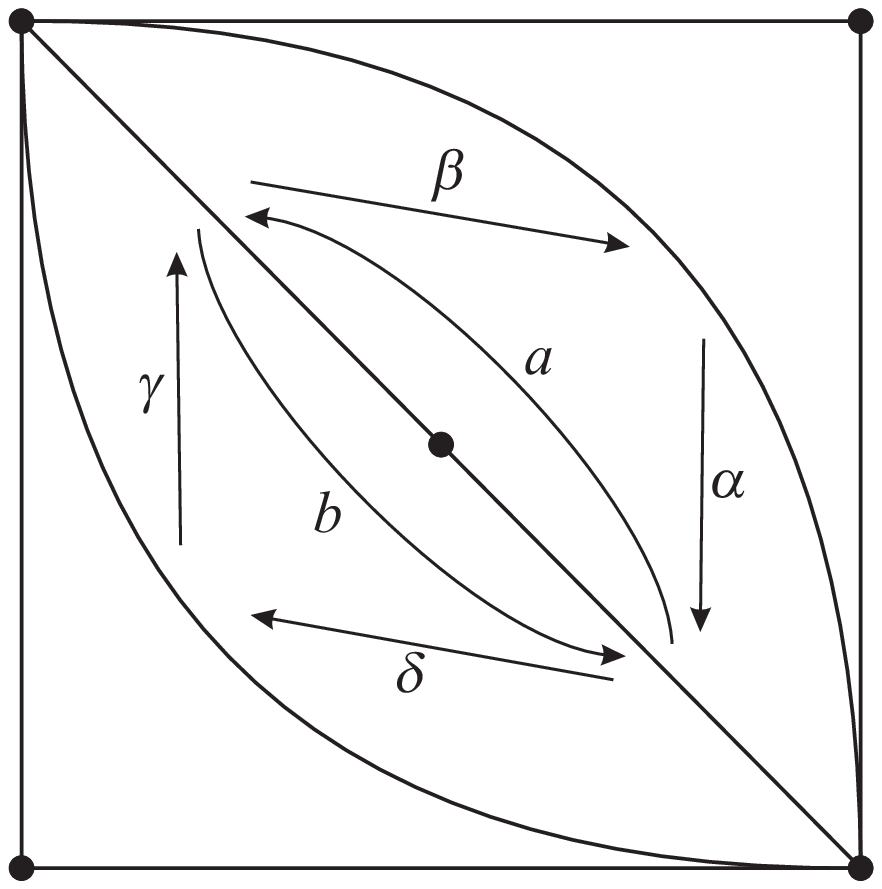}
        \end{figure}
The unreduced signed adjacency quiver $\widehat{Q(\tau)}$ has been drawn on the triangulation. We have $\widehat{S}(\tau)=a\alpha\beta+\gamma\delta b+xab$, where $x\in\field$ is the non-zero scalar attached to the puncture. The $R$-algebra isomorphism $\varphi:\runredatau\rightarrow\runredatau$ whose action on the arrows is given by $a\mapsto a-\frac{\gamma\delta}{x}$, $b\mapsto b-\frac{\alpha\beta}{x}$, and the identity on the rest of the arrows, is a right-equivalence between $\unredastau$ and $(\widehat{A}(\tau),xab-\frac{\gamma\delta\alpha\beta}{x})$. Therefore, $\astau=(A(\tau),-\frac{\gamma\delta\alpha\beta}{x})$.
\end{ex}

The last two examples have something in common, namely, in both of them there is a puncture that is incident to exactly two arcs of the corresponding ideal triangulation. The reduction procedure applied in these examples is general, as we will see in a moment. Now, of course there could be many 2-cycles in $\widehat{Q(\tau)}$, but there are only finitely many of them and each of them has a non-zero scalar multiple appearing as a term of $\widehat{S}(\tau)$, which means that the underlying quiver of the reduced part of $(\widehat{Q(\tau)},\widehat{S}(\tau))$ coincides with the signed adjacency quiver $Q(\tau)$. Moreover, the reduction process leading to $S(\tau)$ can be split into steps, so that one takes care of the 2-cycles one by one.

\begin{figure}[!h]
                \caption{}\label{expressiondigon1}
                \centering
                \includegraphics[scale=.35]{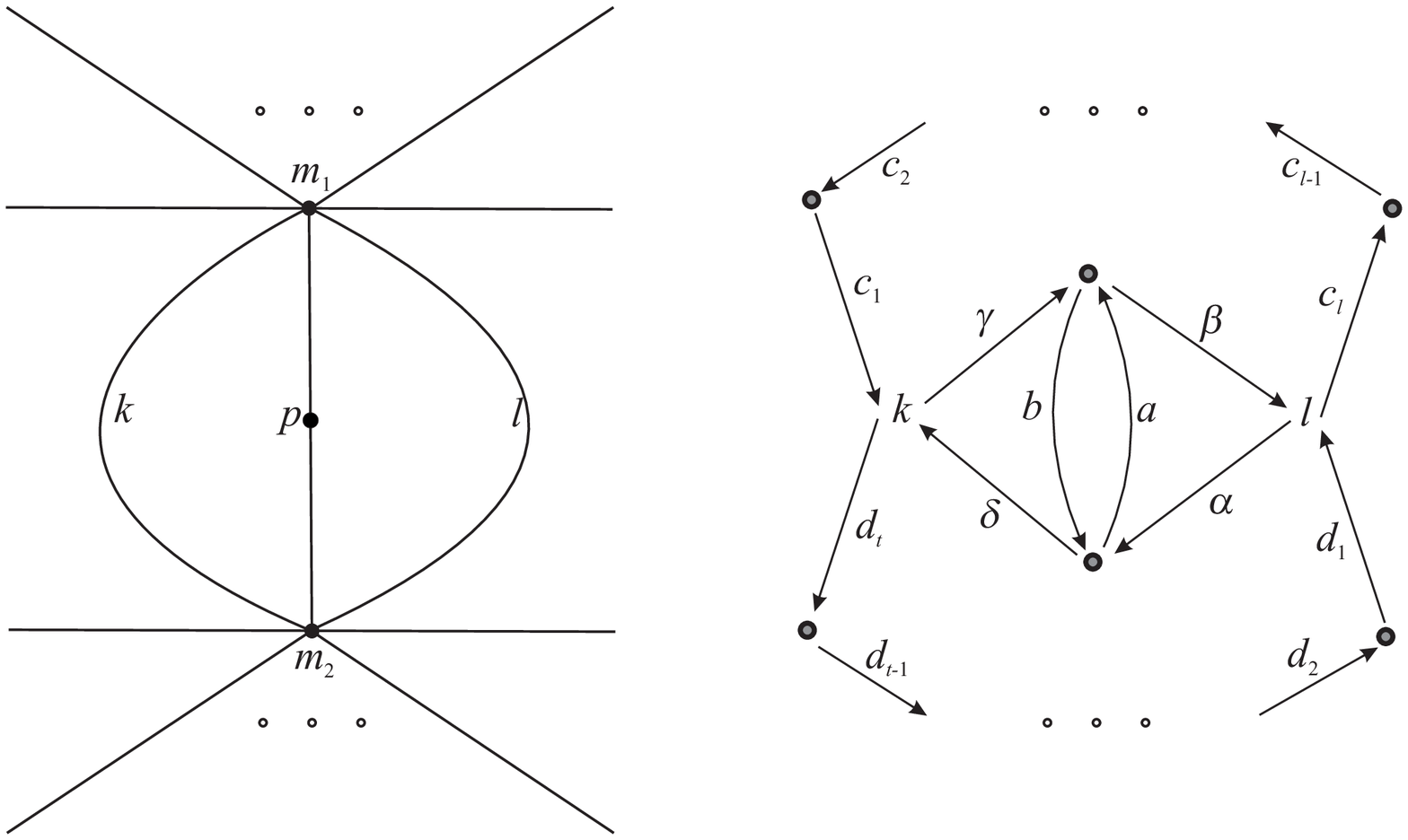}
        \end{figure}

Take an ideal triangulation $\tau$ of $\surf$, and assume that the puncture $p$ is incident to exactly two arcs of $\tau$ as in Figure \ref{expressiondigon1},
% EXPRESSION FOR REDUCED POTENTIAL ON DIGON
assume also that $j$ and $k$ are indeed arcs of $\tau$ (and not part of the boundary of $\surfnoM$). First, consider the case where none of $j$ or $k$ encloses a self-folded triangle. Then the unreduced potential of $\tau$ is $\widehat{S}(\tau)=\gamma\delta b+a\alpha\beta+x_pab+x_{m_1}\beta\gamma c_1\ldots c_l+x_{m_2}\delta\alpha d_1\ldots d_t+\sptau$, where $\sptau\in\runredatau$ is a potential involving none of the arrows $a,b,\alpha,\beta,\gamma,\delta$. Just as in Example \ref{squaredigonexample}, the $R$-algebra isomorphism $\varphi:\runredatau\rightarrow\runredatau$ whose action on the arrows is given by $a\mapsto a-\frac{\gamma\delta}{x_p}$, $b\mapsto b-\frac{\alpha\beta}{x_p}$, and the identity on the rest of the arrows, is a right-equivalence between $\unredastau$ and $(\widehat{A}(\tau),x_pab-\frac{\alpha\beta\gamma\delta}{x_p}+x_{m_1}\beta\gamma c_1\ldots c_l+x_{m_2}\delta\alpha d_1\ldots d_t+\sptau)$. Since $\sptau$ does not involve any of the arrows $a,b,\alpha,\beta,\gamma,\delta$, this implies that the term $-\frac{\alpha\beta\gamma\delta}{x_p}$ will appear, up to right-equivalence, as a term of $\stau$. In other words, $x_pab+a\alpha\beta+\gamma\delta b$ is replaced by $-\frac{\alpha\beta\gamma\delta}{x_p}$ in the reduction process.

\begin{figure}[!h]
                \caption{}\label{expressiondigon2}
                \centering
                \includegraphics[scale=.35]{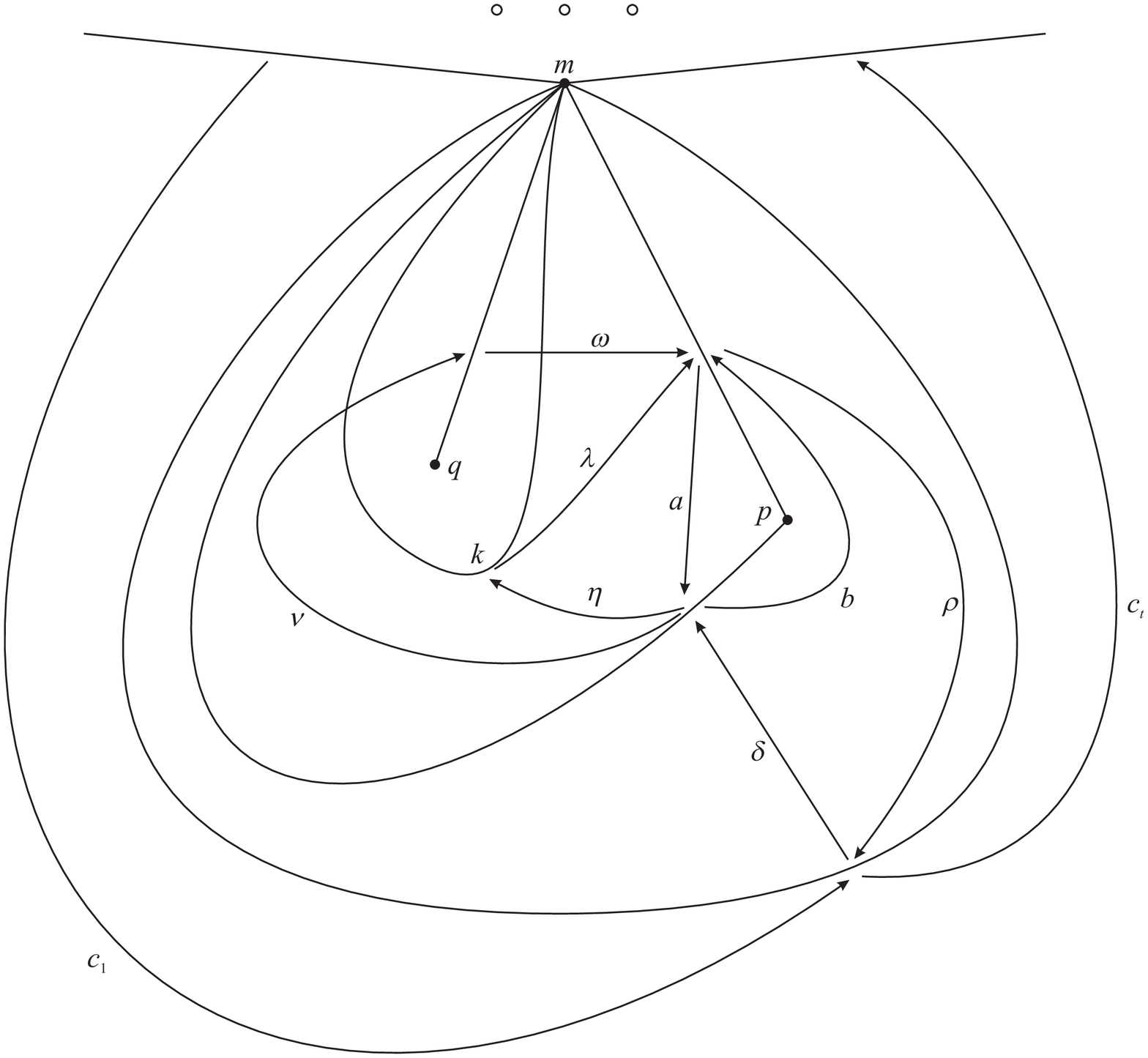}
        \end{figure}

Now, assume that the puncture $p$ is incident to exactly two arcs of $\tau$ and that $k$ encloses a self-folded triangle, see Figure \ref{expressiondigon2}.
% EXPRESSION FOR REDUCED POTENTIAL ON DIGON, ENCLOSING SELF-FOLDED TRIANGLE
Here the unreduced potential is $\widehat{S}(\tau)=a\lambda\eta+\delta\rho b-\frac{a\omega\nu}{x_q}+x_pab+x_m\rho\omega\nu\delta c_1\ldots c_t+\sptau$, where $\sptau\in\runredatau$ is a potential that does not involve any of the arrows $a,b,\delta,\eta,\lambda,\rho,\nu,\omega$. Similarly to Example \ref{lessuglyexample}, the $R$-algebra isomorphism $\varphi:\runredatau\rightarrow\runredatau$ whose action on the arrows is given by $a\mapsto a-\frac{\delta\rho}{x_p}$, $b\mapsto b-\frac{\lambda\eta}{x_p}+\frac{\omega\nu}{x_px_q}$, and the identity on the rest of the arrows, is a right-equivalence between $\unredastau$ and $(\widehat{A}(\sigma), x_pab-\frac{\delta\rho\lambda\eta}{x_p}+\frac{\delta\rho\omega\nu}{x_px_q}+x_m\rho\omega\nu\delta c_1\ldots c_t+\sptau)$. This, together with the fact that $\sptau$ involves none of the arrows $a,b,\delta,\eta,\lambda,\rho,\nu,\omega$, implies that the potential $-\frac{\delta\rho\lambda\eta}{x_p}+\frac{\delta\rho\omega\nu}{x_px_q}$ will appear, up to right-equivalence, as a summand of $\stau$. That is, $x_pab+a\lambda\eta-\frac{a\omega\nu}{x_q}+\delta\rho b$ is replaced by $-\frac{\delta\rho\lambda\eta}{x_p}+\frac{\delta\rho\omega\nu}{x_px_q}$ in the reduction process.

The following lemma says that, using the operation of restriction of QPs, all the QPs we have associated to triangulations can be obtained from QPs associated to triangulations of surfaces without boundary. It will relatively simplify the proof of our first main result.

\begin{lema}\label{allarerestrictions} For every QP of the form $\astau$ there exists an ideal triangulation $\sigma$ of a surface with empty boundary with the following properties:
\begin{itemize}\item $\sigma$ contains all the arcs of $\tau$;
\item the restriction of $\assigma$ to $\tau$ is $\astau$ (except for the fact that the arcs in $\sigma\setminus\tau$ belong to the vertex set of the restriction $(A(\sigma)|_\tau,S(\sigma)|_\tau)$ as isolated vertices, but do not belong to the vertex set of $\astau$; see Remark \ref{remarkisolated}).
\end{itemize}
\end{lema}

\begin{proof} Let $\tau$ be an ideal triangulation of a surface $\surf$ with non-empty boundary. Each boundary component $b$ of $\surfnoM$ is homeomorphic to a circle. Let $m_b$ be the number of marked points lying on $b$. If $m_b=1$ then we can glue $\surfnoM$ and a triangulated twice-punctured monogon along $b$; whereas if $m_b>1$, we can glue $\surfnoM$ and a triangulated $m_b$-gon  along $b$. After doing this for each boundary component of $\surfnoM$, we will end up with an ideal triangulation $\sigma$ of a surface with empty boundary and possessing the desired properties.
\end{proof}

We are ready to state and prove our first main result. Recall from Section \ref{backontriangs} that we write $\sigma=f_i(\tau)$ to indicate that the ideal triangulation $\sigma$ can be obtained from the ideal triangulation $\tau$ by flipping the arc $i$.

\begin{teo}\label{flip<->mutation} Let $\tau$ and $\sigma$ be ideal triangulations of $\surf$. If $\sigma=f_i(\tau)$, then $\mu_i(\atau,\stau)$ and $(\asigma,\ssigma)$ are right-equivalent QPs.
\end{teo}

\begin{proof} By Proposition \ref{resmut=mutres} and Lemma \ref{allarerestrictions}, we can assume, without loss of generality, that the boundary of $\surfnoM$ is empty. We are going to consider several cases, taking on account the configurations that $\tau$ and $\sigma$ can present around the arc $i$ to be flipped. Before proceeding to the case-by-case check, we describe these cases, and to do this, we recall (a slight modification of) the ``puzzle-piece decomposition" mentioned in Remark 4.2 of \cite{FST}.

Any ideal triangulation of $\surf$ can be obtained by means of the following procedure: Consider the four ``puzzle pieces" shown in Figure \ref{puzzlepieces} (a triangle, a once-punctured digon enclosing a self-folded triangle, a twice-punctured monogon enclosing two self-folded triangles, and a once-punctured digon with arcs connecting the puncture to both vertices of the digon).
% PUZZLE PIECE DECOMPOSITION
        \begin{figure}[!h]
                \caption{Puzzle pieces}\label{puzzlepieces}
                \centering
                \includegraphics[scale=.35]{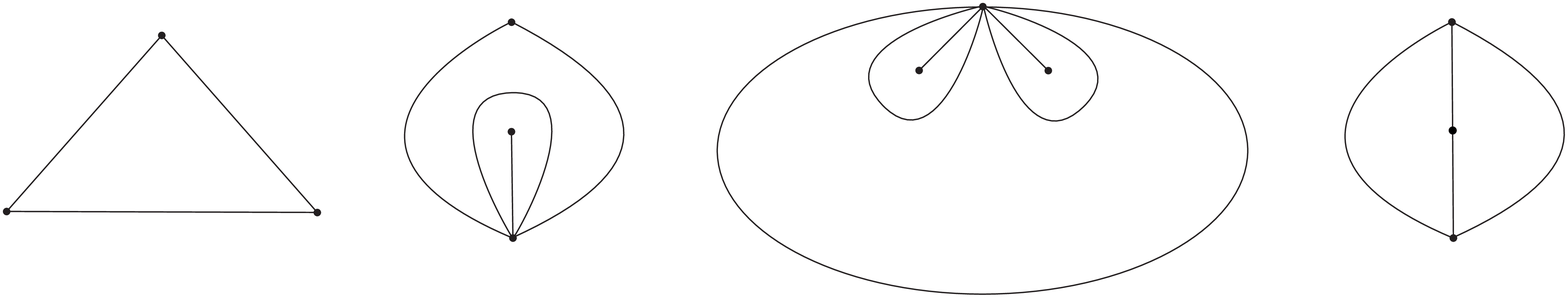}
        \end{figure}
Take several copies of these pieces, assign an orientation to each of the outer sides of these copies and fix a partial matching of these outer sides, never matching two sides of the same copy. (In order to obtain a connected surface, any two puzzle pieces in the collection must be connected via matched pairs). Furthermore, if two sides of a triangle are oriented and matched with two oriented sides of another triangle as shown in Figure \ref{extrapiece}, then replace the pair of triangles with a digon as indicated in Figure \ref{extrapiece} (thus obtaining a new partial matching).
% DIGONS ARE PUZZLE PIECES BY THEIR OWN RIGHT
        \begin{figure}[!h]
                \caption{}\label{extrapiece}
                \centering
                \includegraphics[scale=.4]{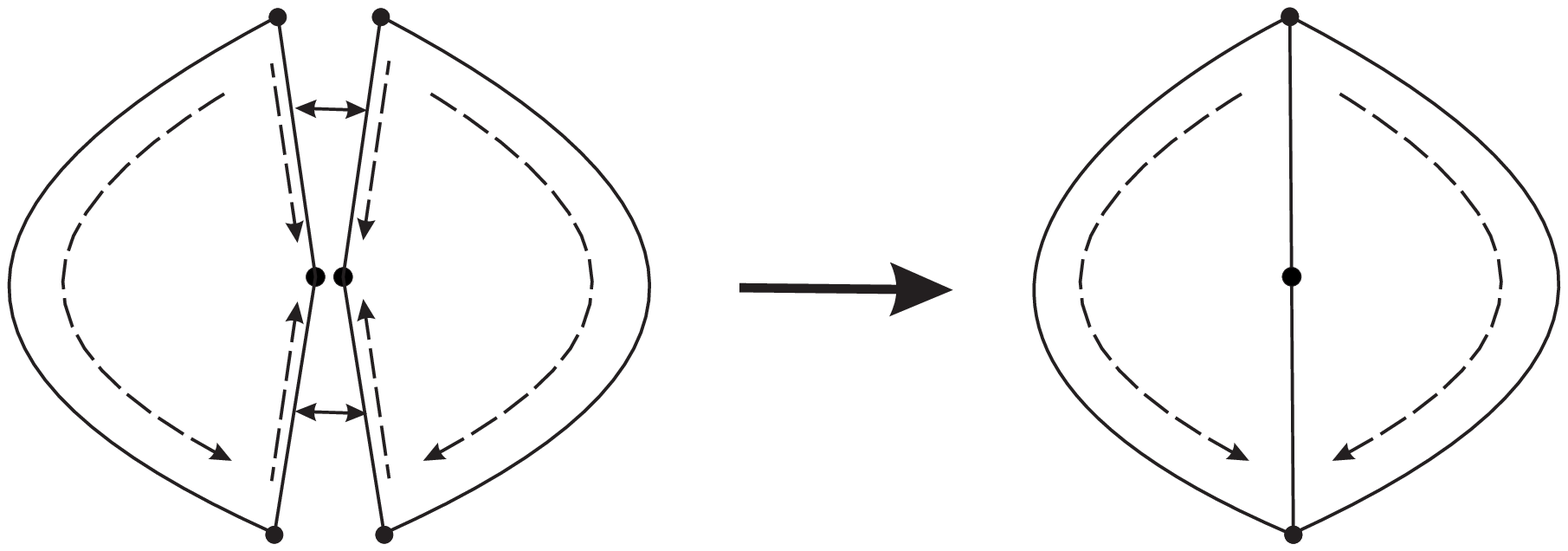}
        \end{figure}
Then glue the puzzle pieces along the matched sides, making sure the orientations match. Though some partial matchings may not lead to an (ideal triangulation of an) oriented surface, any ideal triangulation $\tau$ of an oriented surface can be obtained from a suitable partial matching. Any partial matching giving raise to $\tau$ is called a \emph{puzzle-piece decomposition} of $\tau$.

Now, if we have such a puzzle-piece decomposition of $\tau$, then each flip occurs either inside a puzzle piece, or involves an arc shared by two puzzle pieces. So the proof of the theorem is done by analyzing how these flips affect the corresponding QPs. Figure \ref{matchings} shows a non-oriented list of possible matchings.
Therefore, the proof of the theorem should be carried on by checking the flips that can occur inside a puzzle piece on the one hand, and on the other hand those indicated in Figure \ref{matchings}, where the sides of the puzzle pieces have to be given an orientation and glued along the bold arc, which represents the arc $i$ to be flipped.
% POSSIBLE MATCHINGS
\begin{figure}[!h]
                \caption{}\label{matchings}
                \centering
                \includegraphics[scale=.40]{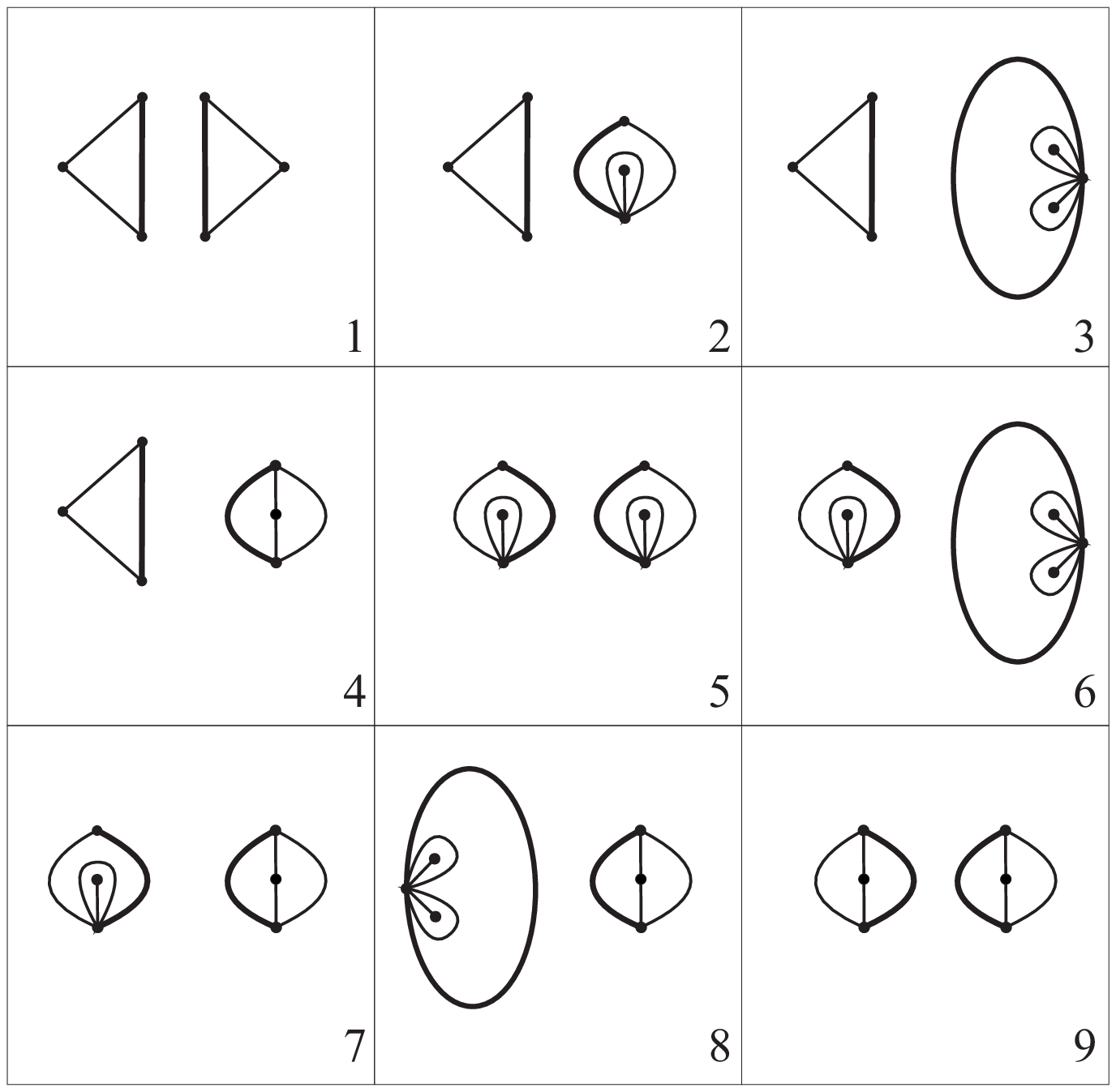}
        \end{figure}

\begin{obs}\label{samemarkedpoint} \begin{enumerate}\item Note that, since QP-mutations are involutive up to right-equivalence, there can be some redundancies when doing this checking: for instance, if this theorem is true in case we flip the loop inside the second puzzle piece of Figure \ref{puzzlepieces}, then it automatically holds also for the flips inside the fourth puzzle piece. Another instance of this redundance is exemplified by the matchings 3 and 5 of Figure \ref{matchings}.
\item Consider the matchings 1, 2, 3, 4, 5, 7 and 9 of Figure \ref{matchings}. Depending on the surface $\surf$ and the ideal triangulation $\tau$ decomposed into puzzle pieces, after gluing the corresponding pair of puzzle pieces one needs to further consider how the vertices on the boundary of the resulting small surface are matched (some of these vertices may represent the same marked point of $\tau$) because different identifications of these vertices may lead to potentials that require slight variations of the right-equivalences we show in the cases analyzed below. The cases we analyze here will have the implicit assumption that different vertices on the boundary of the small surface obtained after gluing the pair of puzzle pieces represent different marked points of $\tau$. We leave to the reader the verification of the cases where different vertices represent the same marked point.
\end{enumerate}
\end{obs}

Cases 1, 3 and 5 below correspond to the matchings 1, 4 and 6 of Figure \ref{matchings}, while Case 3 corresponds to a flip inside the fourth puzzle piece of Figure \ref{puzzlepieces}. Each of the cases below is divided into three stages: in the first stage, we sketch in a Figure the local configuration of $\tau$ around the arc $i$ to be flipped. This configuration will be located at the left of the respective Figure. On the right side of the Figure, the reader will find the corresponding local configuration of the signed adjacency quiver $Q(\tau)$ (this is done also in order to fix some notation for the several arrows we will have to keep track of). In the second stage, we apply to $\qstau$ the QP-mutation with respect to the arc $i$. In the third stage we flip the arc $i$ to obtain $\sigma$ and compare the corresponding QP with the QP obtained in the second stage. The local configuration that $\sigma$ presents around the flip of $i$ is also sketched at the left of the respective Figure, and on the right side of the Figure, the local configuration of the signed adjacency quiver $Q(\sigma)$.

Having said all this, let us proceed to the case-by-case verification. For the rest of the proof, we shall label each puncture $p$ with the scalar $x_p$; we will also make use, without mentioning it, of the reduction process described after Example \ref{squaredigonexample}.

\begin{case}\label{generic1} (First matching of Figure \ref{matchings}, with only one side of each triangle matched to a side of the other one)  Assume that, around the arc $i$, $\tau$ looks like the configuration in Figure \ref{caseone}, with $l,m,n,t>1$ and
$\alpha,\beta,\gamma,\delta,\varepsilon,\eta\notin\{a_1,\ldots,a_l,b_1,\ldots,b_m,c_1,\ldots,c_n,d_1,\ldots,d_t\}$.
% FIRST GENERIC CASE
        \begin{figure}[!h]
                \caption{Case \ref{generic1}, configuration of $\tau$ around $i$}\label{caseone}
                \centering
                \includegraphics[scale=.35]{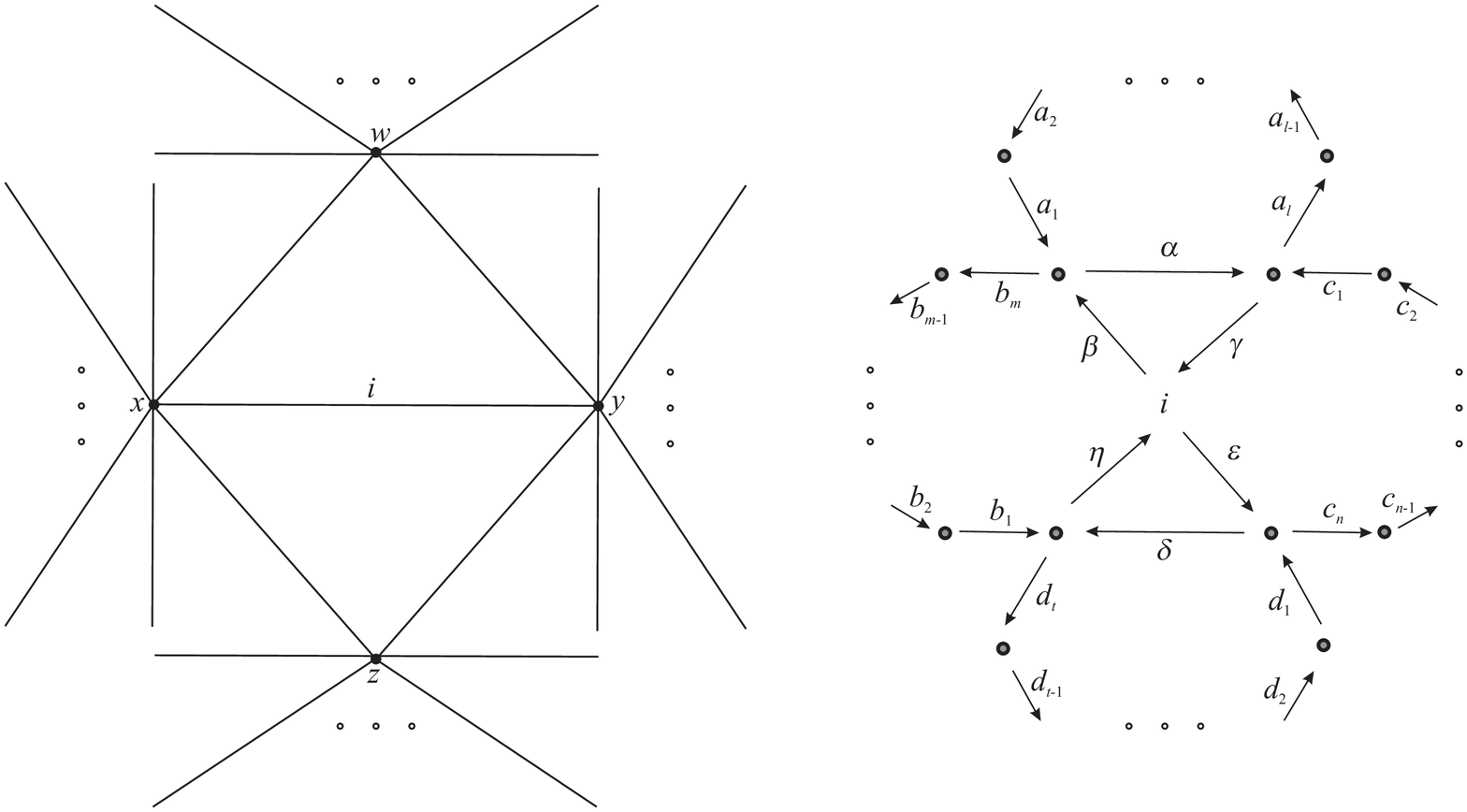}
        \end{figure}
Let us abbreviate $a=a_1\ldots a_l$, $b=b_1\ldots b_m$, $c=c_1\ldots c_n$, $d=d_1\ldots d_t$. Then
$$
\stau=\alpha\beta\gamma+\delta\varepsilon\eta+w\alpha a+x\beta\eta b+y\varepsilon\gamma c+z\delta d+\sptau,$$
with $\sptau\in\ratau$ involving none of the arrows $\alpha,\beta,\gamma,\delta,\varepsilon,\eta$. If we perform the premutation $\premuti$ on $\astau$, we get $\tildeastau$, where $\tildeatau$ is the arrow span of the quiver shown in Figure \ref{premutcaseone}
% PREMUTATED FIRST GENERIC CASE
        \begin{figure}[!h]
                \caption{Case \ref{generic1}, QP-mutation process $\muti\qstau$}\label{premutcaseone}
                \centering
                \includegraphics[scale=.35]{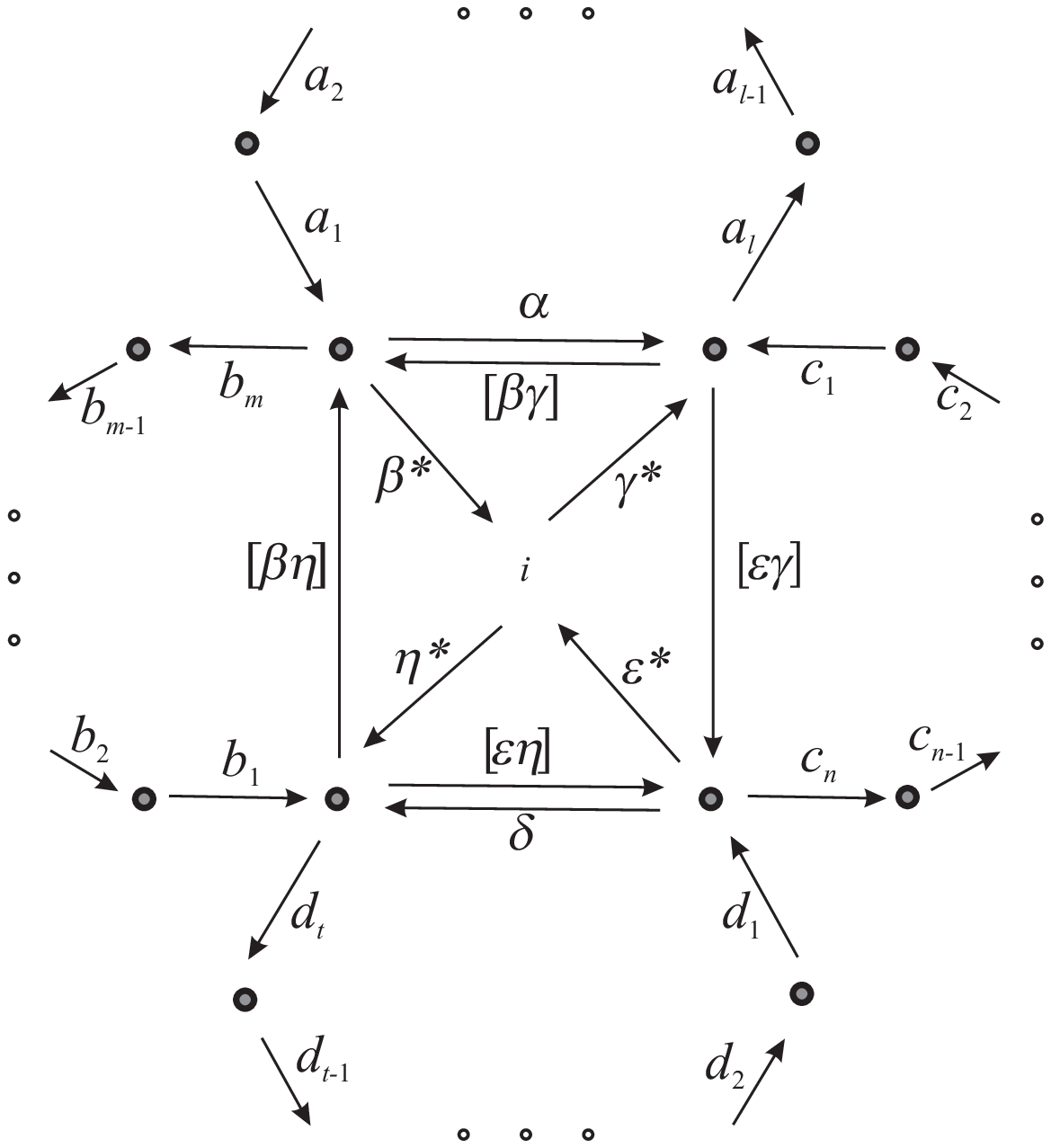}
        \end{figure}
and $\tildestau=\alpha[\beta\gamma]+\delta[\varepsilon\eta]+w\alpha a+x[\beta\eta] b+y[\varepsilon\gamma] c+z\delta
d+\sptau+\gamma^*\beta^*[\beta\gamma]+\eta^*\varepsilon^*[\varepsilon\eta]+\eta^*\beta^*[\beta\eta]+\gamma^*\varepsilon^*[\varepsilon\gamma]\in
\rtildeatau$. The $R$-algebra automorphism $\varphi$ of $\rtildeatau$ whose action on the arrows is given by
$$\alpha\mapsto\alpha-\gamma^*\beta^*,\ [\beta\gamma]\mapsto[\beta\gamma]-wa,\ \delta\mapsto\delta-\eta^*\varepsilon^*,\
[\varepsilon\eta]\mapsto[\varepsilon\eta]-zd,$$
and the identity in the rest of the arrows, sends $\tildestau$ to
$$
\varphi(\tildestau)=\alpha[\beta\gamma]+\delta[\varepsilon\eta]-w\gamma^*\beta^*a+x[\beta\eta] b+y[\varepsilon\gamma] c-z\eta^*\varepsilon^*
d+\sptau+\eta^*\beta^*[\beta\eta]+\gamma^*\varepsilon^*[\varepsilon\gamma].$$
Therefore, the reduced part $\muti\astau$ of $(\tildeatau,\varphi(\tildestau))$ is (up to right-equivalence) the QP on the arrow span $\overatau$ obtained from $\tildeatau$ by deleting the arrows $\alpha,[\beta\gamma],\delta$ and $[\varepsilon\eta]$, with $\varphi(\tildestau)-\alpha[\beta\gamma]-\delta[\varepsilon\eta]$ as its potential.

On the other hand, $\sigma=f_i(\tau)$ and its quiver $\qsigma$ look as Figure \ref{flippedcaseone},
% FLIPPED FIRST GENERIC CASE
        \begin{figure}[!h]
                \caption{Case \ref{generic1}, flip $\sigma=f_i(\tau)$}\label{flippedcaseone}
                \centering
                \includegraphics[scale=.35]{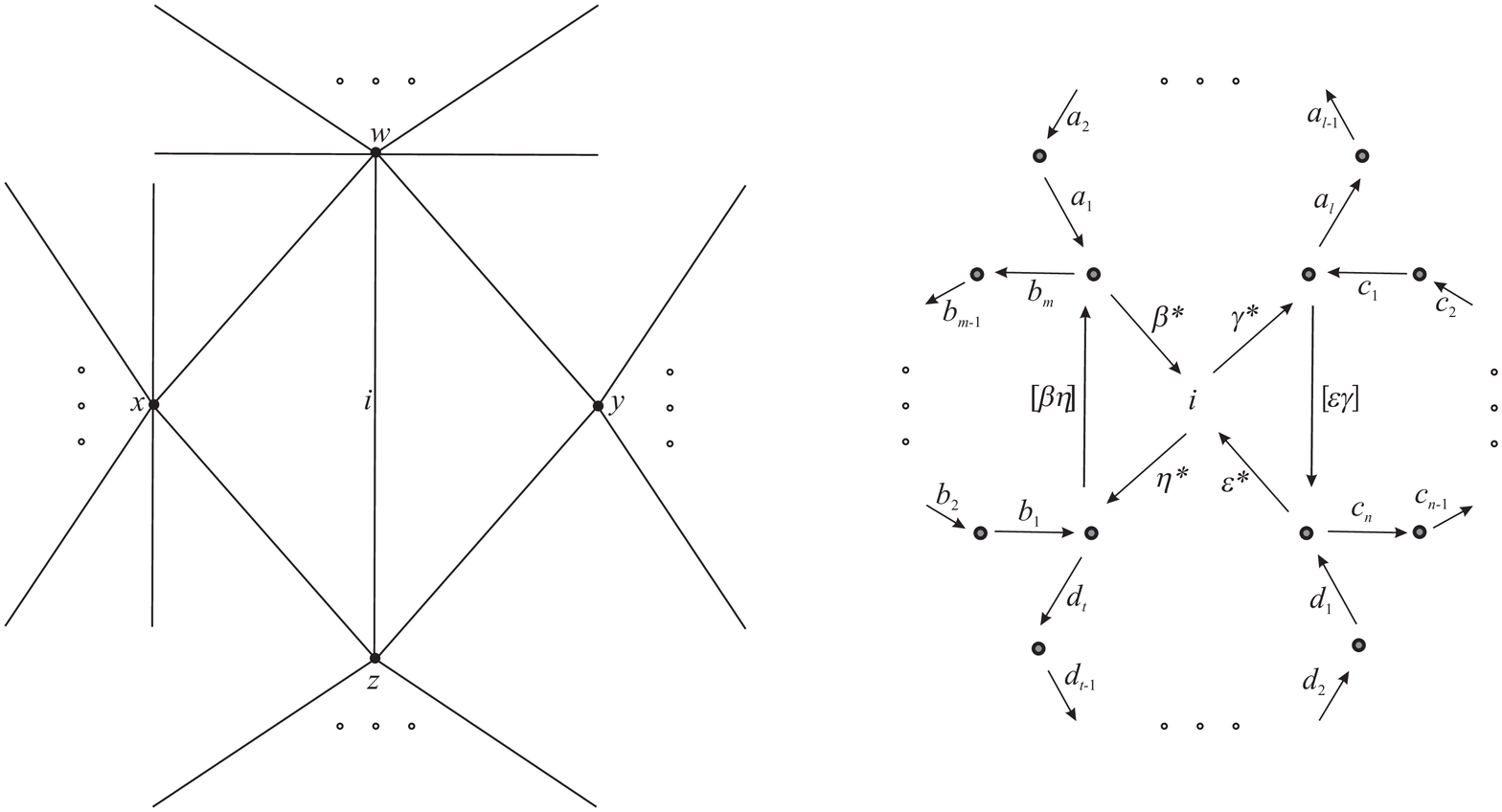}
        \end{figure}
and $\ssigma=w\gamma^*\beta^*a+x[\beta\eta]b+y[\varepsilon\gamma]c+z\eta^*\varepsilon^*d+\eta^*\beta^*[\beta\eta]+
\gamma^*\varepsilon^*[\varepsilon\gamma]+\spsigma$, with $\spsigma=\sptau$. Thus the $R$-algebra isomorphism
$\psi:\roveratau\rightarrow\rasigma$ whose action on the arrows is given by
$$
\beta^*\mapsto-\beta^*,\ \eta^*\mapsto-\eta^*,$$ and the identity
in the rest of the arrows, is a right-equivalence between $\muti\astau$ and $\assigma$.
\end{case}

\begin{case}\label{genericfive}(First matching of Figure \ref{matchings}, with exactly two sides of each triangle matched) Assume that, around the arc $i$, $\tau$ looks like the configuration shown in Figure \ref{casefive}.
% FIFTH GENERIC CASE
        \begin{figure}[!h]
                \caption{Case \ref{genericfive}, configuration of $\tau$ around $i$}\label{casefive}
                \centering
                \includegraphics[scale=.35]{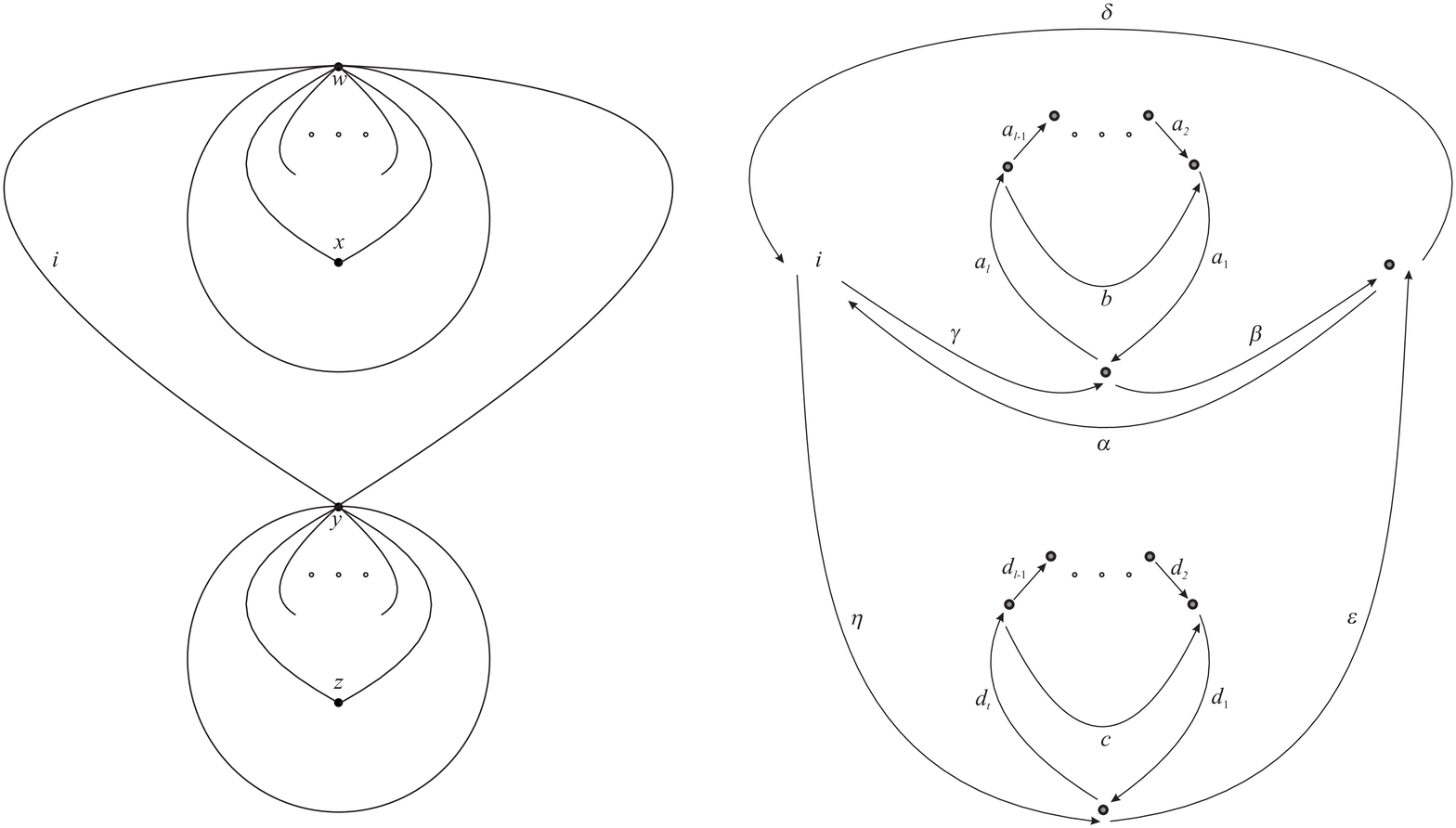}
        \end{figure}
Then
$$
\stau=\beta\gamma\alpha+\epsilon\eta\delta+ba_la_1+cd_td_1+w\beta a\gamma\delta+y\varepsilon d\eta\alpha+\sptau,
$$
with $\sptau\in\ra$ involving none of the arrows $\alpha$, $\beta$, $\gamma$, $\delta$, $\varepsilon$, $\eta$. If we perform the premutation $\premuti$ on $\astau$, we get $\tildeastau$, where $\tildeatau$ is the arrow span of the quiver shown in Figure \ref{premutcasefive},
% PREMUTATED FIFTH GENERIC CASE
        \begin{figure}[!h]
                \caption{Case \ref{genericfive}, QP-mutation process $\muti\astau$}\label{premutcasefive}
                \centering
                \includegraphics[scale=.35]{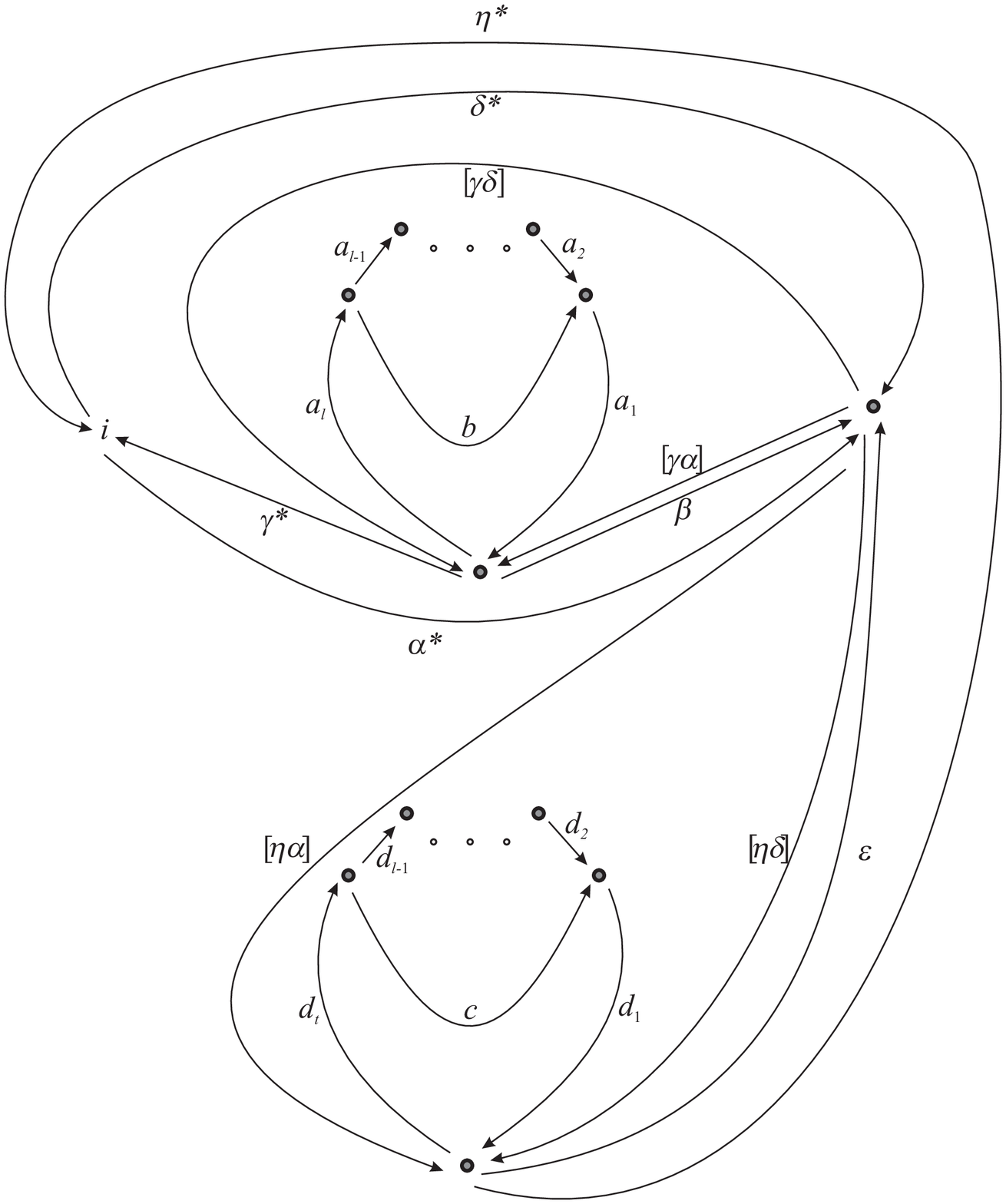}
        \end{figure}
and $\stau=\beta[\gamma\alpha]+\varepsilon[\eta\delta]+ba_la_1+cd_td_1+w\beta a[\gamma\delta]+y\varepsilon d[\eta\alpha]+\sptau+\alpha^*\gamma^*[\gamma\alpha]+\delta^*\eta^*[\eta\delta]+\delta^*\gamma^*[\gamma\delta]+\alpha^*\eta^*[\eta\alpha]
\in\rtildeatau$. The $R$-algebra automorphism $\varphi$ of $\rtildeatau$ whose action on the arrows is given by
$$
\beta\mapsto\beta-\alpha^*\gamma^*,\ [\gamma\alpha]\mapsto[\gamma\alpha]-wa[\gamma\delta],\ \varepsilon\mapsto\varepsilon-\delta^*\eta^*,\ [\eta\delta]\mapsto[\eta\delta]-yd[\eta\alpha],
$$
and the identity in the rest of the arrows, sends $\tildestau$ to
$$
\varphi(\tildestau)=\beta[\gamma\alpha]+\varepsilon[\eta\delta]+ba_la_1+cd_td_1-w\alpha^*\gamma^*a[\gamma\delta]-y\delta^*\eta^*d[\eta\alpha]+
\sptau+\delta^*\gamma^*[\gamma\delta]+\alpha^*\eta^*[\eta\alpha].
$$
Therefore, the reduced part $\muti\astau$ of $\tildeastau$ is (up to right-equivalence) the QP on the arrow span $\overatau$ obtained from $\tildeatau$ by deleting the arrows $\beta$, $[\gamma\alpha]$, $\varepsilon$ and $[\eta\delta]$, with $\varphi(\tildestau)-\beta[\gamma\alpha]+\varepsilon[\eta\delta]$ as its potential.

On the other hand $\sigma=f_i(\tau)$ and its quiver $\qsigma$ look as Figure \ref{flippedcasefive},
% FLIPPED FIFTH GENERIC CASE
        \begin{figure}[!h]
                \caption{Case \ref{genericfive}, flip $\sigma=f_i(\tau)$}\label{flippedcasefive}
                \centering
                \includegraphics[scale=.35]{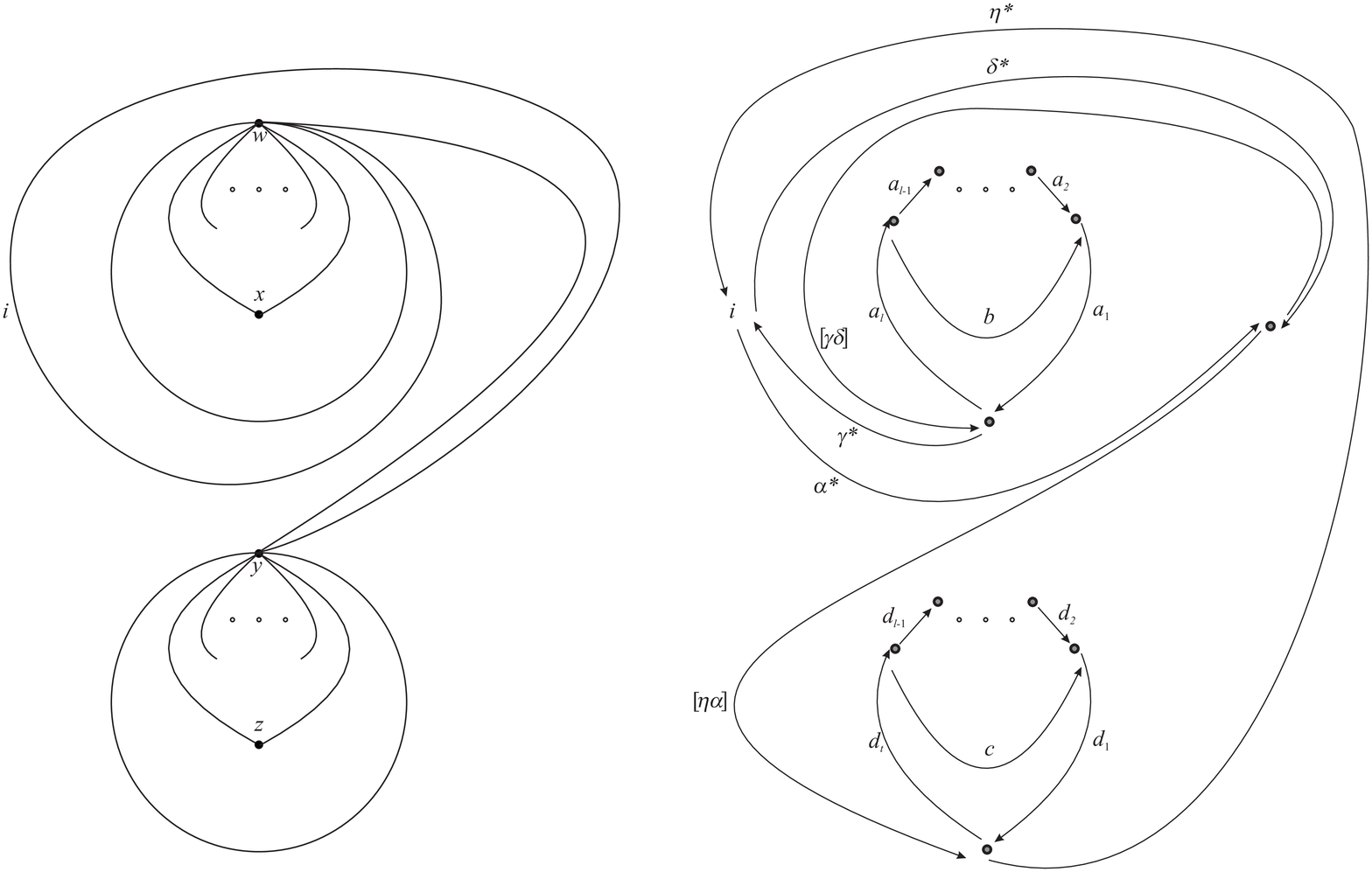}
        \end{figure}
and $S(\sigma)=ba_la_1+cd_td_1+\delta^*\gamma^*[\gamma\delta]+\alpha^*\eta^*[\eta\alpha]+w\alpha^*\gamma^*a[\gamma\delta]+
y\delta^*\eta^*d[\eta\alpha]+\spsigma$, with $\spsigma=\sptau$. Thus the $R$-algebra isomorphism $\psi:\roveratau\rightarrow\rasigma$ whose action on the arrows is given by
$$
\gamma^*\mapsto-\gamma^*,\ \delta^*\mapsto-\delta^*,
$$
and the identity in the rest of the arrows, is a right-equivalence between $\muti\astau$ and $\assigma$.
\end{case}

\begin{case}\label{new2} (Second matching of Figure \ref{matchings}) Assume that, around the arc $i$, $\tau$ looks like the configuration in Figure \ref{newcasetwo}, with the arc $k$ not enclosing a self-folded triangle.
% FIRST GENERIC CASE
        \begin{figure}[!h]
                \caption{Case \ref{new2}, configuration of $\tau$ around $i$}\label{newcasetwo}
                \centering
                \includegraphics[scale=.35]{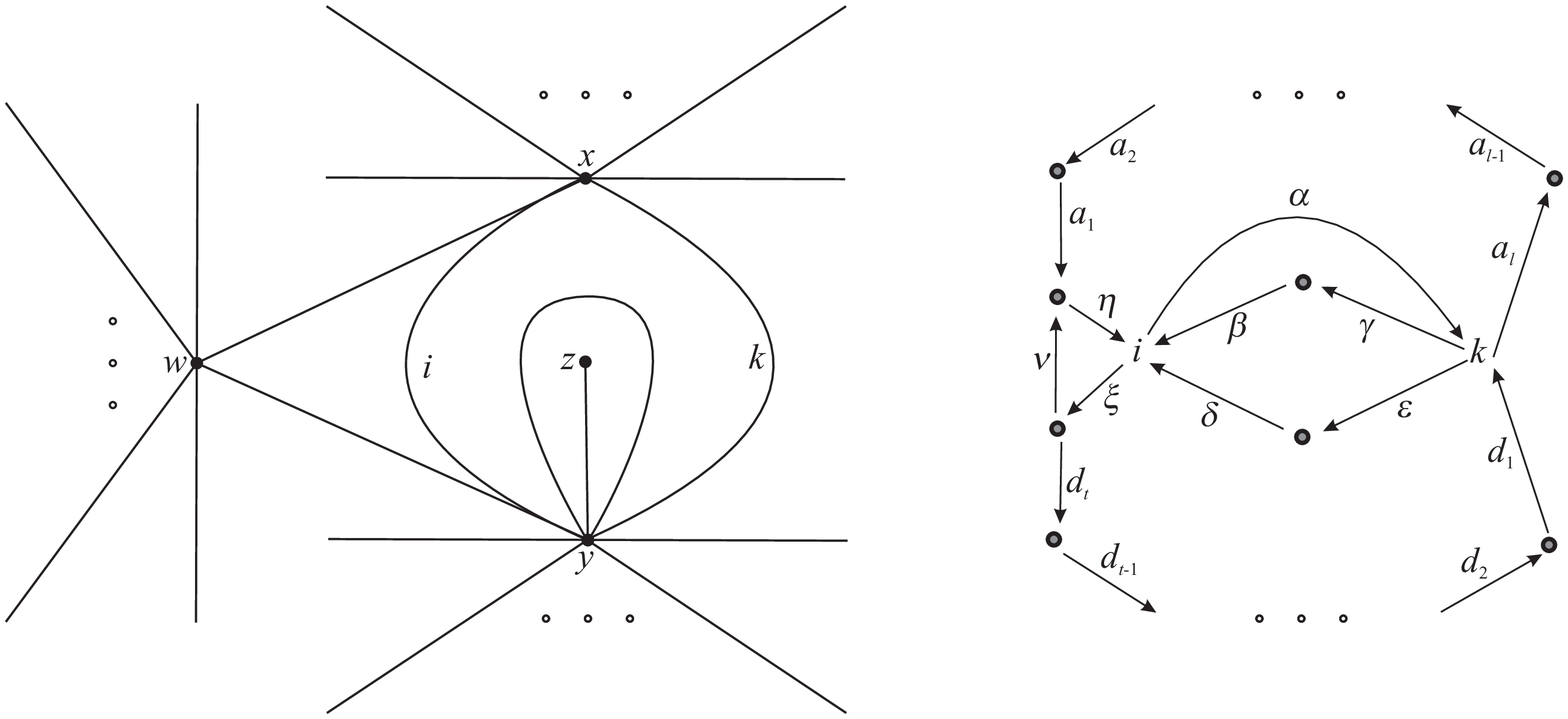}
        \end{figure}
Let us abbreviate $a=a_1\ldots a_l$, $b=b_1\ldots b_m$, $d=d_1\ldots d_t$. Then
$$
\stau=\alpha\beta\gamma+\eta\nu\xi-\frac{\alpha\delta\varepsilon}{z}+x\alpha\eta a+w\nu b+y\xi\delta\varepsilon d\sptau,
$$
with $\sptau\in\ratau$ involving none of the arrows $\alpha,\beta,\gamma,\delta,\varepsilon,\eta,\nu,\xi$. If we perform the premutation $\premuti$ on $\astau$, we get $\tildeastau$, where $\tildeatau$ is the arrow span of the quiver shown in Figure \ref{newpremutcasetwo}
% PREMUTATED FIRST GENERIC CASE
        \begin{figure}[!h]
                \caption{Case \ref{new2}, QP-mutation process $\muti\qstau$}\label{newpremutcasetwo}
                \centering
                \includegraphics[scale=.35]{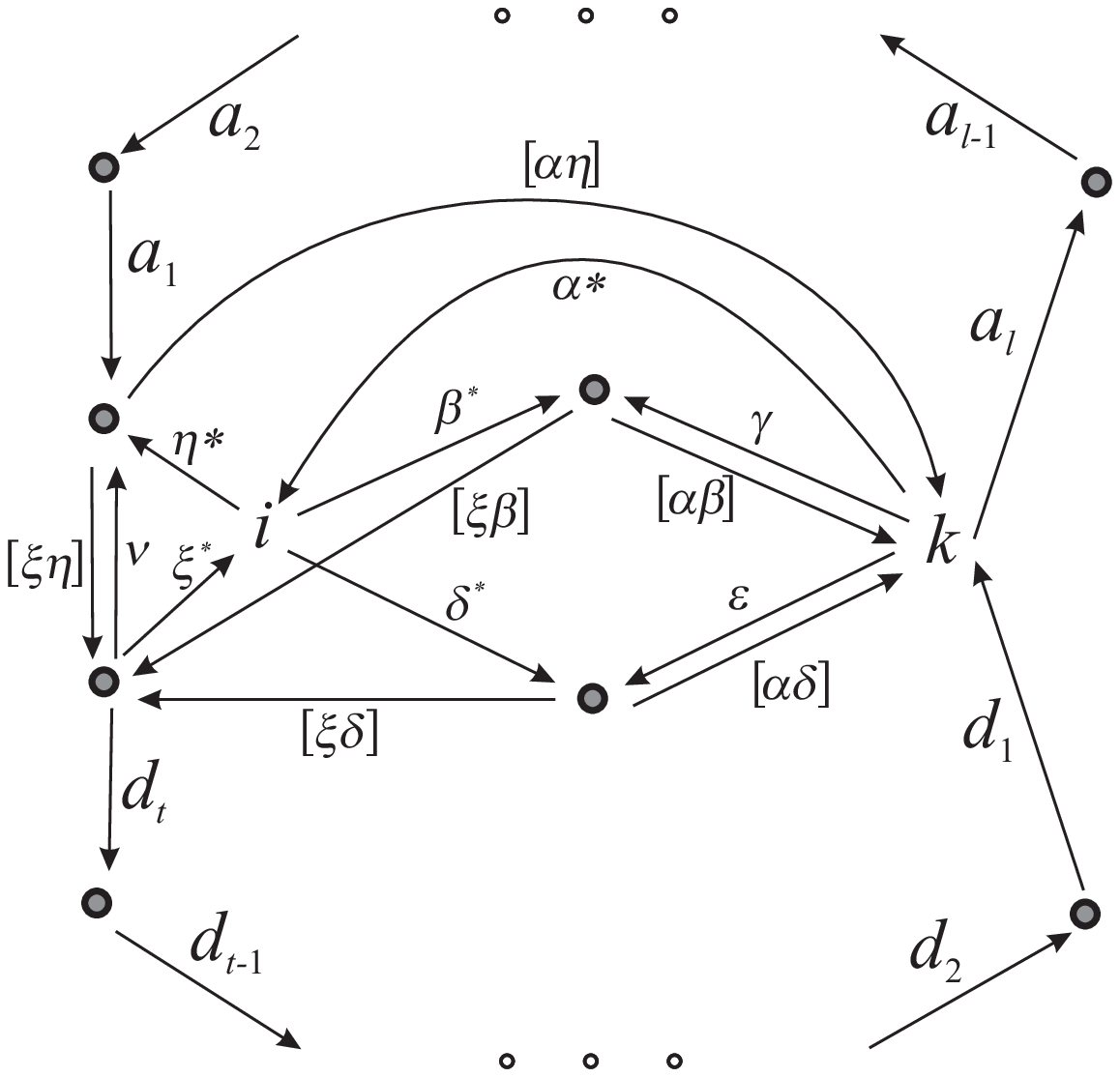}
        \end{figure}
and $\tildestau=[\alpha\beta]\gamma+[\xi\eta]\nu-\frac{\alpha\delta}{z}+x[\alpha\eta]a+wb\nu+yd[\xi\delta]\varepsilon+\sptau+
[\alpha\beta]\beta^*\alpha^*+[\alpha\eta]\eta^*\alpha^*+[\alpha\delta]\delta^*\alpha^*+[\xi\beta]\beta^*\xi^*+[\xi\delta]\delta^*+
[\xi\eta]\eta^*\xi\in\rtildeatau$. The $R$-algebra automorphism $\varphi$ of $\rtildeatau$ whose action on the arrows is given by
$$
\gamma\mapsto\gamma-\beta^*\alpha^*, \ [\xi\eta]\mapsto[\xi\eta]-wb, \ \nu\mapsto\nu-\eta^*\xi^*,
\ [\alpha\delta]\mapsto[\alpha\delta]-yzd[\xi\delta], \ \varepsilon\mapsto\varepsilon-z\delta^*\alpha^*,
$$
and the identity in the rest of the arrows, sends $\tildestau$ to
$$
\varphi(\tildestau)=[\alpha\beta]\gamma+[\xi\eta]\nu-\frac{[\alpha\delta]\varepsilon}{z}+x[\alpha\eta]a+\sptau+[\alpha\eta]\eta^*\alpha^*
+yzd[\xi\delta]\delta^*\alpha^*+[\xi\beta]\beta^*\xi^*+[\xi\delta]\delta^*\xi^*-wb\eta^*\xi^*.
$$
Therefore, the reduced part $\muti\astau$ of $(\tildeatau,\varphi(\tildestau))$ is (up to right-equivalence) the QP on the arrow span $\overatau$ obtained from $\tildeatau$ by deleting the arrows $\gamma,\varepsilon,\nu,[\alpha\beta],[\alpha\delta]$ and $[\xi\eta]$, with $\varphi(\tildestau)-[\alpha\beta]\gamma-[\xi\eta]\nu+\frac{[\alpha\delta]\varepsilon}{z}$ as its potential.

On the other hand, $\sigma=f_i(\tau)$ and its quiver $\qsigma$ look as Figure \ref{newflippedcasetwo},
% FLIPPED FIRST GENERIC CASE
        \begin{figure}[!h]
                \caption{Case \ref{new2}, flip $\sigma=f_i(\tau)$}\label{newflippedcasetwo}
                \centering
                \includegraphics[scale=.35]{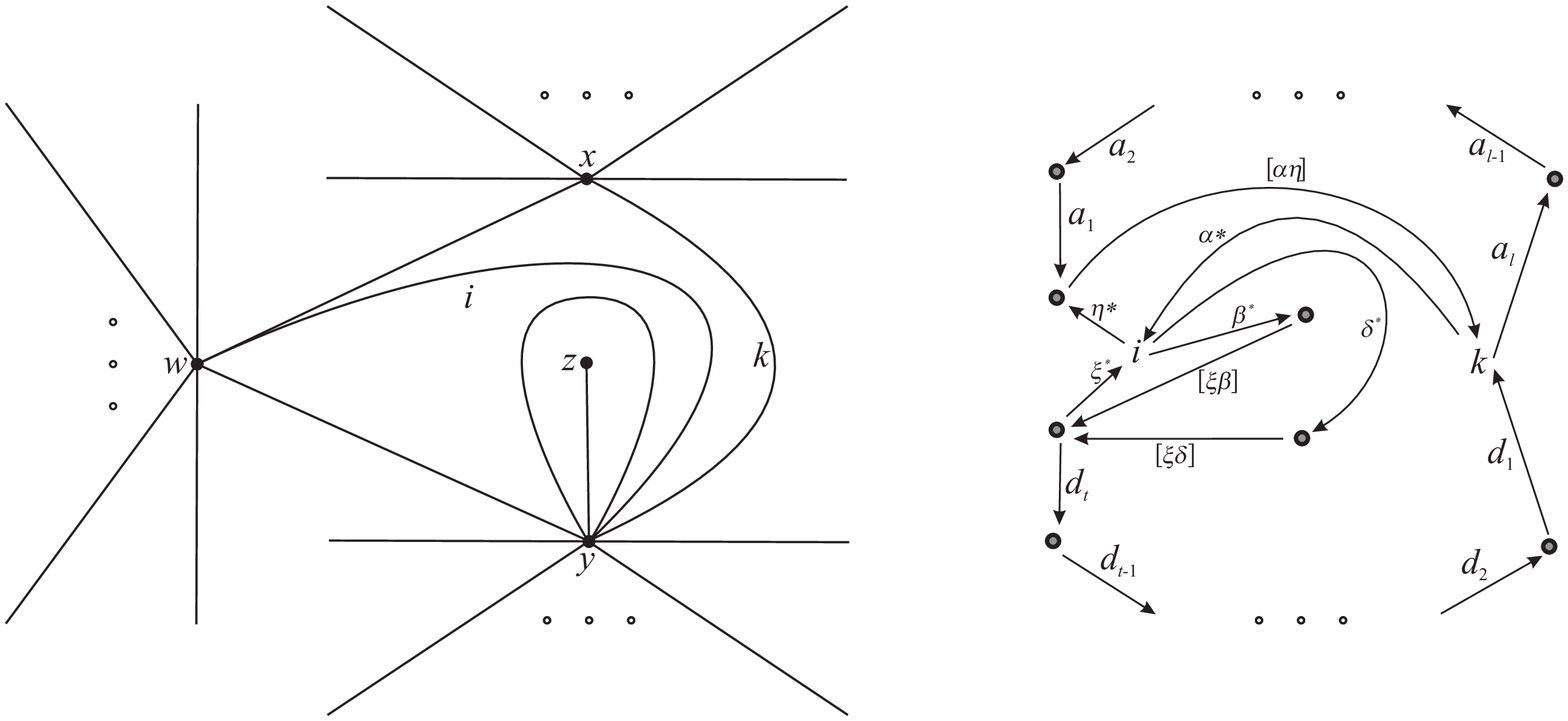}
        \end{figure}
and $\ssigma=[\alpha\eta]\eta^*\alpha^*+[\xi\beta]\beta^*\xi^*-\frac{[\xi\delta]\delta^*\xi^*}{z}+x[\alpha\eta]a+w\eta^*\xi^*b+
y[\xi\delta]\delta^*\alpha^*d+\spsigma$, with $\spsigma=\sptau$. Thus the $R$-algebra isomorphism
$\psi:\roveratau\rightarrow\rasigma$ whose action on the arrows is given by
$$
\alpha^*\mapsto-\alpha^*, \ \delta^*\mapsto-\delta^*, \ \eta^*\mapsto-\eta^*, \ [\xi\delta]\mapsto\frac{[\xi\delta]}{z},
$$
and the identity
in the rest of the arrows, is a right-equivalence between $\muti\astau$ and $\assigma$.
\end{case}

\begin{case}\label{new3} (Third and fifth matchings of Figure \ref{matchings}) Assume that, around the arc $i$, $\tau$ looks like the configuration in Figure \ref{newcasethree}, with $l>1$ and none of the arcs $j$ and $k$ enclosing a self-folded triangle.
% FIRST GENERIC CASE
        \begin{figure}[!h]
                \caption{Case \ref{new3}, configuration of $\tau$ around $i$}\label{newcasethree}
                \centering
                \includegraphics[scale=.3]{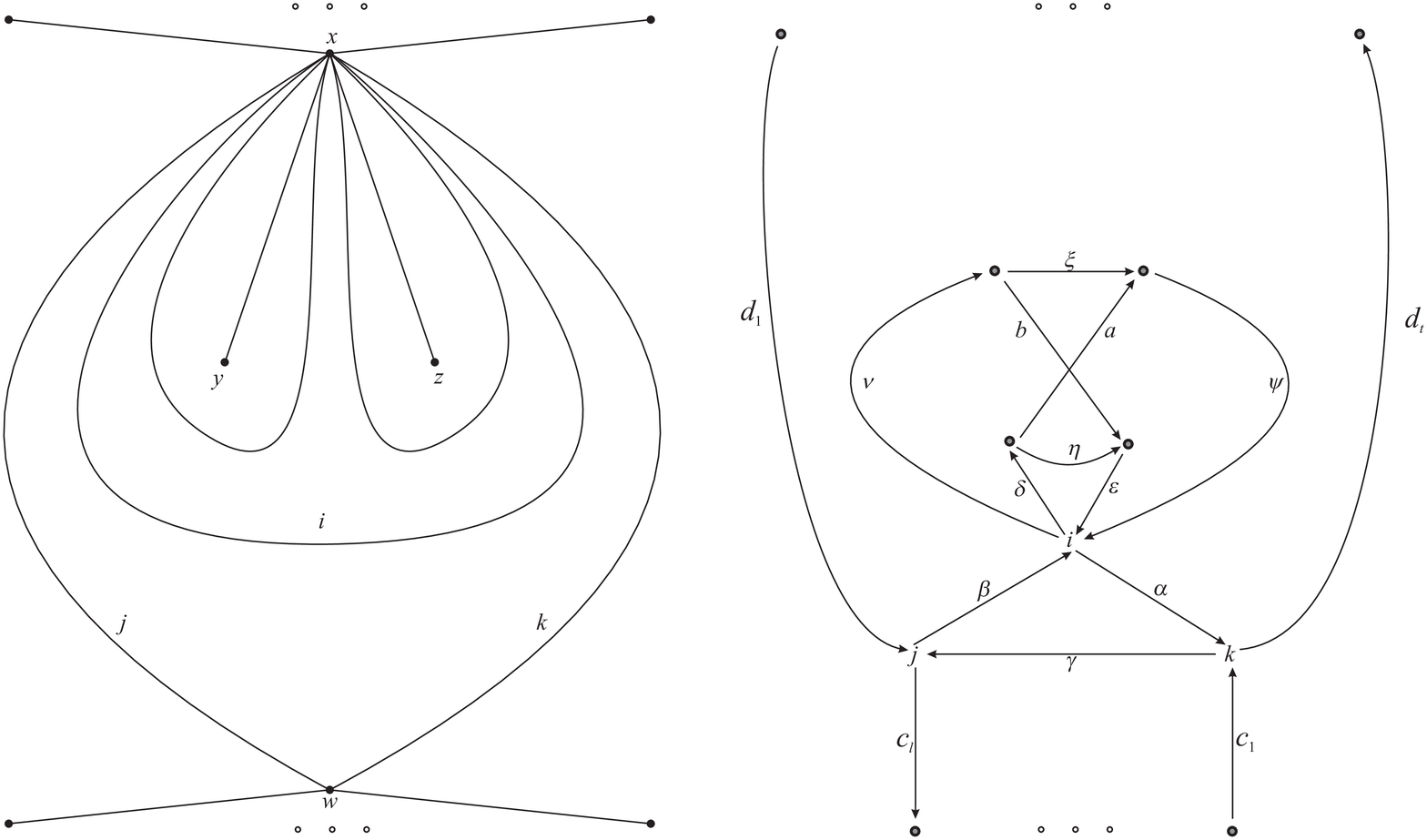}
        \end{figure}
Let us abbreviate $c=c_1\ldots c_l$, $d=d_1\ldots d_t$. Then
$$
\stau=\alpha\beta\gamma+\delta\eta\varepsilon+\frac{\nu\psi\xi}{yz}-\frac{\nu\varepsilon b}{y}-\frac{\delta\psi a}{z}
+x\nu\beta d\alpha\psi\xi+wc\gamma+\sptau,
$$
with $\sptau\in\ratau$ involving none of the arrows $\alpha,\beta,\gamma,\delta,\varepsilon,\eta,\nu,\psi,\xi,a,b$. If we perform the premutation $\premuti$ on $\astau$, we get $\tildeastau$, where $\tildeatau$ is the arrow span of the quiver shown in Figure \ref{newpremutcasethree}
% PREMUTATED FIRST GENERIC CASE
        \begin{figure}[!h]
                \caption{Case \ref{new3}, QP-mutation process $\muti\qstau$}\label{newpremutcasethree}
                \centering
                \includegraphics[scale=.3]{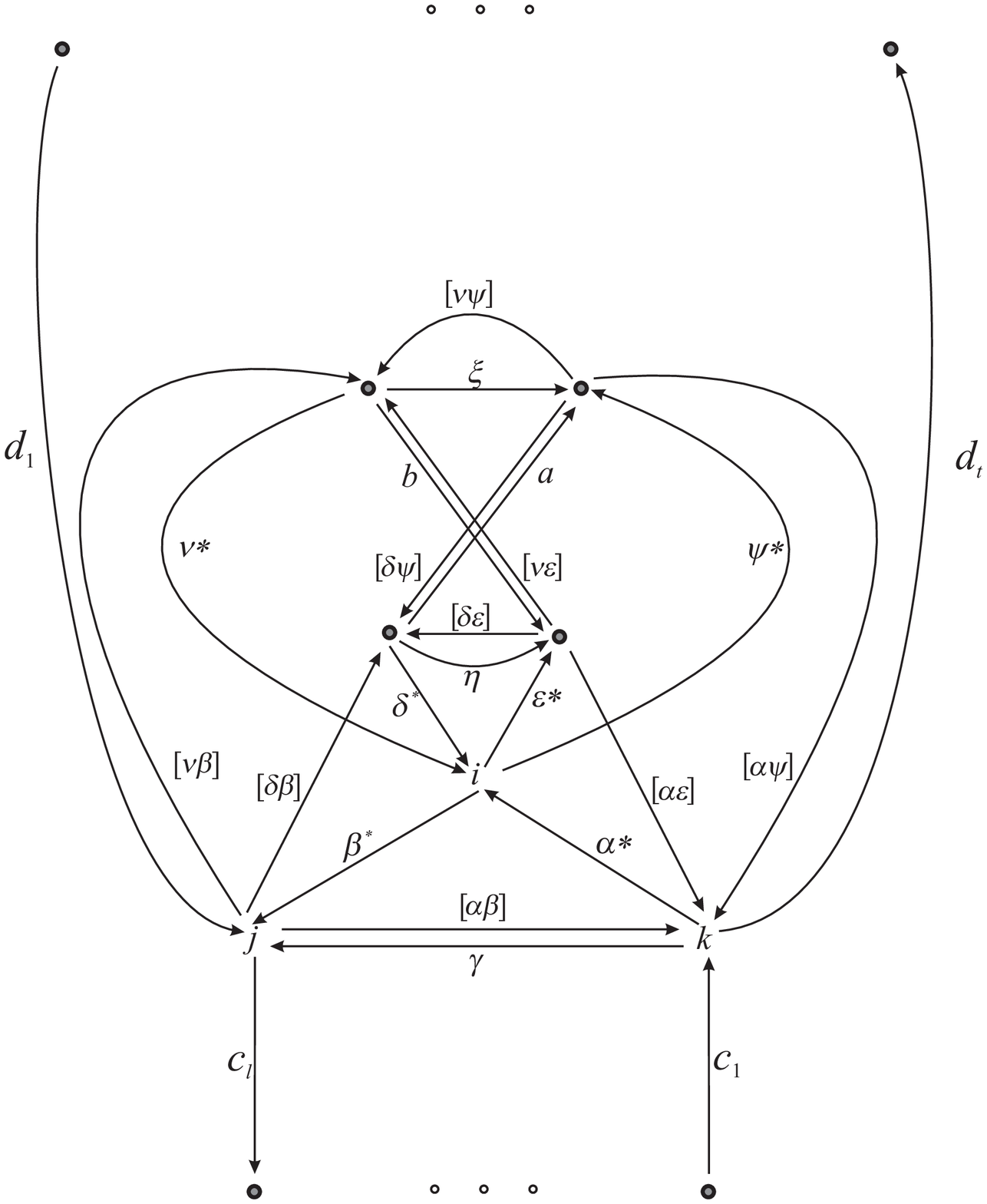}
        \end{figure}
and $\tildestau=[\alpha\beta]\gamma+[\delta\varepsilon]\eta+\frac{[\nu\psi]\xi}{yz}-\frac{[\nu\varepsilon]b}{y}-\frac{[\delta\psi]a}{z}+
x[\nu\beta]d[\alpha\psi]\xi+wc\gamma+\sptau+[\alpha\beta]\beta^*\alpha^*+[\alpha\varepsilon]\varepsilon^*\alpha^*+[\alpha\psi]\psi^*\alpha^*
+[\delta\beta]\beta^*\delta^*+[\delta\varepsilon]\varepsilon^*\delta^*+[\delta\psi]\psi^*\delta^*
+[\nu\beta]\beta^*\nu^*+[\nu\varepsilon]\varepsilon^*\nu^*+[\nu\psi]\psi^*\nu^*\in\rtildeatau$. The $R$-algebra automorphism $\varphi$ of $\rtildeatau$ whose action on the arrows is given by
$$
[\alpha\beta]\mapsto[\alpha\beta]-wc, \ \gamma\mapsto\gamma-\beta^*\alpha^*, \ \eta\mapsto\eta-\varepsilon^*\delta^*, \
[\nu\psi]\mapsto[\nu\psi]-xyz[\nu\beta]d[\alpha\psi],
$$
$$
\xi\mapsto\xi-yz\psi^*\nu^*, \ b\mapsto b+y\varepsilon^*\nu^*, \
a\mapsto a+z\psi^*\delta^*,
$$
and the identity in the rest of the arrows, sends $\tildestau$ to
$$
\varphi(\tildestau)=[\alpha\beta]\gamma+[\delta\varepsilon]\eta+\frac{[\nu\psi]\xi}{yz}-\frac{[\nu\varepsilon]b}{y}-\frac{[\delta\psi]a}{z}+\sptau
-wc\beta^*\alpha^*+[\alpha\varepsilon]\varepsilon^*\alpha^*+[\alpha\psi]\psi^*\alpha^*+[\delta\beta]\beta^*\delta^*+
$$
$$
+[\nu\beta]\beta^*\nu^*
-xyz[\nu\beta]d[\alpha\psi]\psi^*\nu^*.
$$
Therefore, the reduced part $\muti\astau$ of $(\tildeatau,\varphi(\tildestau))$ is (up to right-equivalence) the QP on the arrow span $\overatau$ obtained from $\tildeatau$ by deleting the arrows $[\alpha\beta],\gamma,[\delta\varepsilon],\eta,[\nu\psi],\xi,[\nu\varepsilon],b,
[\delta\psi]$ and $a$, with $\varphi(\tildestau)-[\alpha\beta]\gamma-[\delta\varepsilon]\eta-\frac{[\nu\psi]\xi}{yz}+\frac{[\nu\varepsilon]}{y}+
\frac{[\delta\psi]a}{z}$ as its potential.

On the other hand, $\sigma=f_i(\tau)$ and its quiver $\qsigma$ look as Figure \ref{newflippedcasethree},
% FLIPPED FIRST GENERIC CASE
        \begin{figure}[!h]
                \caption{Case \ref{new3}, flip $\sigma=f_i(\tau)$}\label{newflippedcasethree}
                \centering
                \includegraphics[scale=.3]{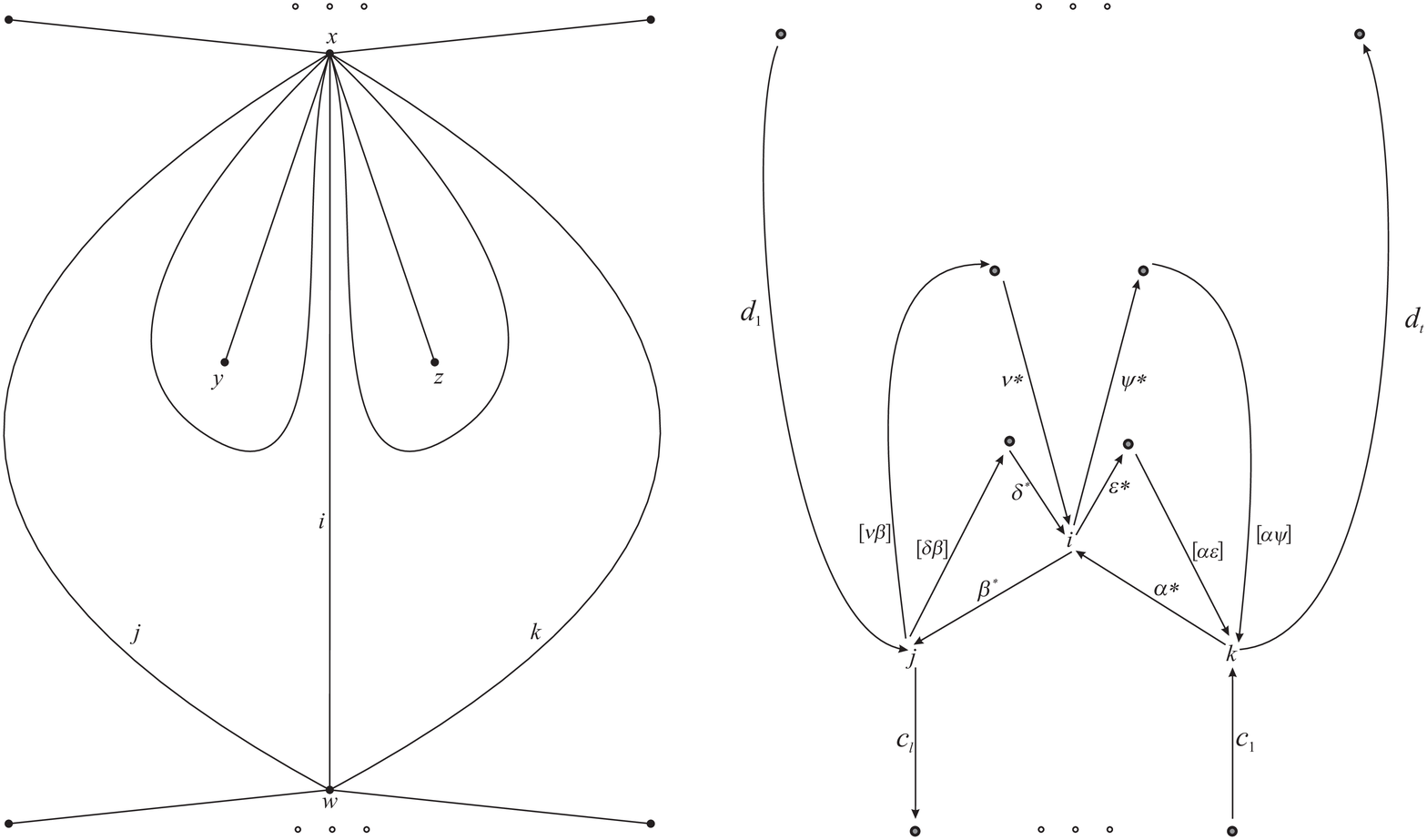}
        \end{figure}
and $\ssigma=[\alpha\varepsilon]\varepsilon^*\alpha^*+[\delta\beta]\beta^*\delta^*-\frac{[\nu\beta]\beta^*\nu^*}{y}
-\frac{[\alpha\psi]\psi^*\alpha^*}{z}+x[\nu\beta]d[\alpha\psi]\psi^*\nu^*+w\beta^*\alpha^*c\spsigma$, with $\spsigma=\sptau$. Thus the $R$-algebra isomorphism
$\psi:\roveratau\rightarrow\rasigma$ whose action on the arrows is given by
$$
\beta^*\mapsto-\beta^* \ [\alpha\psi]\mapsto-\frac{[\alpha\psi]}{z}, \ [\nu\beta]\mapsto\frac{\nu\beta}{y},
$$
and the identity
in the rest of the arrows, is a right-equivalence between $\muti\astau$ and $\assigma$.
\end{case}

\begin{case}\label{generic2} (Fourth matching of Figure \ref{matchings}) Assume that, around the arc $i$, $\tau$ looks like the configuration shown in Figure \ref{casetwo}, with $m,n>0$, $t>1$.
% SECOND GENERIC CASE
        \begin{figure}[!h]
                \caption{Case \ref{generic2}, configuration of $\tau$ around $i$}\label{casetwo}
                \centering
                \includegraphics[scale=.35]{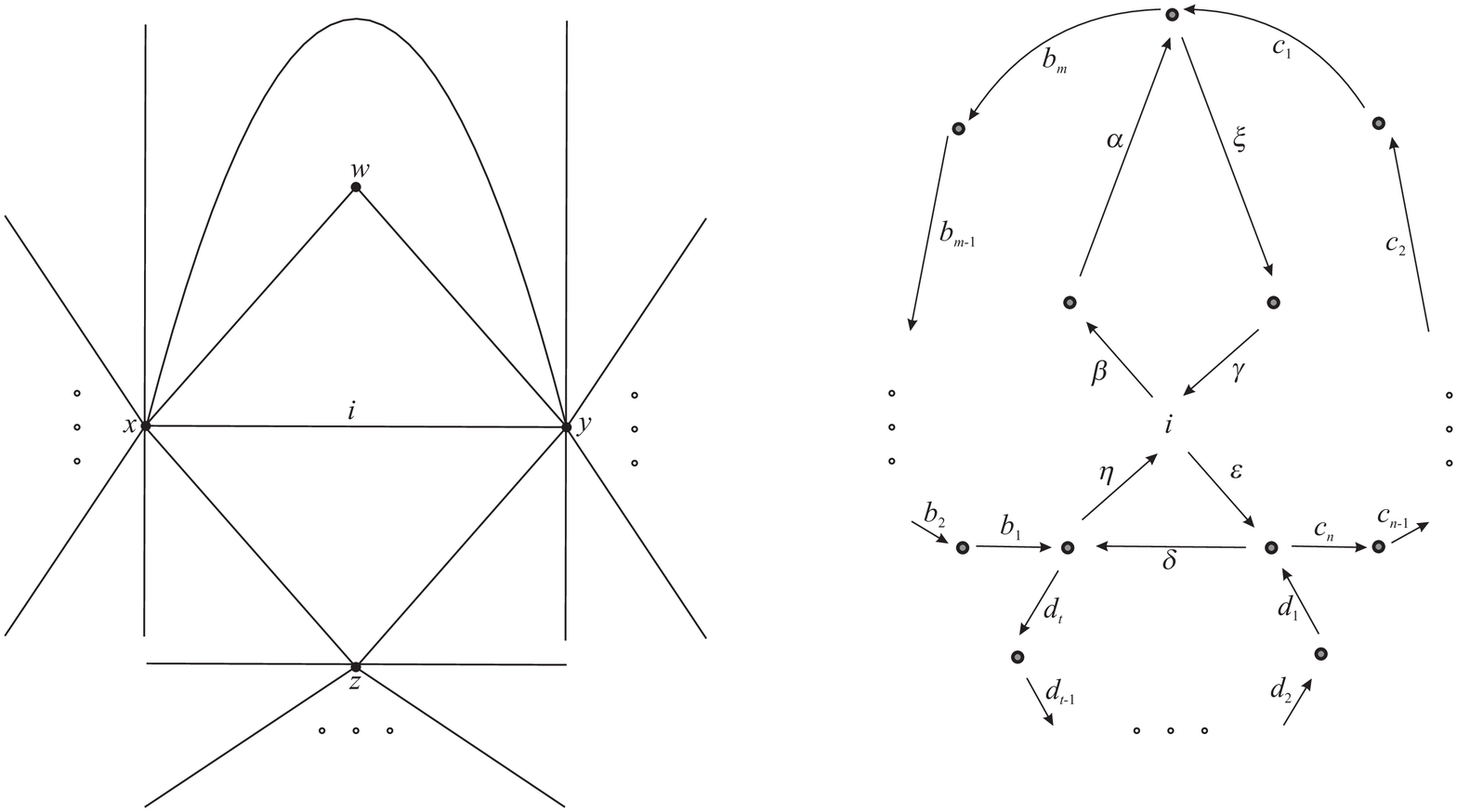}
        \end{figure}
Then (again abbreviating $b=b_1\ldots b_m, c=c_1\ldots c_n, d=d_1\ldots d_t$)
$$
\stau=\delta\varepsilon\eta-\frac{\xi\alpha\beta\gamma}{w}+x\alpha\beta\eta b+y\varepsilon\gamma\xi c+z\delta d+\sptau,
$$
with $\sptau\in\ratau$ involving none of the arrows $\alpha,\beta,\gamma,\delta,\varepsilon,\eta,\xi$. If we perform the premutation $\premuti$ on $\astau$, we get $\tildeastau$, where $\tildeatau$ is the arrow span of the quiver shown in Figure \ref{premutcasetwo},
% PREMUTATED SECOND GENERIC CASE
        \begin{figure}[!h]
                \caption{Case \ref{generic2}, QP-mutation process $\muti\qstau$}\label{premutcasetwo}
                \centering
                \includegraphics[scale=.35]{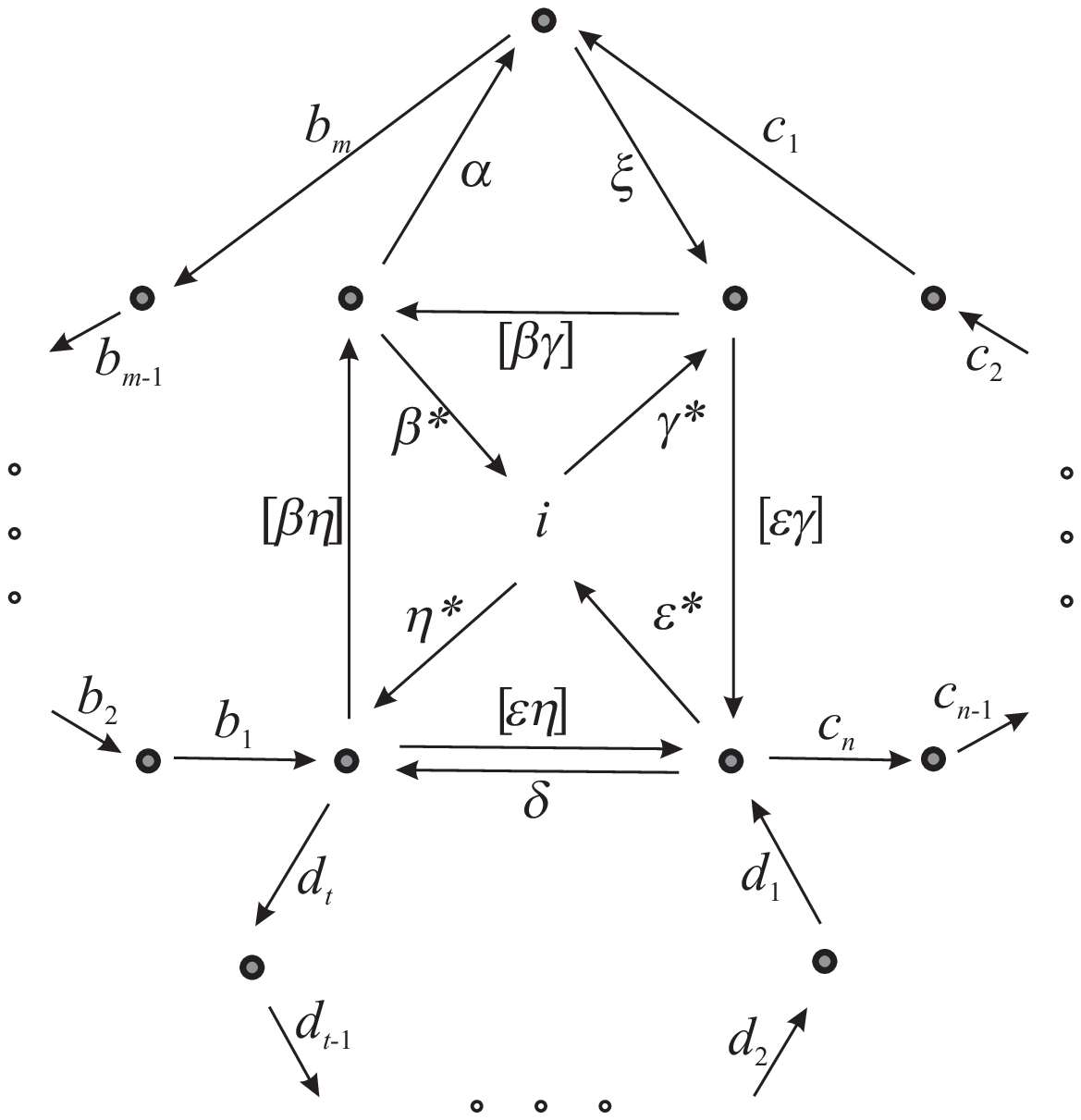}
        \end{figure}
and $\tildestau=\delta[\epsilon\eta]-\frac{\xi\alpha[\beta\gamma]}{w}+x\alpha[\beta\eta]b+y[\epsilon\gamma]\xi c+z\delta
d+\sptau+\eta^*\varepsilon^*[\varepsilon\eta]+\gamma^*\varepsilon^*[\varepsilon\gamma]+\gamma^*\beta^*[\beta\gamma]+\eta^*\beta^*[\beta\eta]\in
\rtildeatau$. The $R$-algebra automorphism $\varphi$ of $\rtildeatau$ whose action on the arrows is given by
$$
\delta\mapsto\delta-\eta^*\varepsilon^*,\ [\varepsilon\eta]\mapsto[\varepsilon\eta]-zd,$$
and the identity in the rest of the arrows, sends $\tildestau$ to
$$
\varphi(\tildestau)=\delta[\varepsilon\eta]-\frac{\xi\alpha[\beta\gamma]}{w}+x\alpha[\beta\eta]b+y[\varepsilon\gamma]\xi c-z\eta^*\varepsilon^*d+\sptau+\gamma^*\varepsilon^*[\varepsilon\gamma]+\gamma^*\beta^*[\beta\gamma]+\eta^*\beta^*[\beta\eta].
$$
Therefore, the reduced part $\muti\astau$ of $\tildeastau$ is (up to right-equivalence) the QP on the arrows span $\overatau$ obtained from $\tildeatau$ by deleting the arrows $\delta$ and $[\varepsilon\eta]$, with $\varphi(\tildestau)-\delta[\varepsilon\eta]$ as its potential.

On the other hand, $\sigma=f_i(\tau)$ and its quiver $\qsigma$ look as Figure \ref{flippedcasetwo},
% FLIPPED SECOND GENERIC CASE
        \begin{figure}[!h]
                \caption{Case \ref{generic2}, flip $\sigma=f_i(\tau)$}\label{flippedcasetwo}
                \centering
                \includegraphics[scale=.35]{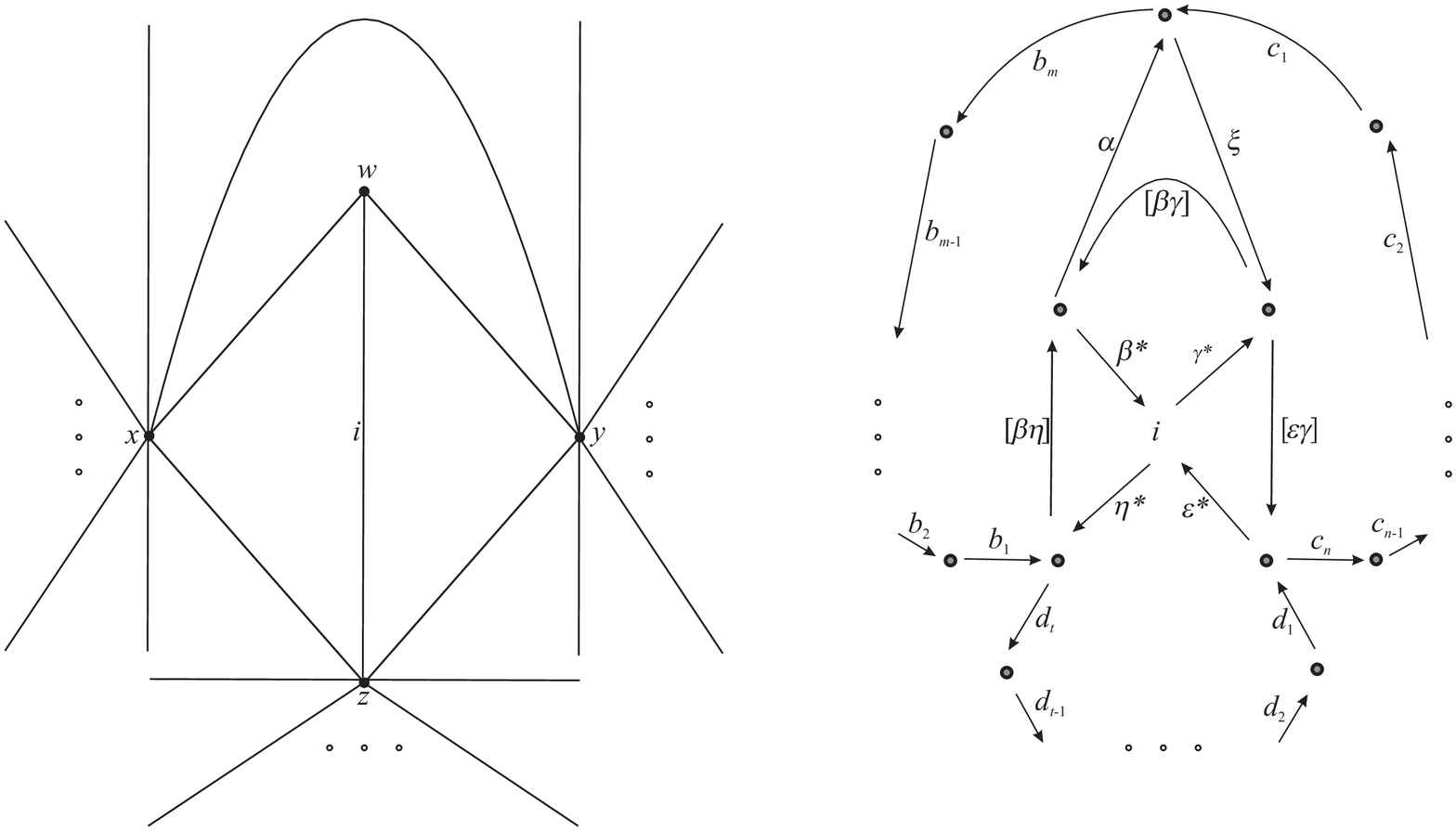}
        \end{figure}
and $\ssigma=\xi\alpha[\beta\gamma]+\eta^*\beta^*[\beta\eta]+\gamma^*\varepsilon^*[\varepsilon\gamma]+w\gamma^*\beta^*[\beta\gamma]+
x\alpha[\beta\eta]b+y[\varepsilon\gamma]\xi c+z\eta^*\varepsilon^*d+\spsigma$, with $\spsigma=\sptau$. Thus, the $R$-algebra isomorphism
$\psi:\roveratau\rightarrow\rasigma$ whose action on the arrows is given by
$$
\alpha\mapsto-\alpha,\ \beta^*\mapsto-\beta^*,\ \gamma^*\mapsto-\gamma^*,\ \varepsilon^*\mapsto-\varepsilon^*,\ [\beta\gamma]\mapsto w[\beta\gamma],\ [\beta\eta]\mapsto-[\beta\eta],$$
and the identity in the rest of the arrows, is a right-equivalence between $\muti\astau$ and $\assigma$.
\end{case}

\begin{case}\label{new6} (Fifth and eighth matchings of Figure \ref{matchings}) Assume that, around the arc $i$, $\tau$ looks like the configuration in Figure \ref{newcasesix}, with the arc $k$ not enclosing a self-folded triangle.
% FIRST GENERIC CASE
        \begin{figure}[!h]
                \caption{Case \ref{new6}, configuration of $\tau$ around $i$}\label{newcasesix}
                \centering
                \includegraphics[scale=.25]{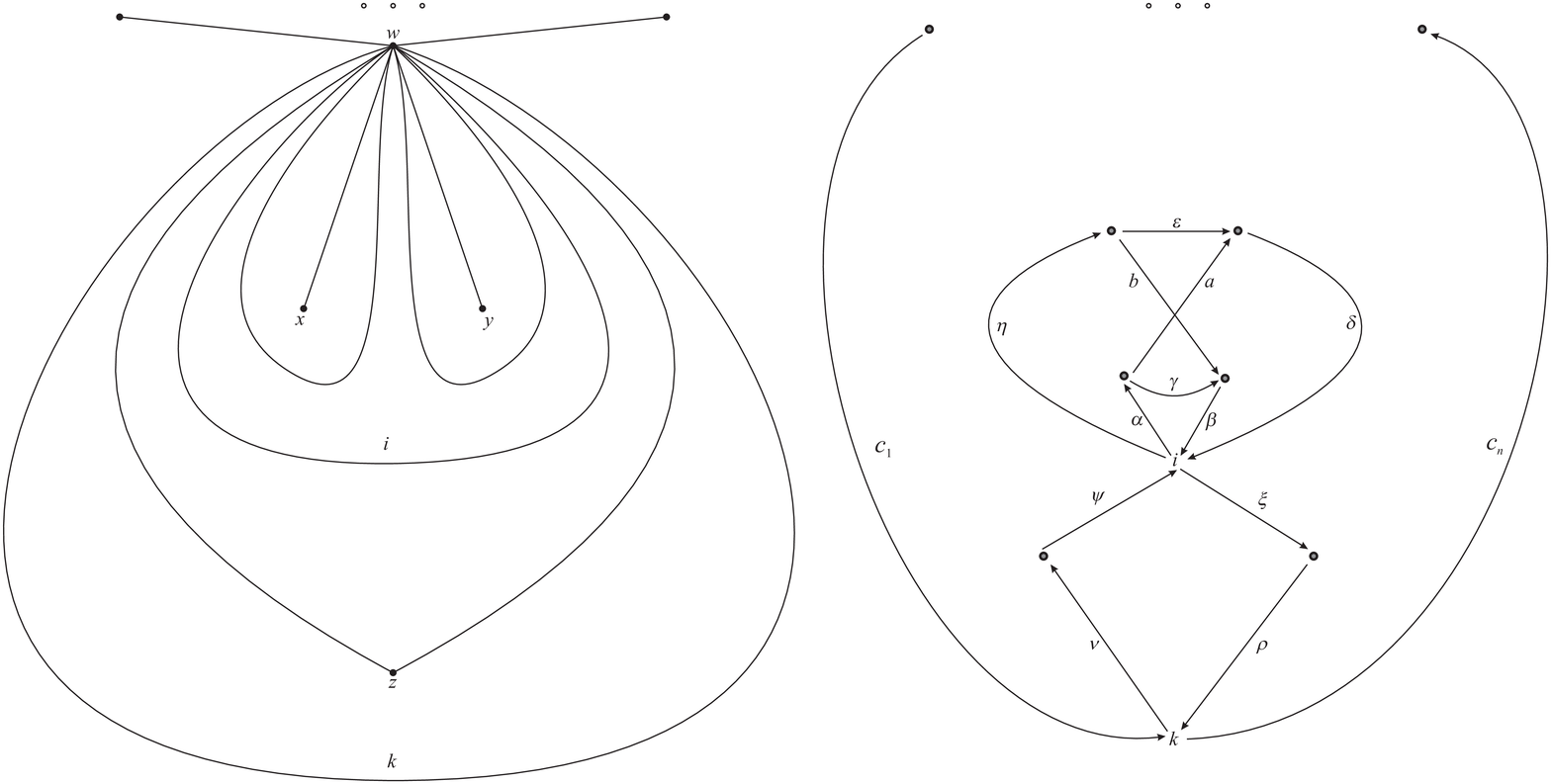}
        \end{figure}
Let us abbreviate $c=c_1\ldots c_n$. Then
$$
\stau=\alpha\beta\gamma+\frac{\eta\delta\varepsilon}{xy}-\frac{\eta\beta b}{x}-\frac{\alpha\delta a}{y}-\frac{\xi\psi\nu\rho}{z}
+w\eta\psi\nu c\rho\xi\delta\varepsilon+\sptau,
$$
with $\sptau\in\ratau$ involving none of the arrows $\alpha,\beta,\gamma,\delta,\varepsilon,\eta,\nu,\rho,\psi,\xi,a,b$. If we perform the premutation $\premuti$ on $\astau$, we get $\tildeastau$, where $\tildeatau$ is the arrow span of the quiver shown in Figure \ref{newpremutcasefour}
% PREMUTATED FIRST GENERIC CASE
        \begin{figure}[!h]
                \caption{Case \ref{new6}, QP-mutation process $\muti\qstau$}\label{newpremutcasesix}
                \centering
                \includegraphics[scale=.3]{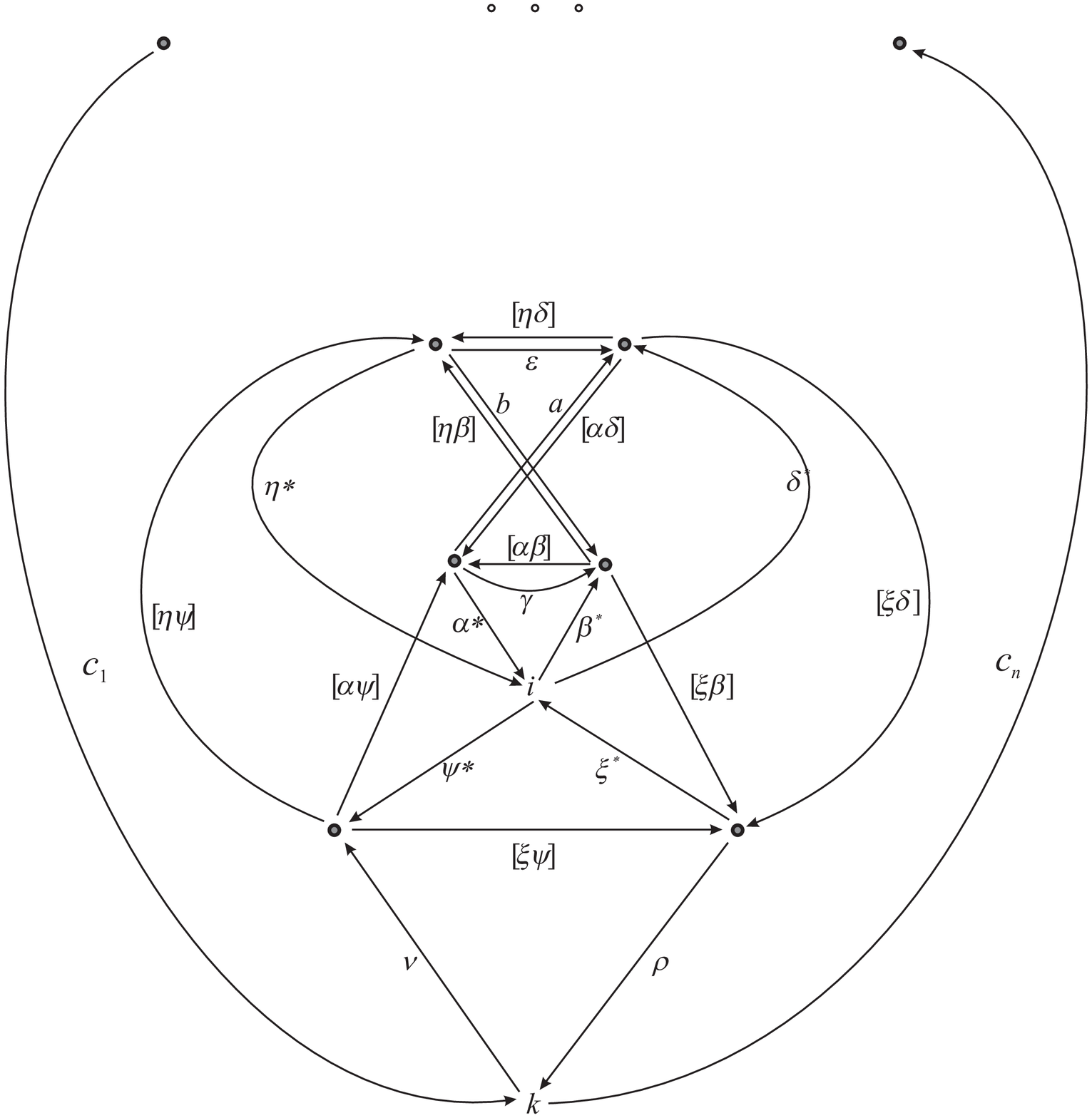}
        \end{figure}
and $\tildestau=[\alpha\beta]\gamma+\frac{[\eta\delta]\varepsilon}{xy}-\frac{[\eta\beta]b}{x}-\frac{[\alpha\delta]a}{y}-\frac{[\xi\psi]\nu\rho}{z}
+w[\eta\psi]\nu c\rho[\xi\delta]\varepsilon+\sptau+[\alpha\beta]\beta^*\alpha*+[\alpha\delta]\delta^*\alpha^*+[\alpha\psi]\psi^*\alpha^*
[\eta\beta]\beta^*\eta*+[\eta\delta]\delta^*\eta^*+[eta\psi]\psi^*\eta^*
[\xi\beta]\beta^*\xi*+[\xi\delta]\delta^*\xi^*+[\xi\psi]\psi^*\xi^*+\in\rtildeatau$. The $R$-algebra automorphism $\varphi$ of $\rtildeatau$ whose action on the arrows is given by
$$
\gamma\mapsto\gamma-\beta^*\alpha^*, \ [\eta\delta]\mapsto[\eta\delta]-wxy[\eta\psi]\nu c\rho[\xi\delta], \ \varepsilon\mapsto\varepsilon-xy\delta^*\eta^*, \ a\mapsto a+y\delta^*\alpha, \ b\mapsto b+x\beta^*\eta^*,
$$
and the identity in the rest of the arrows, sends $\tildestau$ to
$$
\varphi(\tildestau)=[\alpha\beta]\gamma+\frac{[\eta\delta]\varepsilon}{xy}-\frac{[\eta\beta]b}{x}-\frac{[\alpha\delta]a}{y}
-\frac{[\xi\psi]\nu\rho}{z}+\sptau+[\alpha\psi]\psi^*\alpha^*-wxy[\eta\psi]\nu c\rho[\xi\delta]\delta^*\eta*+[\eta\psi]\psi^*\eta^*+
$$
$$
+[\xi\beta]\beta^*\xi^*+[\xi\delta]\delta^*\xi^*+[\xi\psi]\psi^*\xi^*.
$$
Therefore, the reduced part $\muti\astau$ of $(\tildeatau,\varphi(\tildestau))$ is (up to right-equivalence) the QP on the arrow span $\overatau$ obtained from $\tildeatau$ by deleting the arrows $[\alpha\beta],\gamma,[\eta\delta],\varepsilon,[\eta\beta],b,[\alpha\delta]$ and $a$, with $\varphi(\tildestau)-[\alpha\beta]\gamma-\frac{[\eta\delta]\varepsilon}{xy}+\frac{[\eta\beta]b}{x}+\frac{[\alpha\delta]a}{y}$ as its potential.

On the other hand, $\sigma=f_i(\tau)$ and its quiver $\qsigma$ look as Figure \ref{newflippedcasesix},
% FLIPPED FIRST GENERIC CASE
        \begin{figure}[!h]
                \caption{Case \ref{new6}, flip $\sigma=f_i(\tau)$}\label{newflippedcasesix}
                \centering
                \includegraphics[scale=.25]{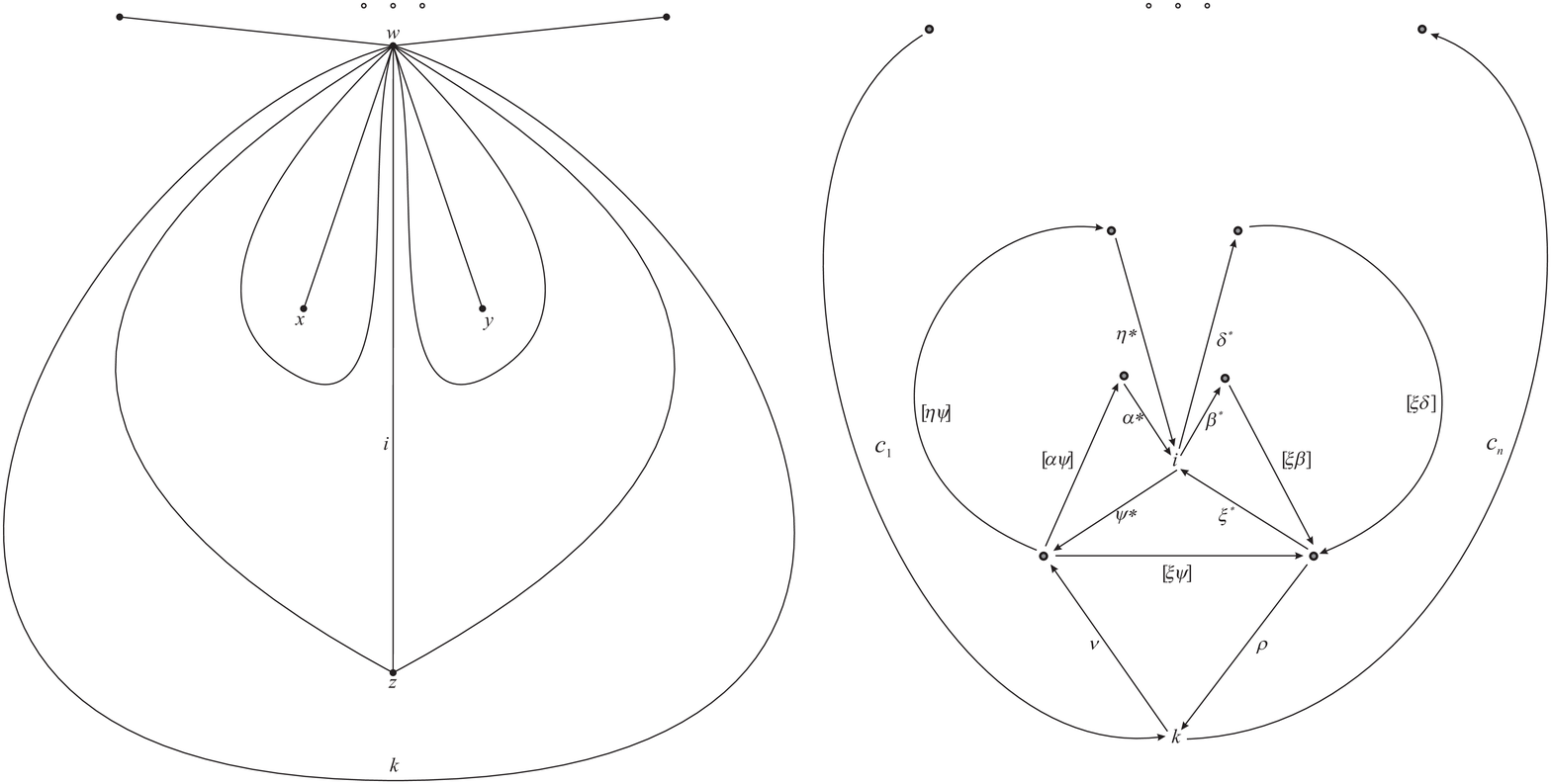}
        \end{figure}
and $\ssigma=[\alpha\psi]\psi^*\alpha^*+[\xi\beta]\beta^*\xi^*+[\xi\psi]\nu\rho-\frac{[\eta\psi]\nu\rho}{x}-\frac{[\xi\delta]\delta^*\xi^*}{y}
+z[\xi\psi]\psi^*\xi^*+w[\eta\psi]\nu c\rho[\xi\delta]\delta^*\eta^*+\spsigma$, with $\spsigma=\sptau$. Thus the $R$-algebra isomorphism
$\psi:\roveratau\rightarrow\rasigma$ whose action on the arrows is given by
$$
\alpha^*\mapsto-\alpha^*, \ \psi^*\mapsto-\psi^*, \ [\eta\psi]\mapsto\frac{[\eta\psi]}{x}, \ [\xi\delta]\mapsto-\frac{[\xi\delta]}{y}, \ [\xi\psi]\mapsto-z[\xi\psi]
$$
and the identity
in the rest of the arrows, is a right-equivalence between $\muti\astau$ and $\assigma$.
\end{case}

\begin{case}\label{generic6} (Sixth matching of Figure \ref{matchings}) Assume that, around $i$, the triangulation $\tau$ looks like the configuration shown in Figure \ref{casesix}.
% SIXTH GENERIC CASE
        \begin{figure}[!h]
                \caption{Case \ref{generic6}, configuration of $\tau$ around $i$}\label{casesix}
                \centering
                \includegraphics[scale=.35]{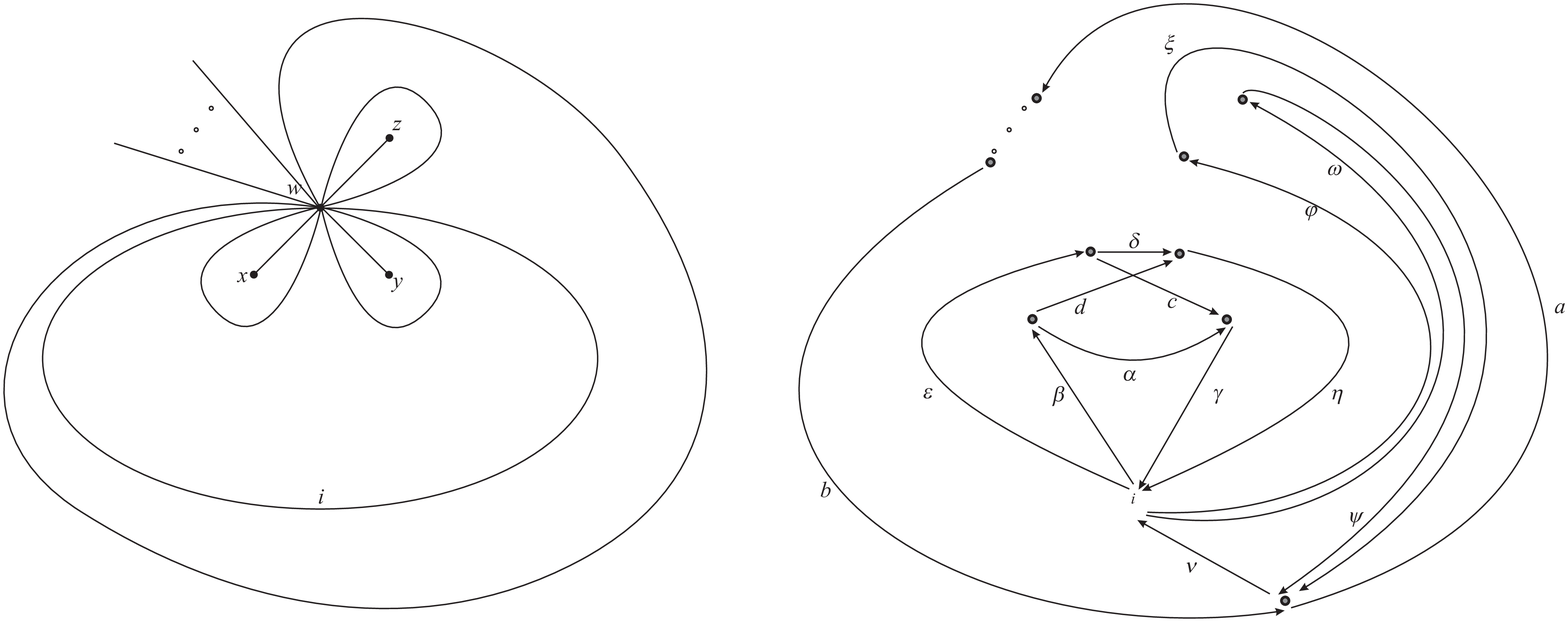}
        \end{figure}
Then
$$
\stau=\alpha\beta\gamma+\frac{\delta\varepsilon\eta}{xy}+\psi\omega\nu-\frac{c\varepsilon\gamma}{x}-\frac{d\beta\eta}{y}-\frac{\xi\varphi\nu}{z}
+w\delta\varepsilon\nu b\ldots a\xi\varphi\eta+\sptau,
$$
with $\sptau\in\ratau$ involving none of the arrows $\alpha$, $\beta$, $\gamma$, $\delta$, $\varepsilon$, $\eta$, $\nu$, $\varphi$, $\psi$, $\omega$, $\xi$, $c$ and $d$. If we perform the premutation $\premuti$ on $\astau$ we get $\tildeastau$, where $\tildeatau$ is the arrow span of the quiver shown in Figure \ref{premutcasesix},
% PREMUTATED SIXTH GENERIC CASE
        \begin{figure}[!h]
                \caption{Case \ref{generic6}, QP-mutation process $\muti\qstau$}\label{premutcasesix}
                \centering
                \includegraphics[scale=.35]{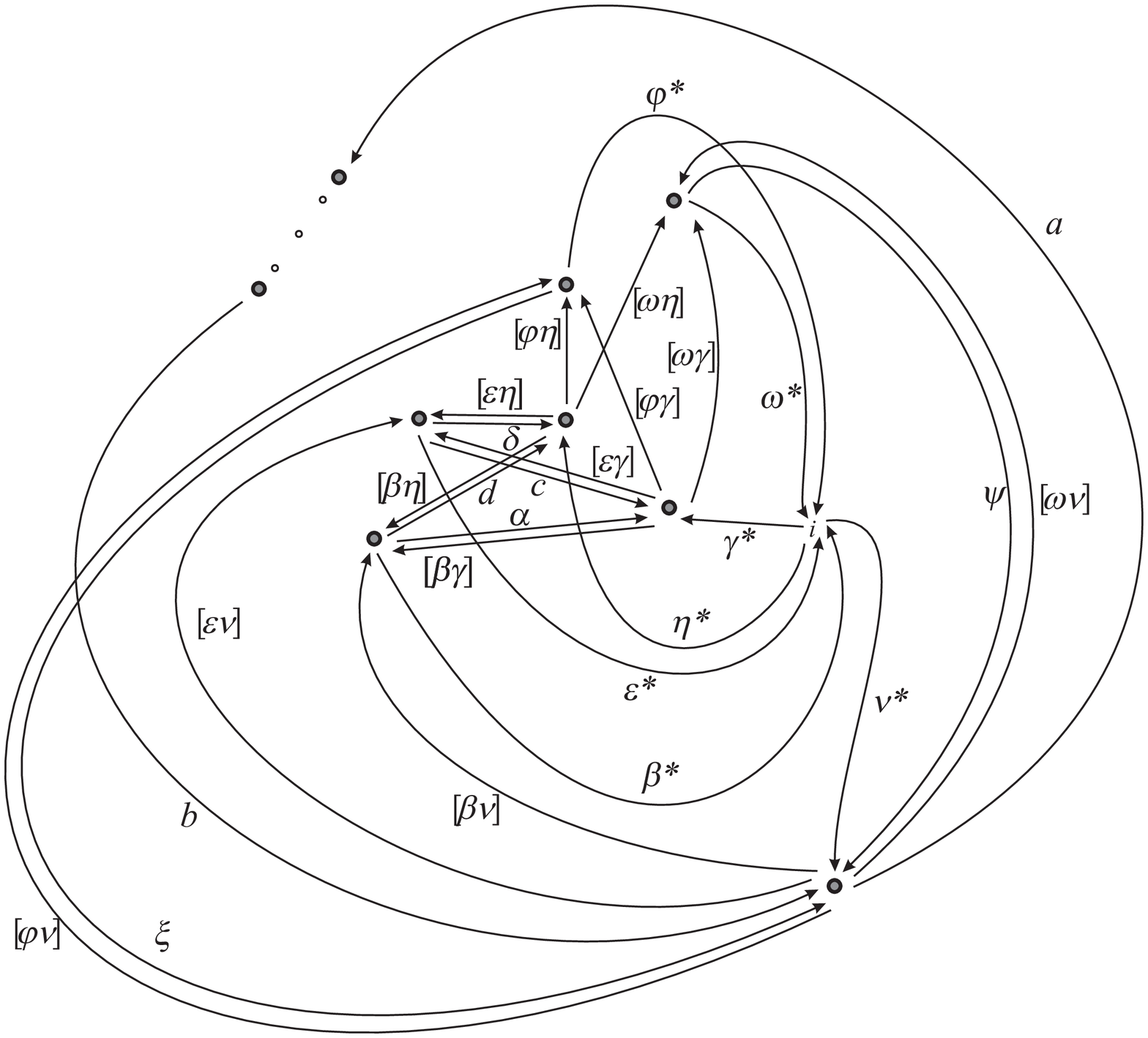}
        \end{figure}
and $\tildestau=\alpha[\beta\gamma]+\frac{\delta[\varepsilon\eta]}{xy}+\psi[\omega\nu]-\frac{c[\varepsilon\gamma]}{x}-\frac{d[\beta\eta]}{y}
-\frac{\xi[\varphi\nu]}{z}+w\delta[\varepsilon\nu]b\ldots a\xi[\varphi\eta]+\sptau+
\gamma^*\beta^*[\beta\gamma]+\gamma^*\varepsilon^*[\varepsilon\gamma]+\gamma^*\varphi^*[\varphi\gamma]+\gamma^*\omega^*[\omega\gamma]+
\eta^*\beta^*[\beta\eta]    +\eta^*\varepsilon^*[\varepsilon\eta]    +\eta^*\varphi^*[\varphi\eta]    +\eta^*\omega^*[\omega\eta]+
\nu^*\beta^*[\beta\nu]      +\nu^*\varepsilon^*[\varepsilon\nu]      +\nu^*\varphi^*[\varphi\nu]      +\nu^*\omega^*[\omega\nu]$. The $R$-algebra automorphism $\varphi_1$ of $\rtildeatau$ whose action on the arrows is given by
$$
\alpha\mapsto\alpha-\gamma^*\beta^*,\ \delta\mapsto\delta-xy\eta^*\varepsilon^*,\ [\varepsilon\eta]\mapsto[\varepsilon\eta]-wxy[\varepsilon\nu]b\ldots a\xi[\varphi\eta],
$$
$$
\psi\mapsto\psi-\nu^*\omega^*,\ c\mapsto c+x\gamma^*\varepsilon^*,\ d\mapsto d+y\eta^*\beta^*,\ \xi\mapsto\xi+z\nu^*\varphi^*,
$$
and the identity in the rest of the arrows, sends $\tildestau$ to
$$
\varphi_1(\tildestau)=\alpha[\beta\gamma]+\frac{\delta[\varepsilon\eta]}{xy}+\psi[\omega\nu]-\frac{c[\varepsilon\gamma]}{x}
-\frac{d[\beta\eta]}{y}-\frac{\xi[\varphi\nu]}{z}+wz\delta[\varepsilon\nu]b\ldots a\nu^*\varphi^*[\varphi\eta]
-wxyz\eta^*\varepsilon^*[\varepsilon\nu]b\ldots a\nu^*\varphi^*[\varphi\eta]
$$
$$
+\sptau+\gamma^*\varphi^*[\varphi\gamma]
+\gamma^*\omega^*[\omega\gamma]-wxy\eta^*\varepsilon^*[\varepsilon\nu]b\ldots a\xi[\varphi\eta]
+\eta^*\varphi^*[\varphi\eta]+\eta^*\omega^*[\omega\eta]+\nu^*\beta^*[\beta\nu]+\nu^*\varepsilon^*[\varepsilon\nu],
$$
which is cyclically equivalent to
$$
T=\alpha[\beta\gamma]+\frac{\delta[\varepsilon\eta]}{xy}+\psi[\omega\nu]-\frac{c[\varepsilon\gamma]}{x}
-\frac{d[\beta\eta]}{y}-\frac{\xi[\varphi\nu]}{z}+wz\delta[\varepsilon\nu]b\ldots a\nu^*\varphi^*[\varphi\eta]
-wxyz\eta^*\varepsilon^*[\varepsilon\nu]b\ldots a\nu^*\varphi^*[\varphi\eta]
$$
$$
+\sptau+\gamma^*\varphi^*[\varphi\gamma]
+\gamma^*\omega^*[\omega\gamma]-wxy\xi[\varphi\eta]\eta^*\varepsilon^*[\varepsilon\nu]b\ldots a
+\eta^*\varphi^*[\varphi\eta]+\eta^*\omega^*[\omega\eta]+\nu^*\beta^*[\beta\nu]+\nu^*\varepsilon^*[\varepsilon\nu],
$$
which in turn is sent to
$$
\varphi_2(T)=\alpha[\beta\gamma]+\frac{\delta[\varepsilon\eta]}{xy}+\psi[\omega\nu]-\frac{c[\varepsilon\gamma]}{x}
-\frac{d[\beta\eta]}{y}-\frac{\xi[\varphi\nu]}{z}-wxyz\eta^*\varepsilon^*[\varepsilon\nu]b\ldots a\nu^*\varphi^*[\varphi\eta]
$$
$$
+\sptau+\gamma^*\varphi^*[\varphi\gamma]
+\gamma^*\omega^*[\omega\gamma]+\eta^*\varphi^*[\varphi\eta]+\eta^*\omega^*[\omega\eta]+\nu^*\beta^*[\beta\nu]+\nu^*\varepsilon^*[\varepsilon\nu]
$$
by the $R$-algebra automorphism $\varphi_2$ of $\rtildeatau$ whose action on the arrows is given by
$$
[\varepsilon\eta]\mapsto[\varepsilon\eta]-wxyz[\varepsilon\nu]b\ldots a\nu^*\varphi^*[\varphi\eta],\
[\varphi\nu]\mapsto[\varphi\nu]-wxyz[\varphi\eta]\eta^*\varepsilon^*[\varepsilon\nu]b\ldots a.
$$
Therefore, the reduced part $\muti\astau$ of $\tildeastau$ is (up to right-equivalence) the QP on the arrows span $\overatau$ obtained from $\tildeatau$ by deleting the arrows $\alpha$, $[\beta\gamma]$, $\delta$, $[\varepsilon\eta]$, $\psi$, $[\omega\nu]$, $c$, $[\varepsilon\gamma]$, $d$, $[\beta\eta]$, $\xi$ and $[\varphi\nu]$, with $\varphi_2(T)-\alpha[\beta\gamma]-\frac{\delta[\varepsilon\eta]}{xy}-\psi[\omega\nu]+\frac{c[\varepsilon\gamma]}{x}
+\frac{d[\beta\eta]}{y}+\frac{\xi[\varphi\nu]}{z}$ as its potential.

On the other hand, $\sigma=f_i(\tau)$ and its quiver $\qsigma$ look as Figure \ref{flippedcasesix},
% FLIPPED SIXTH GENERIC CASE
        \begin{figure}[!h]
                \caption{Case \ref{generic6}, flip $\sigma=f_i(\tau)$}\label{flippedcasesix}
                \centering
                \includegraphics[scale=.35]{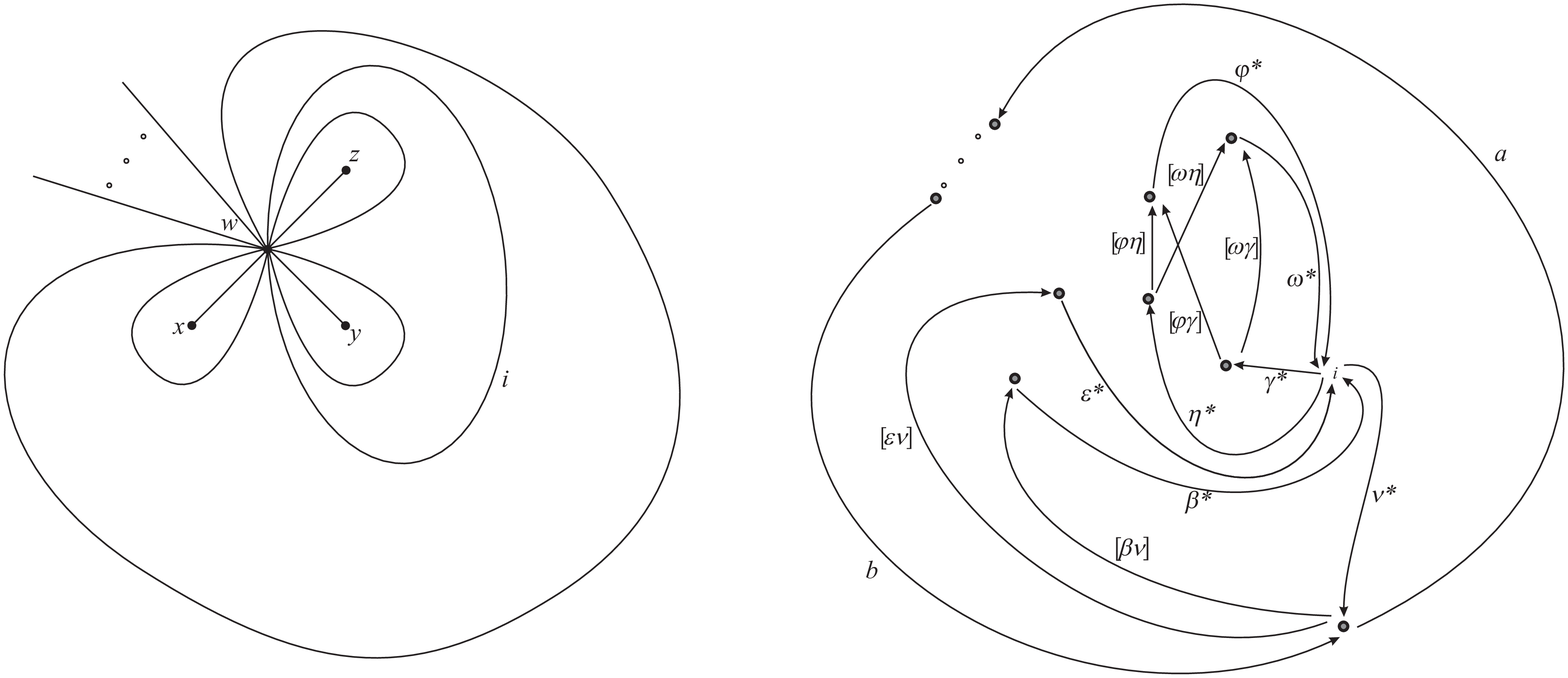}
        \end{figure}
and $\ssigma=\gamma^*\omega^*[\omega\gamma]+\frac{\eta^*\varphi^*[\varphi\eta]}{yz}+\nu^*\beta^*[\beta\nu]
-\frac{\nu^*\varepsilon^*[\varepsilon\nu]}{x}-\frac{\eta^*\omega^*[\omega\eta]}{y}-\frac{\gamma^*\varphi^*[\varphi\gamma]}{z}
+w\eta^*\varepsilon^*[\varepsilon\nu]b\ldots a\nu^*\varphi^*[\varphi\eta]+\spsigma$, with $\spsigma=\sptau$. Thus, the $R$-algebra isomorphism
$\psi:\roveratau\rightarrow\rasigma$ whose action on the arrows is given by
$$
\beta^*\mapsto-x\beta^*,\ \eta^*\mapsto-\frac{\eta^*}{y},\ \nu^*\mapsto-\frac{\nu^*}{x},\ \varphi^*\mapsto-\frac{\varphi^*}{z}
$$
and the identity in the rest of the arrows, is a right-equivalence between $\muti\astau$ and $\assigma$.
\end{case}

\begin{case}\label{new4} (Seventh matching of Figure \ref{matchings}) Assume that, around the arc $i$, $\tau$ looks like the configuration in Figure \ref{newcasefour}, with none of the arcs $j$ and $k$ enclosing a self-folded triangle.
% FIRST GENERIC CASE
        \begin{figure}[!h]
                \caption{Case \ref{new4}, configuration of $\tau$ around $i$}\label{newcasefour}
                \centering
                \includegraphics[scale=.25]{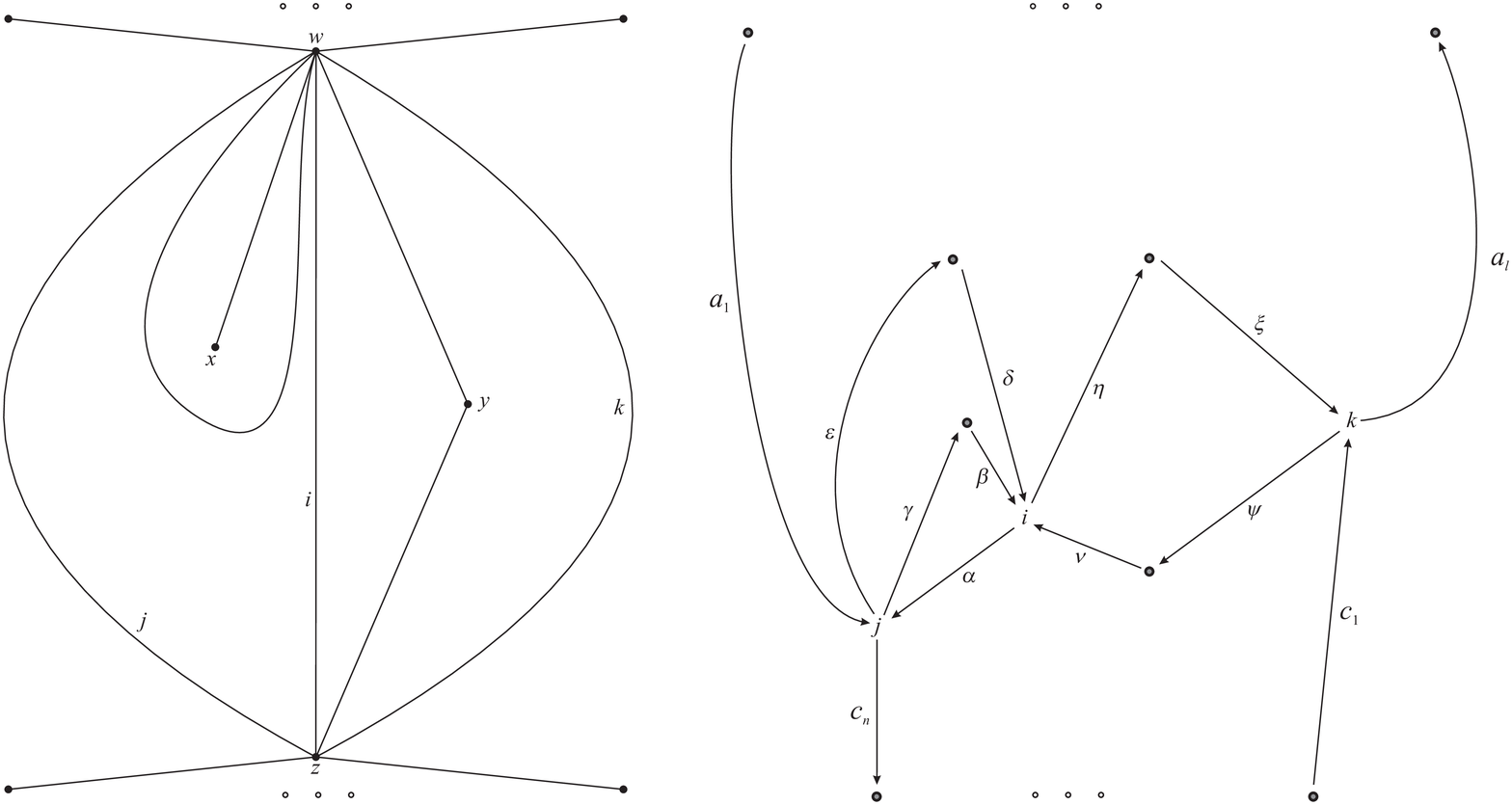}
        \end{figure}
Let us abbreviate $a=a_1\ldots a_l$, $c=c_1\ldots c_n$. Then
$$
\stau=\alpha\beta\gamma-\frac{\alpha\delta\varepsilon}{x}-\frac{\eta\nu\psi\xi}{y}+w\eta\delta\varepsilon a\xi+z\alpha\nu\psi c\sptau,
$$
with $\sptau\in\ratau$ involving none of the arrows $\alpha,\beta,\gamma,\delta,\varepsilon,\eta,\nu,\psi,\xi$. If we perform the premutation $\premuti$ on $\astau$, we get $\tildeastau$, where $\tildeatau$ is the arrow span of the quiver shown in Figure \ref{newpremutcasefour}
% PREMUTATED FIRST GENERIC CASE
        \begin{figure}[!h]
                \caption{Case \ref{new4}, QP-mutation process $\muti\qstau$}\label{newpremutcasefour}
                \centering
                \includegraphics[scale=.25]{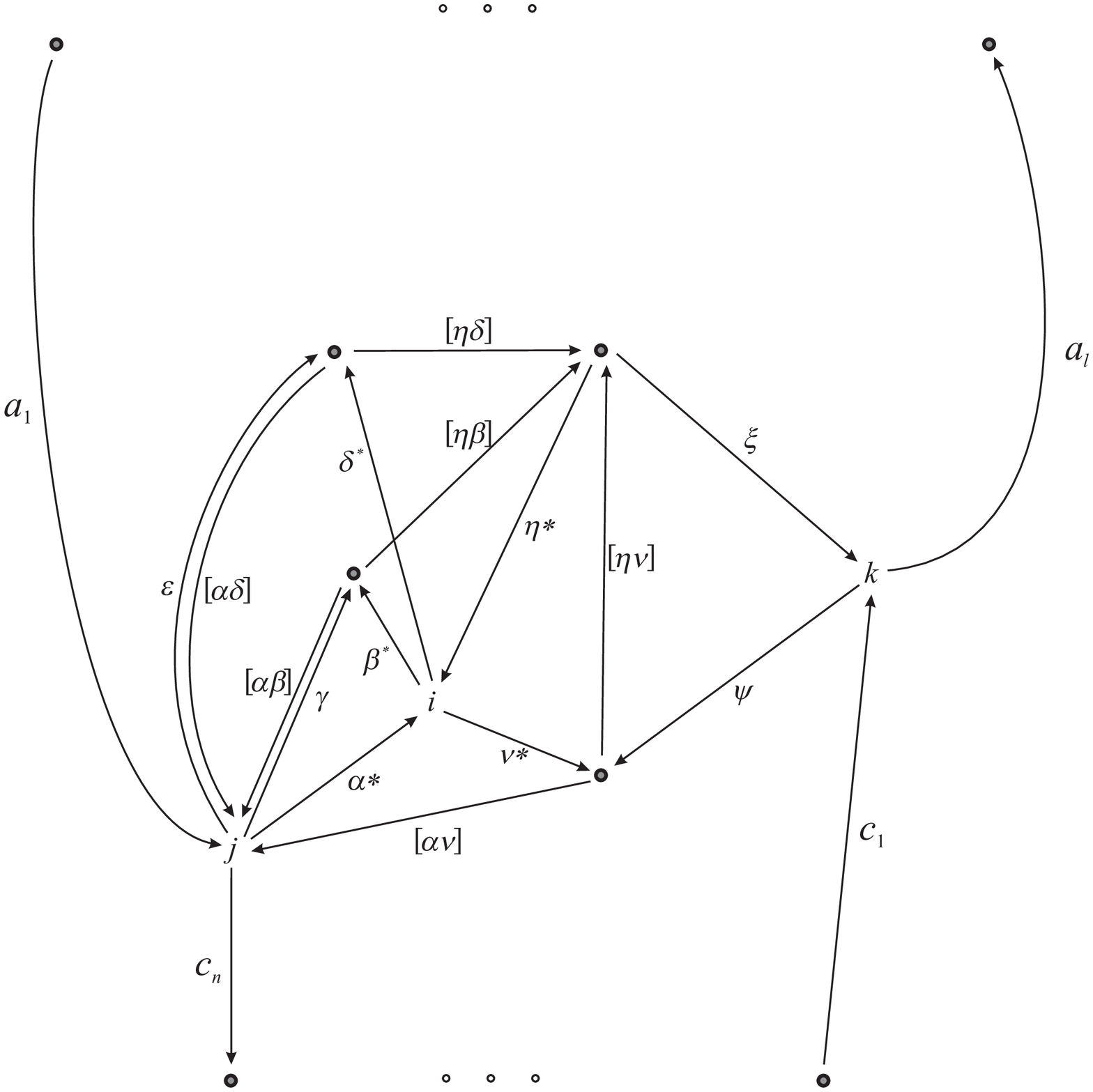}
        \end{figure}
and $\tildestau=[\alpha\beta]\gamma-\frac{[\alpha\delta]\varepsilon}{x}-\frac{[\eta\nu]\psi\xi}{y}+w[\eta\delta]\varepsilon a\xi+z[\alpha\nu]\psi c+\sptau+[\alpha\beta]\beta^*\alpha^*+[\alpha\delta]\delta^*\alpha^*+[\alpha\nu]\nu^*\alpha^*+[\eta\beta]\beta^*\eta^*+
[\eta\delta]\delta^*\eta^*+[\eta\nu]\nu^*\eta^*\in\rtildeatau$. The $R$-algebra automorphism $\varphi$ of $\rtildeatau$ whose action on the arrows is given by
$$
\gamma\mapsto\gamma-\beta^*\alpha^*, \ [\alpha\delta]\mapsto[\alpha\delta]+wxa\xi[\delta\eta], \ \varepsilon\mapsto\varepsilon+x\delta^*\alpha^*,
$$
and the identity in the rest of the arrows, sends $\tildestau$ to
$$
\varphi(\tildestau)=[\alpha\beta]\gamma-\frac{[\alpha\delta]\varepsilon}{x}-\frac{[\eta\nu]\psi\xi}{y}+z[\alpha\nu]\psi c+\sptau+wxa\xi[\eta\delta]\delta^*\alpha^*+[\alpha\nu]\nu^*\alpha^*+[\eta\beta]\beta^*\eta^*+[\eta\delta]\delta^*\eta^*+[\eta\nu]\nu^*\eta^*.
$$
Therefore, the reduced part $\muti\astau$ of $(\tildeatau,\varphi(\tildestau))$ is (up to right-equivalence) the QP on the arrow span $\overatau$ obtained from $\tildeatau$ by deleting the arrows $[\alpha\beta],\gamma,[\alpha\delta]$ and $\varepsilon$, with $\varphi(\tildestau)-[\alpha\beta]\gamma+\frac{[\alpha\delta]\varepsilon}{x}$ as its potential.

On the other hand, $\sigma=f_i(\tau)$ and its quiver $\qsigma$ look as Figure \ref{newflippedcasefour},
% FLIPPED FIRST GENERIC CASE
        \begin{figure}[!h]
                \caption{Case \ref{new4}, flip $\sigma=f_i(\tau)$}\label{newflippedcasefour}
                \centering
                \includegraphics[scale=.25]{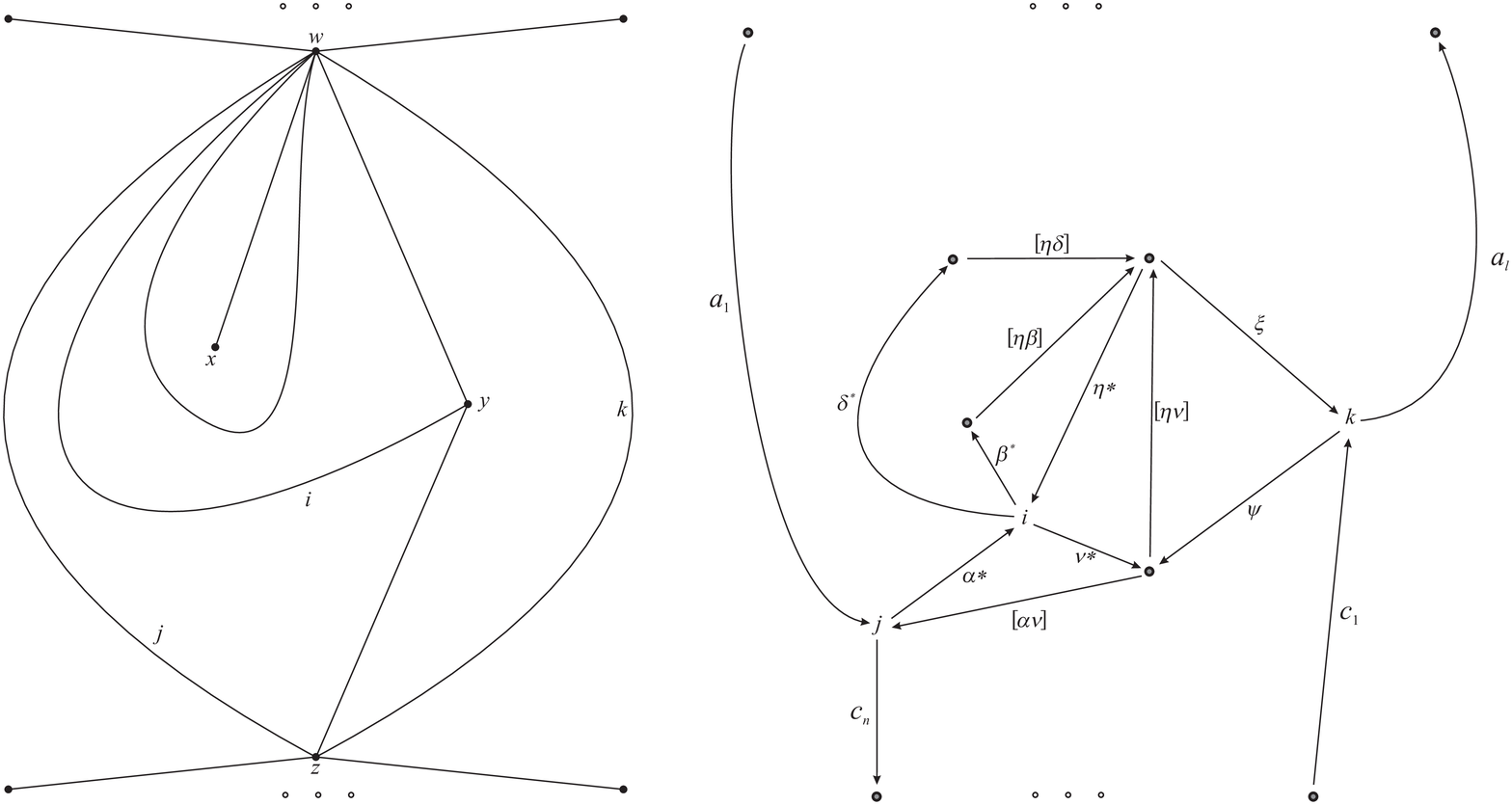}
        \end{figure}
and $\ssigma=[\alpha\nu]\nu^*\alpha^*+[\eta\nu]\psi\xi+[\eta\beta]\beta^+\eta^*-\frac{[\eta\delta]\delta^*\eta^*}{x}+y[\eta\nu]\nu^*\eta^*+
wa\xi[\eta\delta]\delta^*\alpha^*+z[\alpha\nu]\psi c+\spsigma$, with $\spsigma=\sptau$. Thus the $R$-algebra isomorphism
$\psi:\roveratau\rightarrow\rasigma$ whose action on the arrows is given by
$$
\beta^*\mapsto-\beta^*, \ \eta\mapsto-\eta^*, [\eta\delta]\mapsto\frac{[\eta\delta]}{x}, \ [\eta\nu]\mapsto-y[\eta\nu],
$$
and the identity
in the rest of the arrows, is a right-equivalence between $\muti\astau$ and $\assigma$.
\end{case}

\begin{case}\label{new5} (Ninth matching of Figure \ref{matchings}) Assume that, around the arc $i$, $\tau$ looks like the configuration in Figure \ref{newcasefive}.
% FIRST GENERIC CASE
        \begin{figure}[!h]
                \caption{Case \ref{new5}, configuration of $\tau$ around $i$}\label{newcasefive}
                \centering
                \includegraphics[scale=.25]{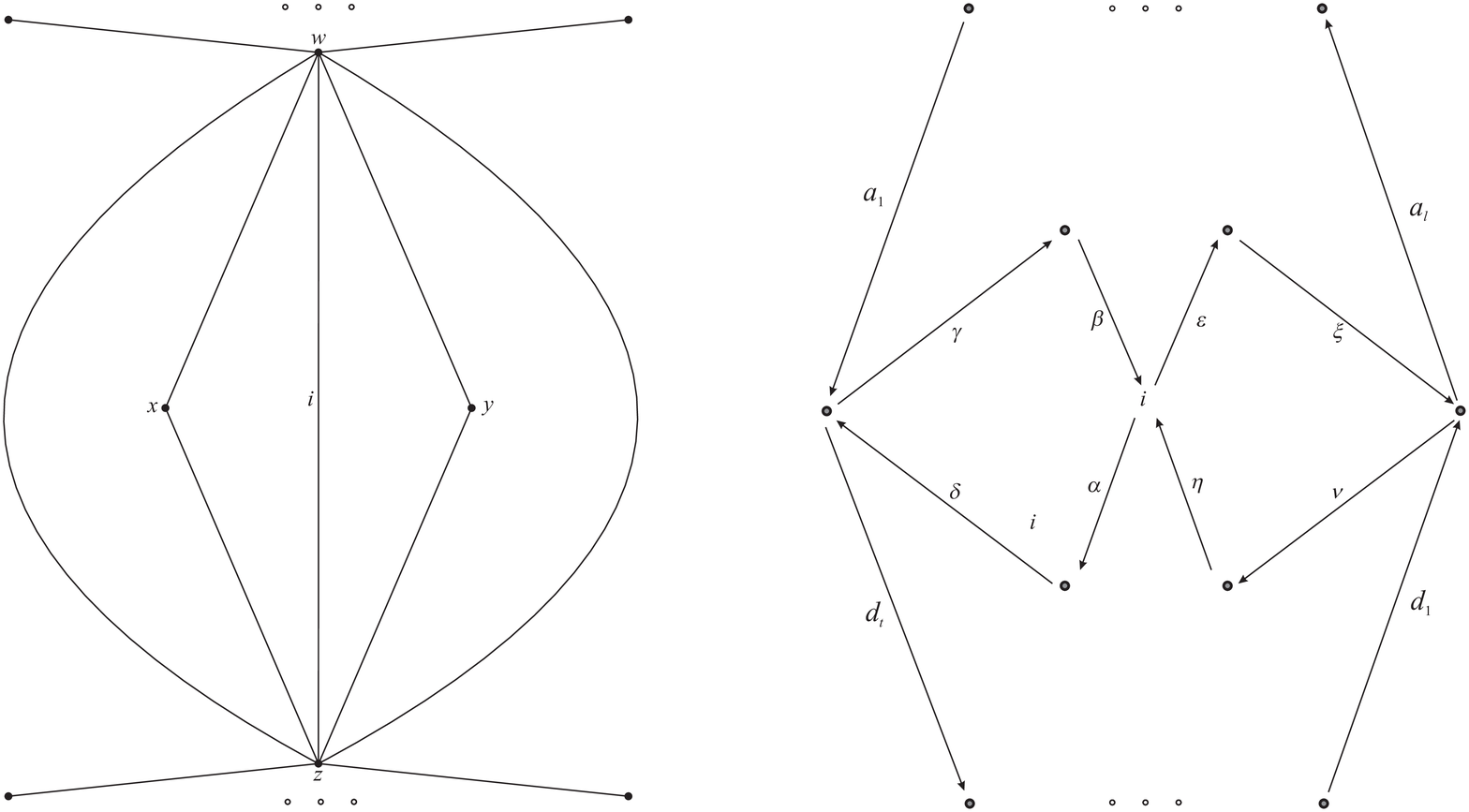}
        \end{figure}
Let us abbreviate $a=a_1\ldots a_l$, $d=d_1\ldots d_t$. Then
$$
\stau=-\frac{\alpha\beta\gamma\delta}{x}-\frac{\varepsilon\eta\nu\xi}{y}+wa\xi\varepsilon\beta\gamma+zd\alpha\delta\alpha\eta\nu+\sptau,
$$
with $\sptau\in\ratau$ involving none of the arrows $\alpha,\beta,\gamma,\delta,\varepsilon,\eta,\nu,\psi,\xi$. If we perform the premutation $\premuti$ on $\astau$, we get $\tildeastau$, where $\tildeatau$ is the arrow span of the quiver shown in Figure \ref{newpremutcasefive}
% PREMUTATED FIRST GENERIC CASE
        \begin{figure}[!h]
                \caption{Case \ref{new5}, QP-mutation process $\muti\qstau$}\label{newpremutcasefive}
                \centering
                \includegraphics[scale=.25]{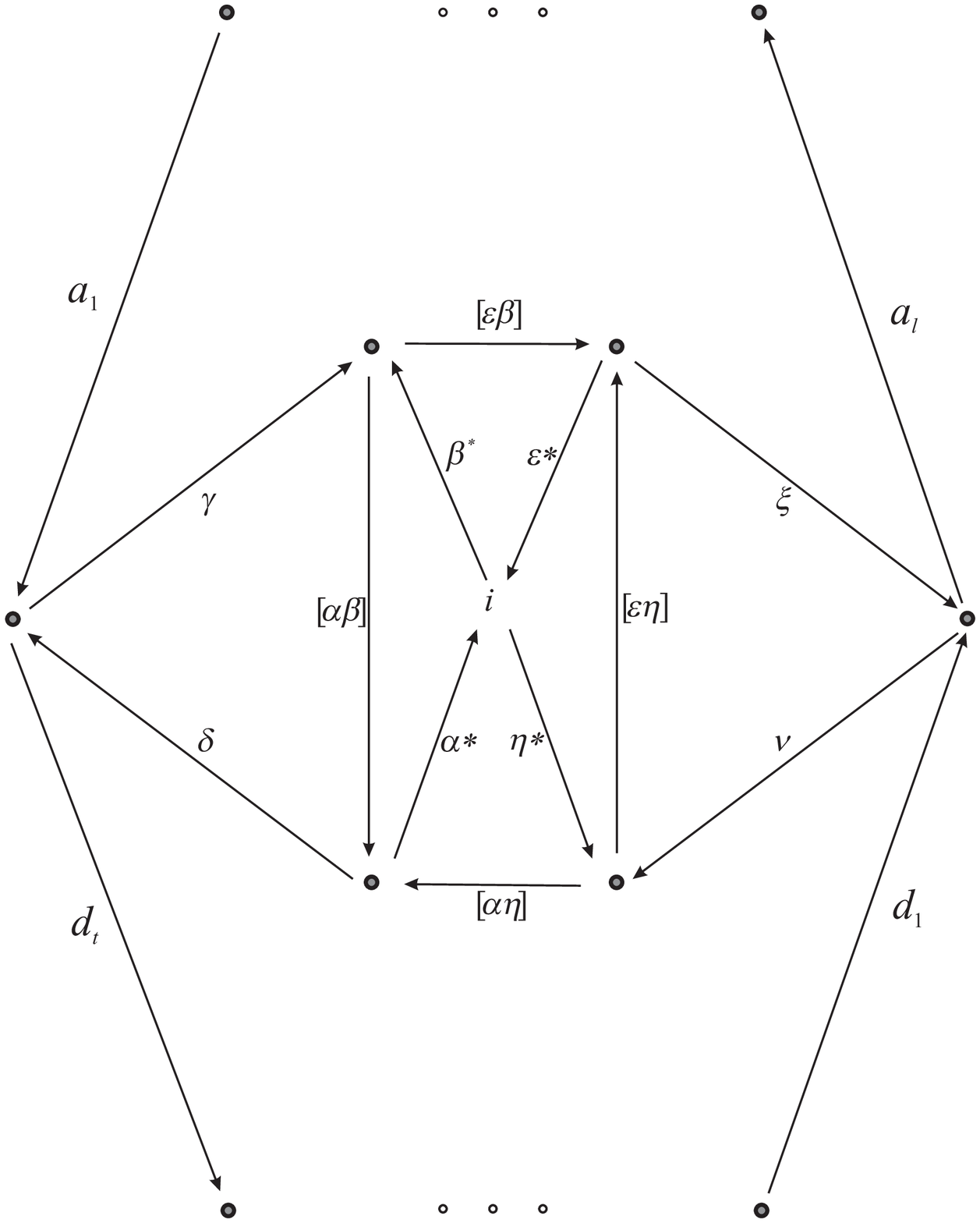}
        \end{figure}
and $\tildestau=-\frac{[\alpha\beta]\gamma\delta}{x}-\frac{[\varepsilon\eta]\nu\xi}{y}+wa\xi[\varepsilon\beta]\gamma+zd\delta[\alpha\eta]\nu+
\sptau+[\alpha\beta]\beta^*\alpha^*+[\alpha\eta]\eta^*\alpha^*+[\varepsilon\beta]\beta^*\varepsilon^*+[\varepsilon\eta]\eta^*\varepsilon^*
\in\rtildeatau$. Since $\tildeastau$ is already reduced, we have $\muti\astau?\tildeastau$.

On the other hand, $\sigma=f_i(\tau)$ and its quiver $\qsigma$ look as Figure \ref{newflippedcasefive},
% FLIPPED FIRST GENERIC CASE
        \begin{figure}[!h]
                \caption{Case \ref{new5}, flip $\sigma=f_i(\tau)$}\label{newflippedcasefive}
                \centering
                \includegraphics[scale=.25]{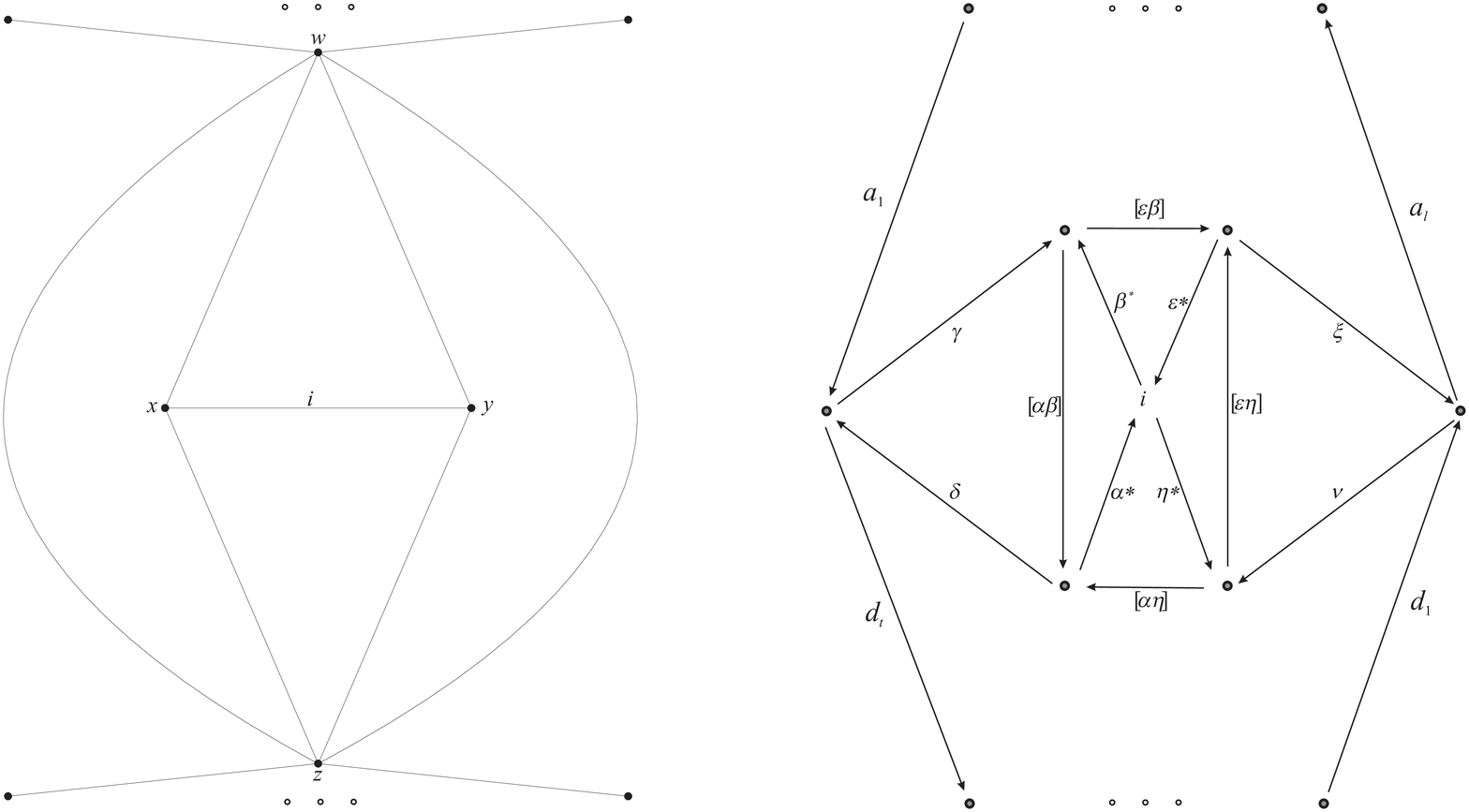}
        \end{figure}
and $\ssigma=[\alpha\eta]\eta^*\alpha^*+[\varepsilon\beta]\beta^*\varepsilon^*+[\alpha\beta]\gamma\delta+[\varepsilon\eta]\nu\xi+
x[\alpha\beta]\beta^*\alpha^*+y[\varepsilon\eta]\eta^*\varepsilon^*+wa\xi[\varepsilon\beta]\gamma+zd\delta[\alpha\eta]\nu+\spsigma$, with $\spsigma=\sptau$. Thus the $R$-algebra isomorphism
$\psi:\roveratau\rightarrow\rasigma$ whose action on the arrows is given by
$$
\alpha^*\mapsto-\alpha^*, \ \eta^*\mapsto-\eta^*, \ [\alpha\beta]\mapsto-x[\alpha\beta], \ [\varepsilon\eta]\mapsto-y[\varepsilon\eta],
$$
and the identity
in the rest of the arrows, is a right-equivalence between $\muti\astau$ and $\assigma$.
\end{case}

\begin{case}\label{genericthree} (Flip inside the fourth puzzle piece of Figure \ref{puzzlepieces}) Assume that, around the arc $i$, $\tau$ looks like the configuration shown in Figure \ref{casethree}, with $l>1$, and none of $j$ and $k$ enclosing a
self-folded triangle.
% THIRD GENERIC CASE
        \begin{figure}[!h]
                \caption{Case \ref{genericthree}, configuration of $\tau$ around $i$}\label{casethree}
                \centering
                \includegraphics[scale=.35]{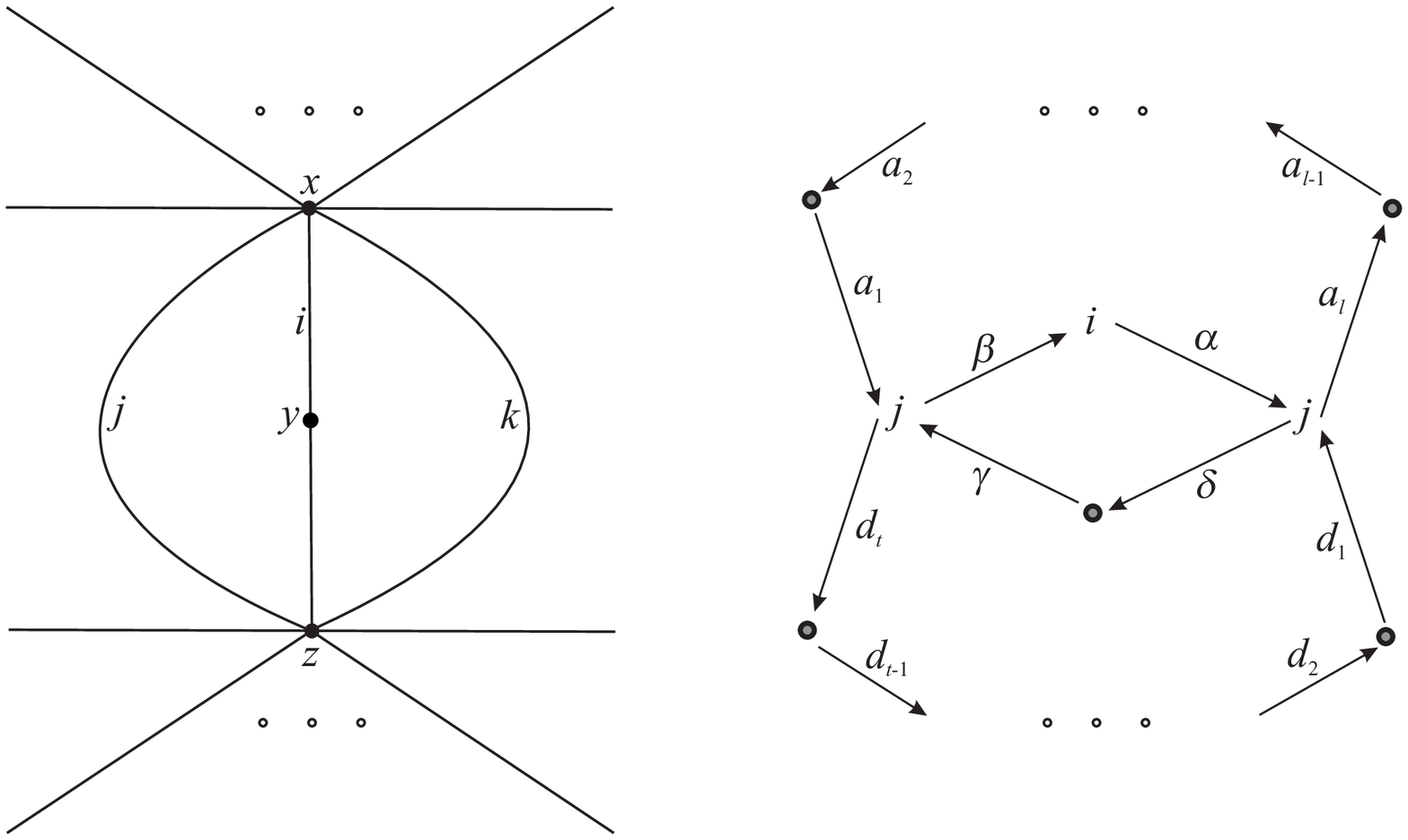}
        \end{figure}
Then
$$
\stau=-\frac{\alpha\beta\gamma\delta}{y}+x\alpha\beta a+z\gamma\delta d+\sptau,$$
with $\sptau$ involving none of the arrows $\alpha,\beta,\gamma,\delta$. If we perform the premutation $\premuti$ on $\astau$, we get $\tildeastau$, where $\tildeatau$ is the arrow span of the quiver shown in Figure \ref{premutcasethree},
% PREMUTATED THIRD GENERIC CASE
        \begin{figure}[!h]
                \caption{Case \ref{genericthree}, QP-mutation process $\muti\qstau$}\label{premutcasethree}
                \centering
                \includegraphics[scale=.35]{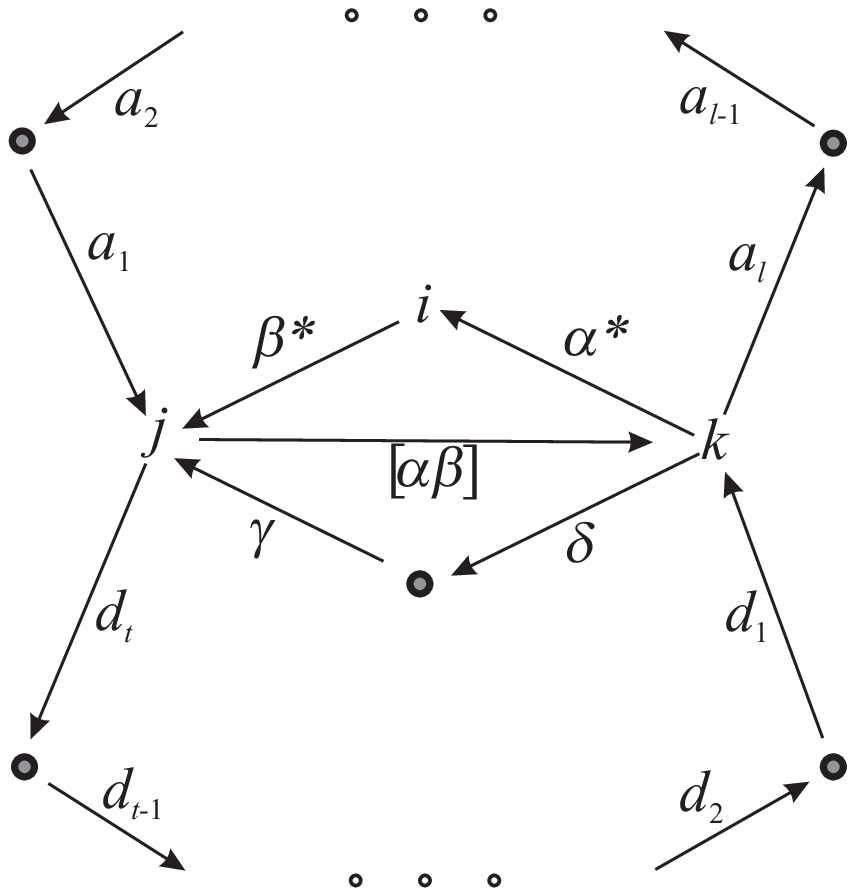}
        \end{figure}
and $\tildestau=-\frac{[\alpha\beta]\gamma\delta}{y}+x[\alpha\beta]a+z\gamma\delta d+\sptau+\beta^*\alpha^*[\alpha\beta]$. Since $\tildeastau$
is already reduced, it coincides with $\muti\astau$.

On the other hand, $\sigma=f_i(\tau)$ and its quiver $\qsigma$ look as Figure \ref{flippedcasethree},
% FLIPPED THIRD GENERIC CASE
        \begin{figure}[!h]
                \caption{Case \ref{genericthree}, flip $\sigma=f_i(\tau)$}\label{flippedcasethree}
                \centering
                \includegraphics[scale=.35]{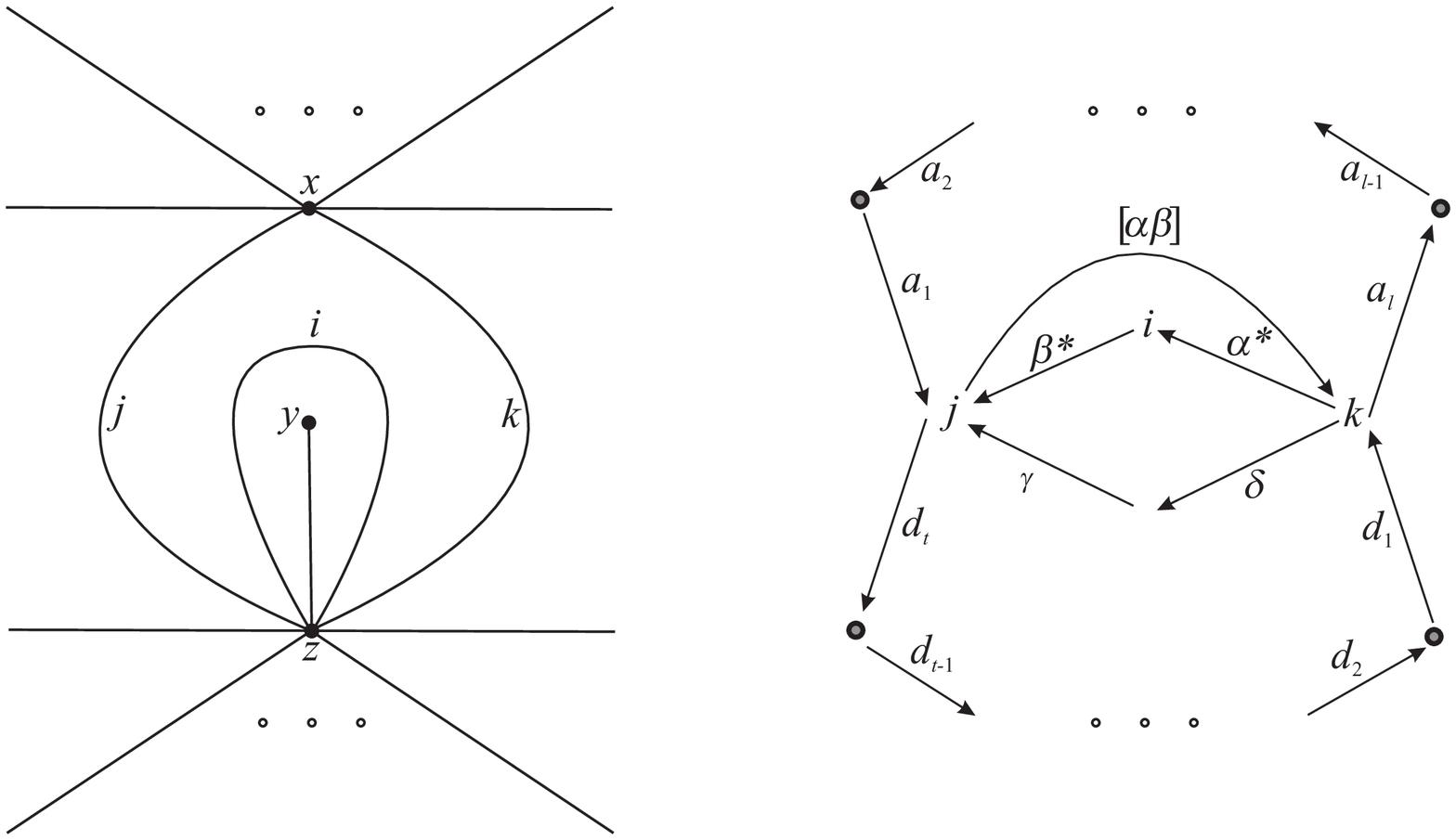}
        \end{figure}
and $\ssigma=\beta^*\alpha^*[\alpha\beta]-\frac{[\alpha\beta]\gamma\delta}{y}+x[\alpha\beta]a+z\gamma\delta d+\spsigma$, with
$\spsigma=\sptau$. Therefore, $\muti\astau=\assigma$.
\end{case}

\begin{case}\label{new1} (Flip inside the fourth puzzle piece of Figure \ref{puzzlepieces}) Assume that, around the arc $i$, $\tau$ looks like the configuration in Figure \ref{newcaseone}.
% NEW FIRST GENERIC CASE
        \begin{figure}[!h]
                \caption{Case \ref{new1}, configuration of $\tau$ around $i$}\label{newcaseone}
                \centering
                \includegraphics[scale=.35]{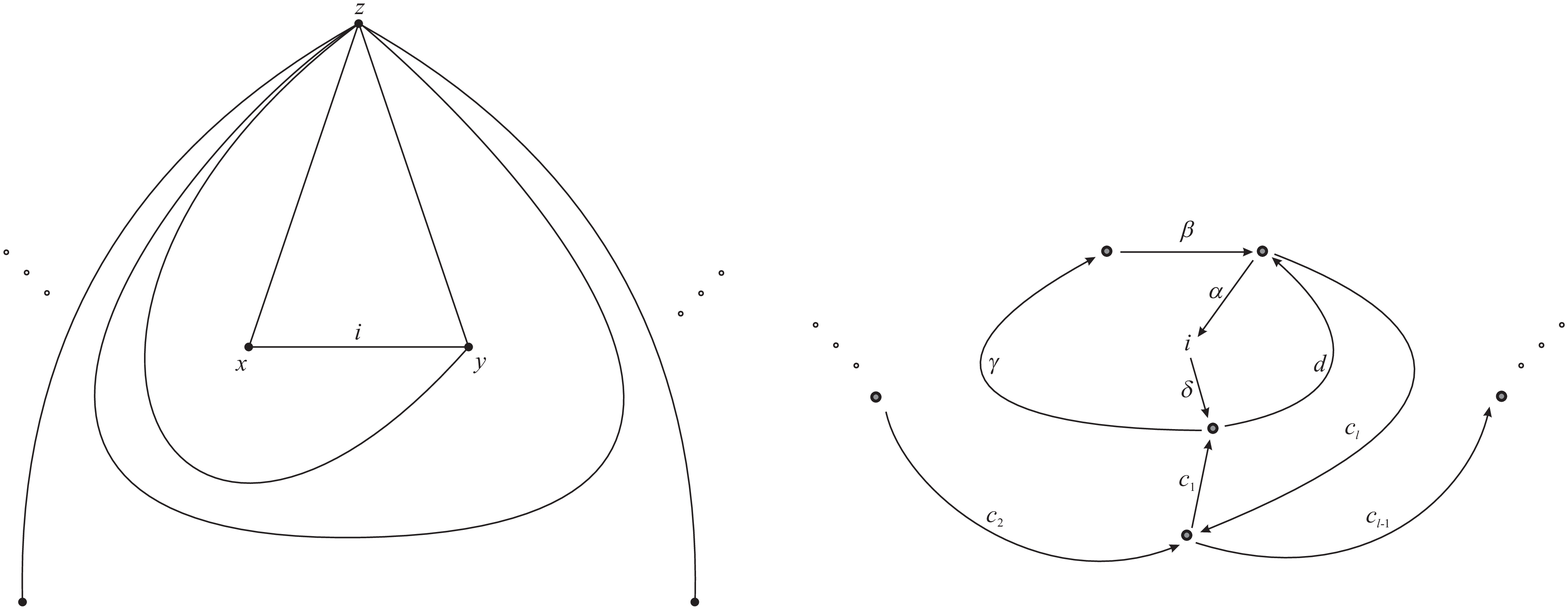}
        \end{figure}
Let us abbreviate $c=c_1\ldots c_l$. Then
$$
\stau=y\delta\alpha d+dc_1c_l-\frac{\delta\alpha\beta\gamma}{x}+zc\beta\gamma+\sptau
$$
with $\sptau\in\ratau$ involving none of the arrows $\alpha,\beta,\gamma,\delta,c_1,c_l,d$. If we perform the premutation $\premuti$ on $\astau$, we get $\tildeastau$, where $\tildeatau$ is the arrow span of the quiver shown in Figure \ref{newpremutcaseone}
% NEW PREMUTATED FIRST GENERIC CASE
        \begin{figure}[!h]
                \caption{Case \ref{new1}, QP-mutation process $\muti\qstau$}\label{newpremutcaseone}
                \centering
                \includegraphics[scale=.35]{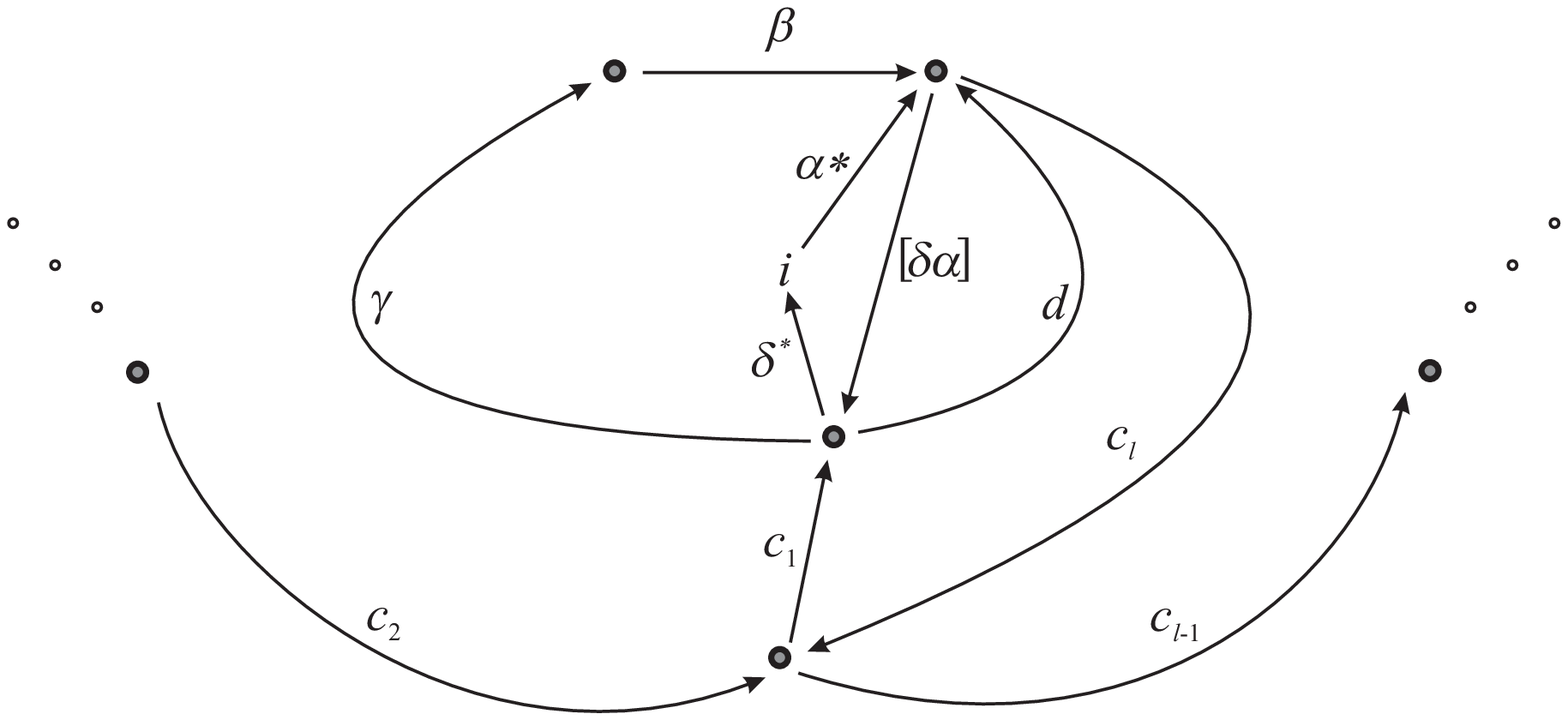}
        \end{figure}
and $\tildestau=y[\delta\alpha]d+dc_1c_l-\frac{[\delta\alpha]\beta\gamma}{x}+zc\beta\gamma+\sptau+[\delta\alpha]\alpha^*\delta^*\in
\rtildeatau$. The $R$-algebra automorphism $\varphi$ of $\rtildeatau$ whose action on the arrows is given by
$$
[\delta\alpha]\mapsto[\delta\alpha]-\frac{c_1c_l}{y}, \ d\mapsto d+\frac{\beta\gamma}{xy}-\frac{\alpha^*\delta^*}{y},
$$
and the identity in the rest of the arrows, sends $\tildestau$ to
$$
\varphi(\tildestau)=y[\delta\alpha]d+\frac{c_1c_l\beta\gamma}{xy}-\frac{c_1c_l\alpha^*\delta^*}{y}+zc\beta\gamma+\sptau.
$$
Therefore, the reduced part $\muti\astau$ of $(\tildeatau,\varphi(\tildestau))$ is (up to right-equivalence) the QP on the arrow span $\overatau$ obtained from $\tildeatau$ by deleting the arrows $[\delta\alpha]$ and $d$, with $\varphi(\tildestau)-y[\delta\alpha]d$ as its potential.

On the other hand, $\sigma=f_i(\tau)$ and its quiver $\qsigma$ look as Figure \ref{newflippedcaseone},
% NEW FLIPPED FIRST GENERIC CASE
        \begin{figure}[!h]
                \caption{Case \ref{new1}, flip $\sigma=f_i(\tau)$}\label{newflippedcaseone}
                \centering
                \includegraphics[scale=.35]{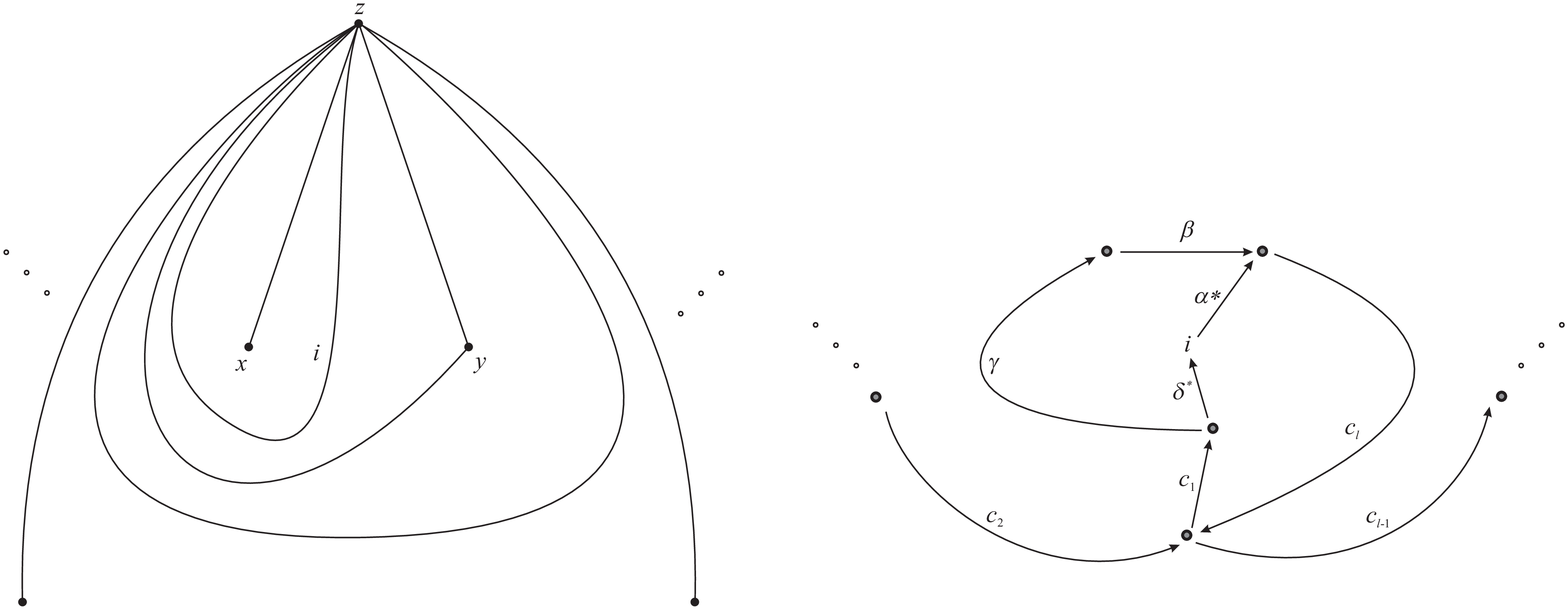}
        \end{figure}
and $\ssigma=\frac{c_1c_l\beta\gamma}{xy}-\frac{c_1c_l\alpha^*\delta^*}{y}+zc\beta\gamma+\spsigma$, with $\spsigma=\sptau$. Thus the identity is
a right-equivalence between $\muti\astau$ and $\assigma$.
\end{case}

\begin{case}\label{new7} (Flip inside the third puzzle piece of Figure \ref{puzzlepieces}) Assume that, around the arc $i$, $\tau$ looks like the configuration in Figure \ref{newcaseseven}.
% FIRST GENERIC CASE
        \begin{figure}[!h]
                \caption{Case \ref{new7}, configuration of $\tau$ around $i$}\label{newcaseseven}
                \centering
                \includegraphics[scale=.3]{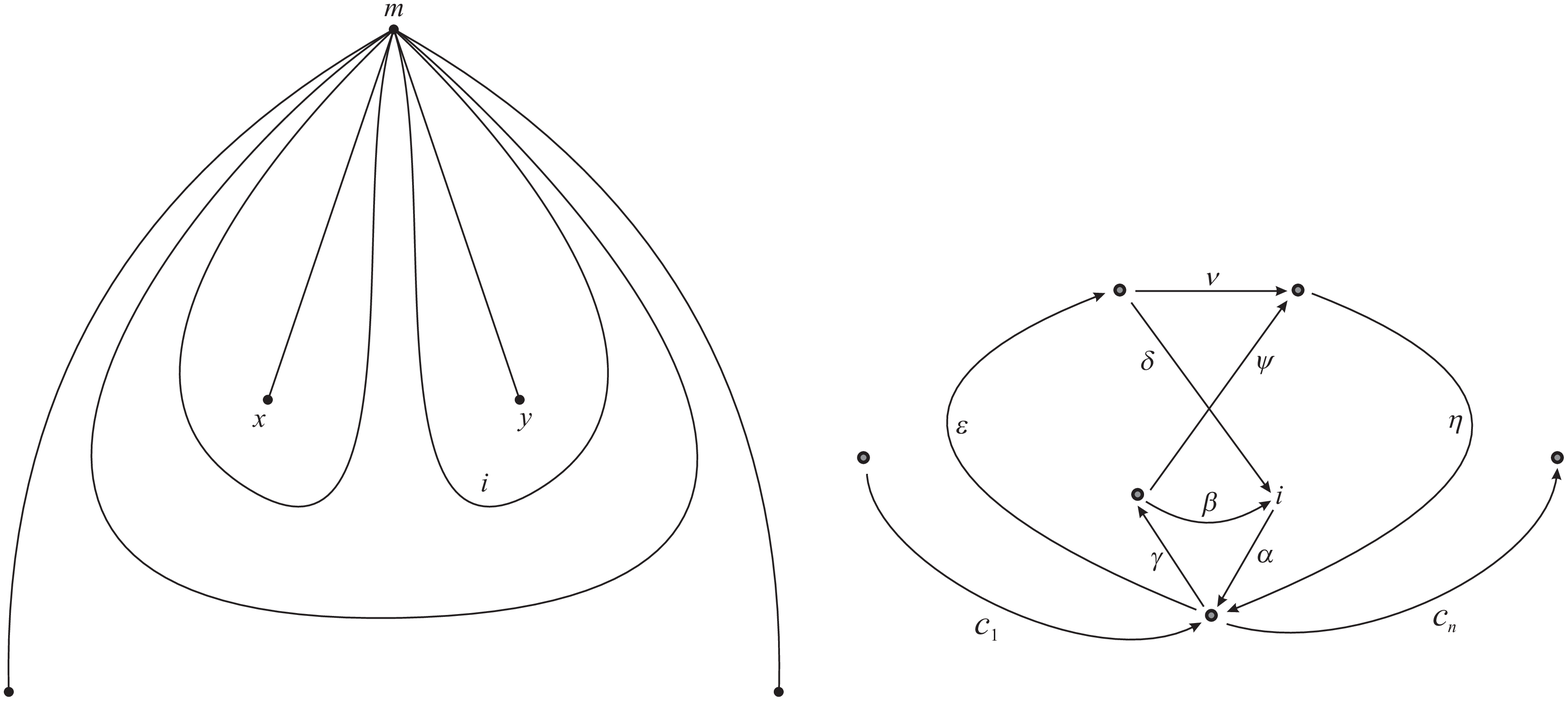}
        \end{figure}
Let us abbreviate $c=c_1\ldots c_n$. Then
$$
\stau=\alpha\beta\gamma-\frac{\alpha\delta\varepsilon}{x}-\frac{\eta\psi\gamma}{y}+\frac{\eta\nu\varepsilon}{xy}+wc\eta\nu\varepsilon+\sptau,
$$
with $\sptau\in\ratau$ involving none of the arrows $\alpha,\beta,\gamma,\delta,\varepsilon,\nu,\psi$. If we perform the premutation $\premuti$ on $\astau$, we get $\tildeastau$, where $\tildeatau$ is the arrow span of the quiver shown in Figure \ref{newpremutcasefive}
% PREMUTATED FIRST GENERIC CASE
        \begin{figure}[!h]
                \caption{Case \ref{new7}, QP-mutation process $\muti\qstau$}\label{newpremutcaseseven}
                \centering
                \includegraphics[scale=.3]{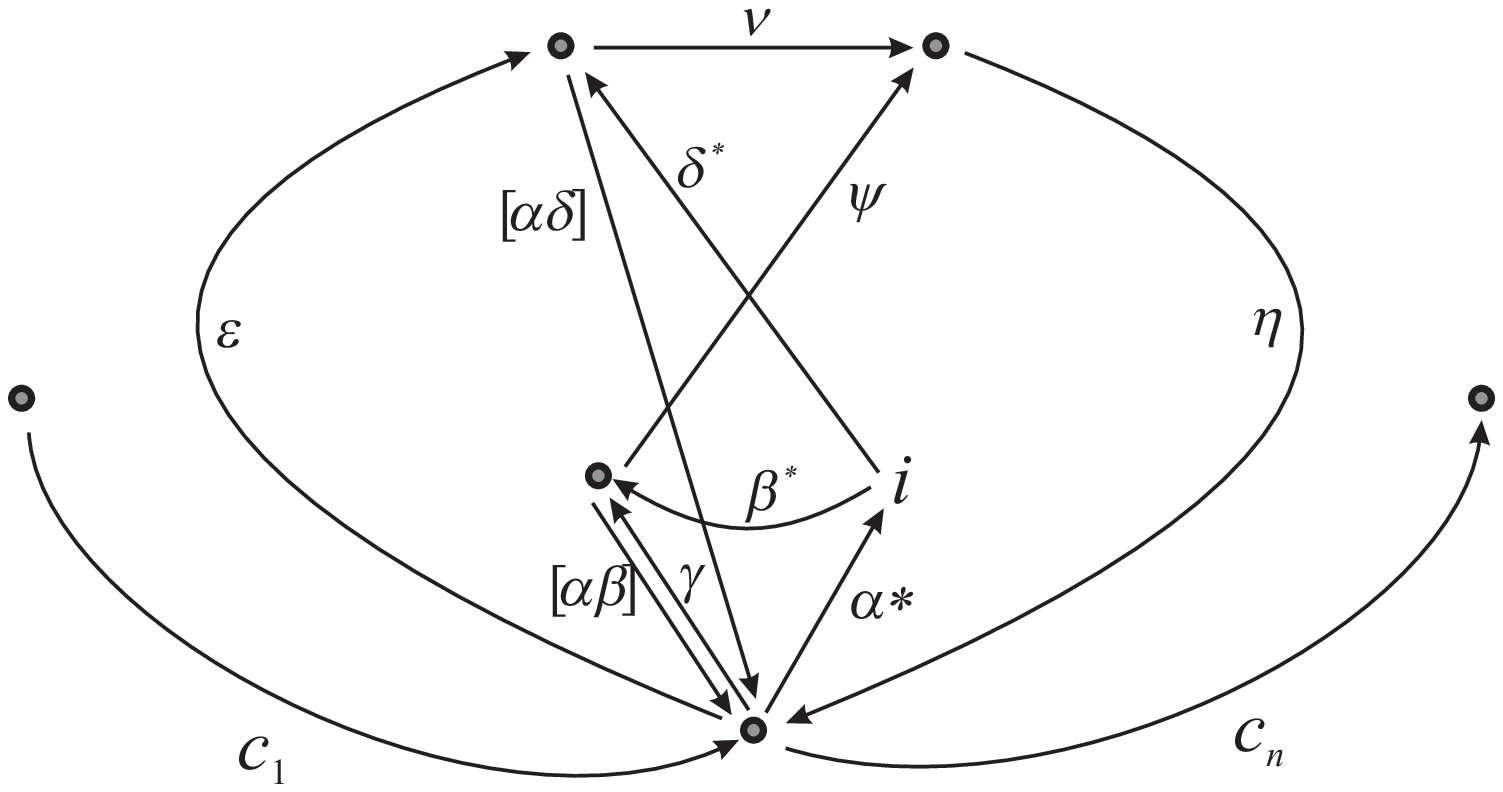}
        \end{figure}
and $\tildestau=[\alpha\beta]\gamma-\frac{[\alpha\delta]\varepsilon}{x}-\frac{\eta\psi\gamma}{y}+\frac{\eta\nu\varepsilon}{xy}+wc\eta\nu\varepsilon+
\sptau+[\alpha\beta]\beta^*\alpha+[\alpha\delta]\delta^*\alpha^*\in\rtildeatau$. The $R$-algebra automorphism $\varphi$ of $\rtildeatau$ whose action on the arrows is given by
$$
[\alpha\beta]\mapsto[\alpha\beta]+\frac{\eta\psi}{y}, \ \gamma\mapsto\gamma-\beta^*\alpha^*, \ [\alpha\delta]\mapsto[\alpha\delta]+\frac{\eta\nu}{y}+wxc\eta\nu, \ \varepsilon\mapsto\varepsilon+x\delta^+\alpha^*,
$$
and the identity in the rest of the arrows, sends $\tildestau$ to
$$
\varphi(\tildestau)=[\alpha\beta]\gamma-\frac{[\alpha\delta]\varepsilon}{x}+\frac{\eta\psi\beta^*\alpha^*}{y}+\frac{\eta\nu\delta^*\alpha^*}{y}
+wxc\eta\nu\delta^*\alpha^*+\sptau.
$$
Therefore, the reduced part $\muti\astau$ of $(\tildeatau,\varphi(\tildestau))$ is (up to right-equivalence) the QP on the arrow span
$\overatau$ obtained from $\tildeatau$ by deleting the arrows $[\alpha\beta],\gamma,[\alpha\delta]$ and $\varepsilon$, with $\varphi(\tildestau)-[\alpha\beta]\gamma+\frac{[\alpha\delta]\varepsilon}{x}$ as its potential.

On the other hand, $\sigma=f_i(\tau)$ and its quiver $\qsigma$ look as Figure \ref{newflippedcaseseven},
% FLIPPED FIRST GENERIC CASE
        \begin{figure}[!h]
                \caption{Case \ref{new7}, flip $\sigma=f_i(\tau)$}\label{newflippedcaseseven}
                \centering
                \includegraphics[scale=.3]{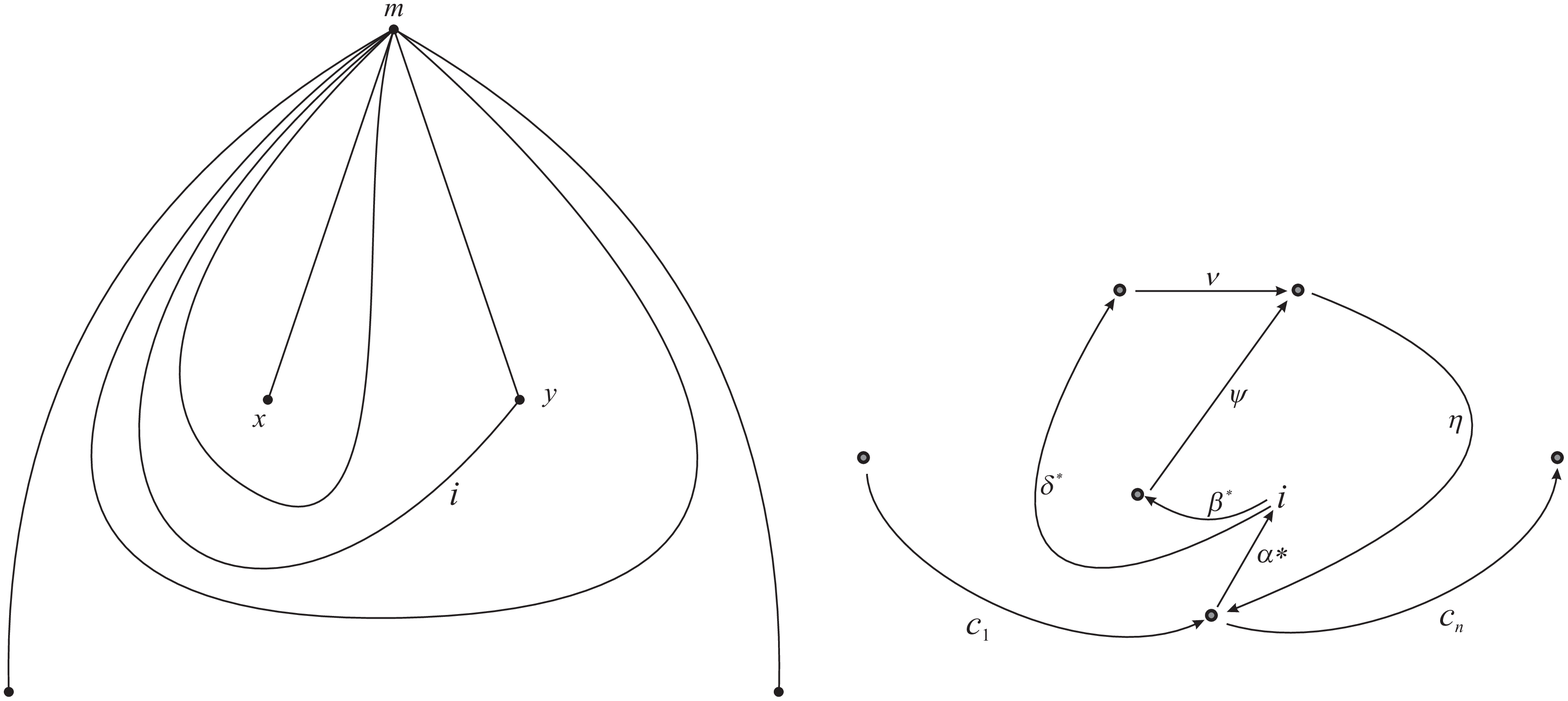}
        \end{figure}
and $\ssigma=-\frac{\alpha^*\eta\psi\beta^*}{y}+\frac{\alpha^*\eta\nu\delta^*}{xy}+w\eta\nu\delta^*\alpha^*c+\spsigma$, with $\spsigma=\sptau$. Thus the $R$-algebra isomorphism
$\psi:\roveratau\rightarrow\rasigma$ whose action on the arrows is given by
$$
\alpha^*\mapsto-\alpha^*, \ \delta^*\mapsto-\frac{\delta^*}{x},
$$
and the identity in the rest of the arrows, is a right-equivalence between $\muti\astau$ and $\assigma$.
\end{case}

A similar analysis in the rest of the cases (see Remark \ref{samemarkedpoint}) finishes the proof of Theorem \ref{flip<->mutation}.
\end{proof}

\section{Non-empty boundary: Rigidity and finite dimension}\label{rigidityfornonemptyboundary}

Our second main result ensures non-degeneracy for the QPs constructed in Definition \ref{QPfortriangulation} provided the boundary of the underlying surface is non-empty.

\begin{teo}\label{noboundrigid} Let $\surf$ be a bordered surface with marked points. If $\surfnoM$ has non-empty boundary, then the QP associated in Definition
\ref{QPfortriangulation} to any ideal triangulation of $\surf$ (under any choice $(x_p)_{p\in\punct}$) is rigid, hence non-degenerate.
\end{teo}

\begin{proof} Changing the notation a little bit, throughout the proof we will assume that $\surf$ has no punctures, in other words, all the marked points in $M$ belong to the boundary. We begin by inductively constructing a sequence of triangulations $\sigma_1,\sigma_2,\ldots$, such that
each $\sigma_n$ is an ideal triangulation of $\surfpn$ for some set $\punctn=\{p_1,\ldots,p_n\}$ of $n$ distinct punctures on $\surfnoM$ and
\begin{equation}\label{qinclusions} \tau\subseteq\sigma_1\subseteq\sigma_2\subseteq\ldots \text{ and }
\qtau\subseteq\qsigmaone\subseteq\qsigmatwo\subseteq\ldots
\end{equation}

Let $\tau$ be any ideal triangulation of $\surf$ (which we assume to have no punctures). To construct $\sigma_1$ we need to choose a puncture in $\surfnoM$. We put $p_1$ inside any non-interior triangle $\triangle_0$ of $\tau$. Then we draw the three arcs emanating from $p_1$ and going to the three
vertices of $\triangle_0$. The result is an ideal triangulation $\sigma_1$ of $\surfpone$. For $n>1$, once $\sigma_{n-1}$ has been constructed, we put $p_n$ inside a non-interior
triangle $\triangle_{n-1}$ of $\sigma_{n-1}$ having $p_{n-1}$ as a vertex, then we draw the three arcs emanating from $p_n$ and going to the
three vertices of $\triangle_{n-1}$. The result is an ideal triangulation of $\surfpn$.

Now, for $n>0$, let $\tau_n=f_{j_1^n}(\sigma_n)$ be the triangulation obtained by the flip of the arc $j_1^n$ of the triangulation $\sigma_n$, see Figure \ref{foneftwosigma}. Note that $\sigma_{n-1}\subseteq\tau_n$ and $\qsigmanminone\subseteq\qtaun$ for $n\geq 0$ if we denote $\sigma_0=\tau_0=\tau$.
% DEFINITION of TAU_n
        \begin{figure}[!h]
                \caption{$\sigma_n$}\label{foneftwosigma}
                \centering
                \includegraphics[scale=.5]{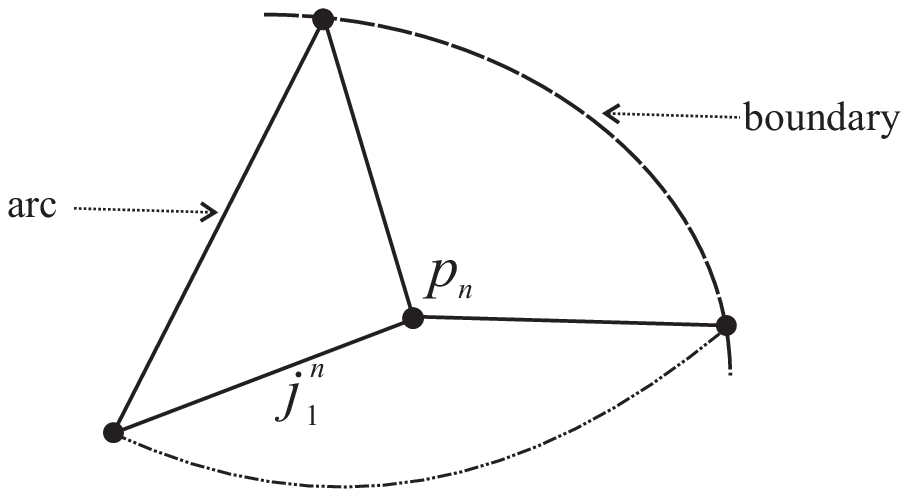}
        \end{figure}

\begin{lema}\label{jacidealspreserved} With the above notation, every potential on $\qsigmanminone$ belonging to $J(S(\sigma_{n-1}))$ is cyclically equivalent to an
element of $J(S(\tau_n))$.
\end{lema}
\begin{proof} For $n=1$ we actually have $J(S(\tau))\subseteq J(S(\tau_1))$. So let us treat the case $n>1$. With the notation of Figure
\ref{indsteprigid},
% INDUCTIVE STEP FOR RIGIDITY
        \begin{figure}[!h]
                \caption{Proving rigidity}\label{indsteprigid}
                \centering
                \includegraphics[scale=.65]{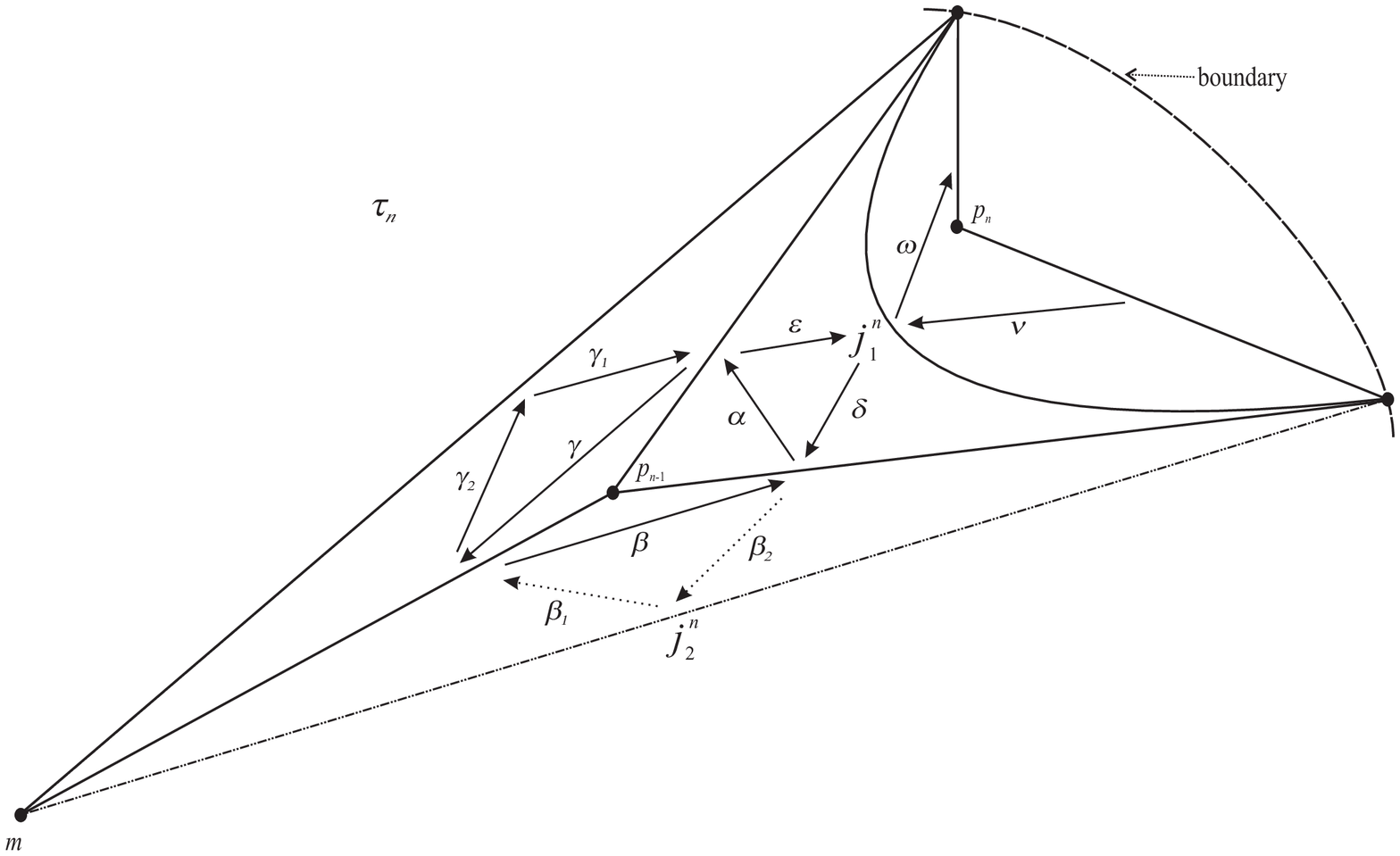}
        \end{figure}
we have $S(\tau_n)=S(\sigma_{n-1})+\alpha\delta\varepsilon$, and hence, for $a\in Q_1(\sigma_{n-1})$, $a\neq\alpha$ (note that $\delta,\varepsilon\notin Q_1(\sigma_{n-1})$), we have
$\partial_a(S(\tau_n))=\partial_a(S(\sigma_{n-1}))$, whereas
$\partial_\alpha(S(\tau_n))=\partial_a(S(\sigma_{n-1}))+\delta\varepsilon=x_{p_{n-1}}\beta\gamma+\delta\varepsilon$,
$\partial_\delta(S(\tau_n))=\varepsilon\alpha$ and $\partial_\varepsilon(S(\tau_n))=\alpha\delta$.

Let $c$ be a cycle of $\qsigmanminone$ that has $\beta\gamma$ as factor. Then one of the following three cases holds:
\begin{itemize}
\item $c$ is cyclically equivalent to a cycle in $\qsigmanminone$ that has $\beta\gamma\alpha$ as factor.
\item $c$ is cyclically equivalent to a cycle in $\qsigmanminone$ that has $\alpha\beta\gamma$ as factor.
\item The segment labeled $j_2^n$ in Figure \ref{indsteprigid} is an arc of $\sigma_{n-1}$ (hence the dotted arrows are indeed arrows of
$\sigma_{n-1}$) and $c$ is cyclically equivalent to a cycle that has $\beta_2\beta\gamma\gamma_1$ as factor.
\end{itemize}
Since $\beta\gamma\alpha=x_{p_{n-1}}^{-1}\partial_\alpha(S(\tau_n))\alpha-x_{p_{n-1}}^{-1}\delta\varepsilon\alpha\in J(S(\tau_n))$ and
$\alpha\beta\gamma=x_{p_{n-1}}^{-1}\alpha\partial_\alpha(S(\tau_n))-x_{p_{n-1}}^{-1}\alpha\delta\varepsilon\in J(S(\tau_n))$, in the first two
cases we see that $c$ is cyclically equivalent to an element of $J(S(\tau_n))$.

In the third case, we have two possibilities: The marked point $m$ either lies on the boundary of $\surfnoM$ or is a puncture of
$\surfpnminone$. If it lies on the boundary, then
$\beta_2\beta\gamma\gamma_1=\beta_2\beta\partial_{\gamma_2}(S(\sigma_{n-1}))=\beta_2\beta\partial_{\gamma_2}(S(\tau_n))\in J(S(\tau_n))$.
Otherwise, if $m$ is a puncture, then $\beta_2\beta\gamma\gamma_1=\beta_2\beta\partial_{\gamma_2}(S(\sigma_{n-1}))-x_m\beta_2\beta\beta_1b=
\beta_2\beta\partial_{\gamma_2}(S(\tau_n))-x_m\beta_2\partial_{\beta_2}(S(\tau_n))b\in J(S(\tau_n))$ for a unique path $b$ in $Q(\sigma_{n-1})$
from the arc $t(\gamma_1)$ to the arc $t(\beta_1)$.

This and the fact that $\partial_a(S(\tau_n))=\partial_a(S(\sigma_{n-1}))$ for $a\in Q_1(\sigma_{n-1})$ imply Lemma \ref{jacidealspreserved}.
\end{proof}

We now return to the proof of Theorem \ref{noboundrigid}. Since rigidity is preserved by QP-mutation, by Proposition \ref{seqofflips} and Theorem \ref{flip<->mutation} it suffices to show that $\qstaun$ is rigid. We prove this by induction on $n\geq 0$.

For $n=0$ we have $\tau_0=\tau$. Each arrow of $\qtau$ appears in at most one term of
$\stau$, and all the terms of $\stau$ are oriented triangles coming from interior triangles of $\tau$. Also, since there are no punctures,
every cycle in $\qtau$ is cyclically equivalent to a cycle that has $ab=\partial_c(S)$ as factor for some oriented cycle $abc$ that appears (up
to cyclical equivalence) as a term of $\stau$. Therefore, $\qstau$ is rigid.

For the inductive step, let $n>0$ and assume that the QP $(Q(\tau_{n-1}),S(\tau_{n-1}))$ associated to the
triangulation $\tau_{n-1}$ of $\surfpnminone$ is rigid; then the QP $\qssigmanminone$ is rigid as well. Take any cycle $c$ in $\qtaun$. If $c$ is contained in $Q(\sigma_{n-1})$, then
$c$ is cyclically equivalent to an element of $J(S(\sigma_{n-1}))$, which in turn is cyclically equivalent to an element of $J(S(\tau_n))$ by
Lemma \ref{jacidealspreserved}. If $c$ is not contained in $Q(\sigma_{n-1})$, then $c$ is cyclically equivalent to a cycle $c'$ that has
$\delta\varepsilon$ as factor, say $c'=d\delta\varepsilon$, with $d$ a path in $Q(\tau_n)$ from $h(\delta)$ to $t(\varepsilon)$. Then
$c'=d\partial_\alpha(S(\tau_n))-x_{p_{n-1}}d\beta\gamma$, and we can keep substituting each factor $\delta\varepsilon$ of $d$ by
$\partial_\alpha(S(\tau_n))-x_{p_{n-1}}\beta\gamma$. After doing this, we see that $c'$ is the sum of an element of $J(S(\tau_n))$ with a scalar
multiple of a cycle $c'' $ in $Q(\sigma_{n-1})$. This cycle $c''$ is cyclically equivalent to an element of $J(S(\tau_n))$ as we have seen
above. We conclude that $\qstaun$ is rigid. This finishes the proof of Theorem \ref{noboundrigid}.
\end{proof}

\begin{obs}\begin{enumerate}\item In terms of trace spaces (see \cite{DWZ}, Definition 3.4), Lemma \ref{jacidealspreserved} says that the inclusion of quivers $\qsigmanminone\hookrightarrow\qtaun$ induces a well defined map between the trace spaces $\trspace(\mathcal{P}(Q(\sigma_{n-1}),S(\sigma_{n-1})))\rightarrow\trspace(\mathcal{P}(Q(\tau_n),S(\tau_n)))$.
\item Since every quiver mutation equivalent to a quiver of the form $Q(\tau)$, with $\tau$ an ideal triangulation of $\surf$, is of the same form, Theorem \ref{noboundrigid} says in particular that in Definition \ref{QPfortriangulation} we have given an explicit construction of a non-degenerate potential for each of the members of the mutation equivalence class of the quivers that arise as signed adjacency quivers of triangulations of surfaces with non-empty boundary.
\end{enumerate}
\end{obs}

\begin{conj}\label{conjnondeg} The QP associated in Definition \ref{QPfortriangulation} to any triangulation of $\surf$ is always non-degenerate, regardless of the emptiness of the boundary of $\surfnoM$.
\end{conj}

\begin{conj}\label{conjnonrigid} If $\surf$ has empty boundary, then the QP associated to any triangulation of $\surf$ is never rigid. In other words, the converse of Theorem \ref{noboundrigid} is true as well.
\end{conj}

To illustrate this two conjectures, let us look at an example.

\begin{ex} Consider the ``canonical" triangulation $\tau$ of the once-punctured torus, see Figure \ref{torusqp}.
% TORUS
        \begin{figure}[!h]
                \caption{}\label{torusqp}
                \centering
                \includegraphics[scale=.65]{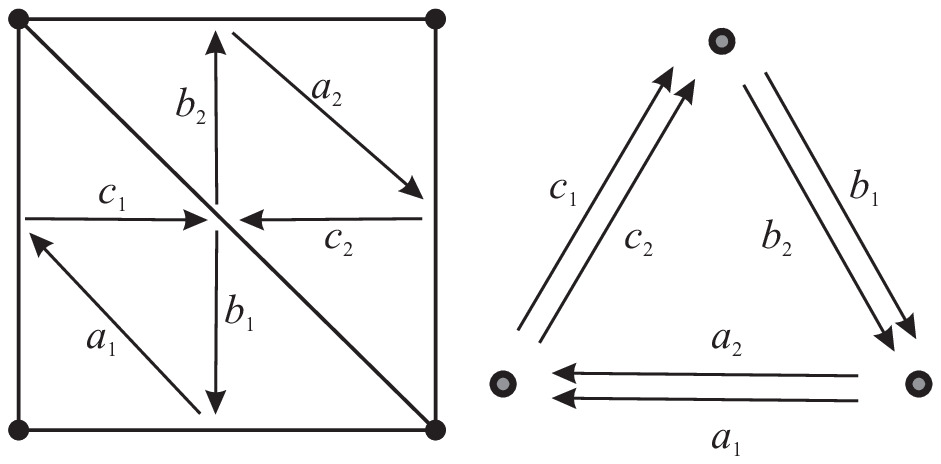}
        \end{figure}
We have $\stau=a_1b_1c_1+a_2b_2c_2+xa_1b_2c_1a_2b_1c_2$. If we mutate $\qstau$ in direction $i=t(b_1)$, we obtain $\overqstau$, where $\overqtau$ is the quiver shown in Figure \ref{torusqpflipped} and $\overstau=c^*_1b^*_2[b_2c_1]+c^*_2b^*_1[b_1c_2]+xc^*_1b^*_1[b_2c_1]c^*_2b^*_2[b_1c_2]$, and we therefore have $\qssigma=\muti\qstau$, where $\sigma=f_{i}(\tau)$.
% TORUS MUTATED
        \begin{figure}[!h]
                \caption{}\label{torusqpflipped}
                \centering
                \includegraphics[scale=.65]{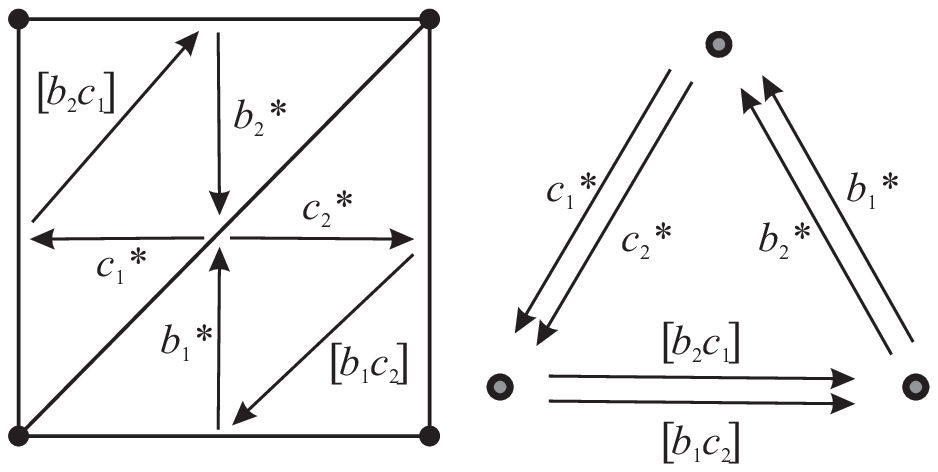}
        \end{figure}
This example shows in particular that Theorem \ref{flip<->mutation} holds for the once-punctured torus and that Definition \ref{QPfortriangulation} gives an explicit non-degenerate potential for the ``double cyclic triangle" $\qtau$. However, it is known (see \cite{DWZ}, example 8.6) that the double cyclic triangle does not admit a rigid potential. In short words, Conjectures \ref{conjnondeg} and \ref{conjnonrigid} hold for the once-punctured torus.
\end{ex}

To finish this section and close this paper, we show that non-empty boundary implies also finite dimension of the Jacobian algebra.

\begin{teo}\label{finitedimensional} If the surface $\surf$ has non-empty boundary, then for any triangulation $\tau$ of $\surf$ the Jacobian algebra $\mathcal{P}\qstau$ is finite dimensional, the ideal $I(\tau)$ of the path algebra $R\langle\qtau\rangle$ generated by $\{\partial_a(\stau)\suchthat a\in Q_1(\tau)\}$ is admissible (that is, it is contained in the square of the ideal generated by the arrows and contains all paths of sufficiently large length), and $\mathcal{P}\qstau$ is isomorphic to the quotient $R\langle\qtau\rangle/I(\tau)$.
\end{teo}

\begin{proof} Assume that $\surf$ has non-empty boundary and no punctures, and let $\tau$ be any triangulation of $\surf$. For each marked point $m$ we have the following configuration in a sufficiently small neighborhood of $m$:
% SMALL NEIGHBORHOOD OF m
        \begin{figure}[!h]
                \caption{}\label{smallneigh}
                \centering
                \includegraphics[scale=.65]{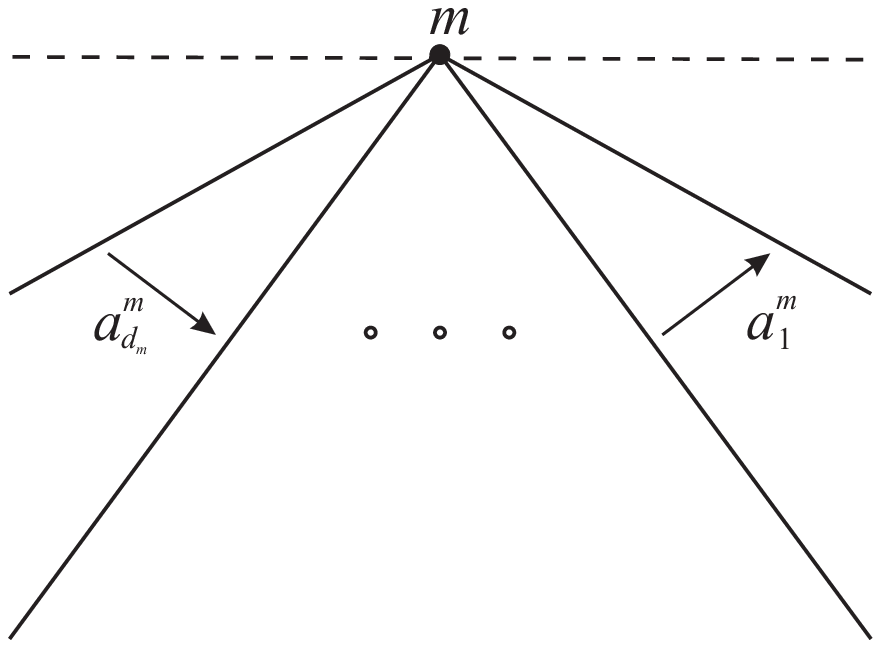}
        \end{figure}\\
where the arrows $a_1^m,\ldots,a_{d_m}^m$ are uniquely determined by $m$. Let $N=\max\{d_m\suchthat m\in M\}$. Then any path in $\qtau$ having length greater than $N$ must have $ab=\partial_c(\stau)$ as a factor for some oriented triangle $abc$ appearing as a term of $\stau$.

Now let $\tau_1$ be the triangulation of $\surfpone$ constructed in the proof of Theorem \ref{noboundrigid} (see Figure \ref{tauone}). Then we have $J(\tau)\subseteq J(\tau_1)$ and $\alpha\delta,\delta\varepsilon,\varepsilon\alpha\in J(\tau_1)$. Therefore any path in $Q(\tau_1)$ of length greater than $N+2$ belongs to $J(\tau_1)$.
% SURFACE WITH BOUNDARY AND ONE PUNCTURE
        \begin{figure}[!h]
                \caption{}\label{tauone}
                \centering
                \includegraphics[scale=.65]{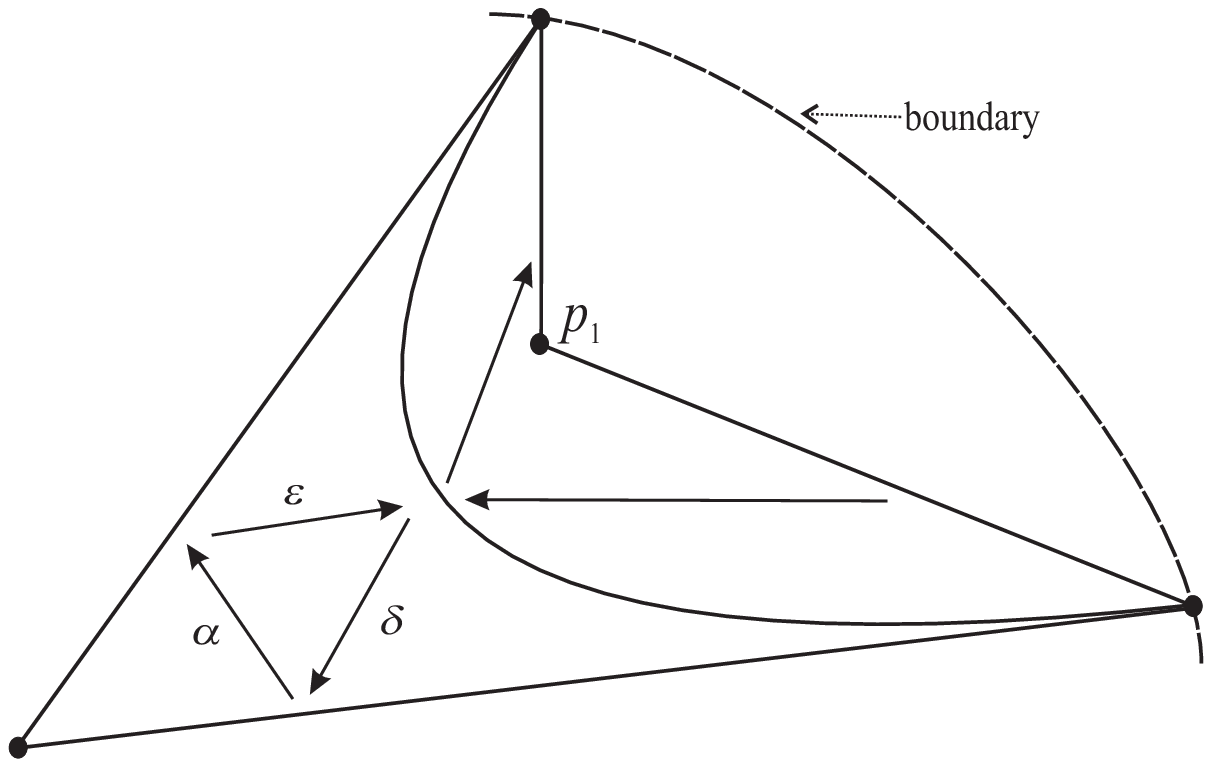}
        \end{figure}
Since finite-dimensionality of the Jacobian algebra is preserved under QP-mutation, we have proved the Theorem for non-empty boundary and at most one puncture.

The proof for the case of more than one puncture is similar to the proof of Lemma \ref{jacidealspreserved}. So let $n>1$ and, with the notation of Theorem \ref{noboundrigid}, assume inductively that the Jacobian algebra $\mathcal{P}(Q(\tau_{n-1}),S(\tau_{n-1}))$ is finite-dimensional. Then $\mathcal{P}(Q(\sigma_{n-1}),S(\sigma_{n-1}))$ is finite-dimensional as well.

With the notation of Figure \ref{indsteprigid}, let $\mathrm{P}$ denote the set of all paths in $Q(\sigma_{n-1})$ that do not start at $t(\gamma)$ and do not end at $h(\beta)$. If $u\in\mathrm{P}$ has $\beta\gamma$ as a factor, that is, if $u=u_1\beta\gamma u_2$ for some paths $u_1,u_2$ in $Q(\sigma_{n-1})$, then one of the following three conditions holds:
\begin{itemize}
\item $u$ has $\beta\gamma\alpha$ as factor.
\item $u$ has $\alpha\beta\gamma$ as factor.
\item The segment labeled $j_2^n$ in Figure \ref{indsteprigid} is an arc of $\sigma_{n-1}$ (hence the dotted arrows are indeed arrows of
$Q(\sigma_{n-1})$) and $u$ has $\beta_2\beta\gamma\gamma_1$ as factor.
\end{itemize}
In any of these three cases, we have $u\in J(S(\tau_{n}))$ just as in the proof of Lemma \ref{jacidealspreserved}.

Since $\mathcal{P}(Q(\sigma_{n-1}),S(\sigma_{n-1}))$ is finite-dimensional, there exists a positive integer $N'$ such that every element of $\mathrm{P}$ whose length is greater than $N'$ belongs to the Jacobian ideal $J(S(\sigma_{n-1}))$. Let $\mathrm{P}_{>N'}$ be the subset of $\mathrm{P}$ consisting of paths of length greater than $N'$. By the previous paragraph and because
$\partial_a(S(\tau_n))=\partial_a(S(\sigma_{n-1}))$ for $a\in Q_1(\sigma_{n-1})$, $a\neq\alpha$ (note that $\delta,\varepsilon\notin Q_1(\sigma_{n-1})$), we have $\mathrm{P}_{>N'}\subseteq J(S(\tau_n))$.

Now take any path $u$ in $Q(\sigma_{n-1})$ of length greater than $N'+2$. If $u$ does not start at $t(\gamma)$ and does not end at $h(\beta)$, then $u\in J(S(\tau_n))$. Otherwise, we can write $u=\beta u'\gamma$ for some path $u'$ in $Q(\sigma_{n-1})$ (remember that $\delta,\varepsilon\notin Q(\sigma_{n-1})$). This path $u'$ has length greater than $N'$, does not start at $t(\gamma)$ and does not end at $h(\beta)$, and hence belongs to $J(S(\tau_{n}))$.

Finally, let $v$ be any path in $Q(\tau_n)$ of length greater than $N'+6$, we claim that $v\in J(S(\tau_n))$. Since $\partial_\alpha(S(\tau_n))=x_{p_{n-1}}\beta\gamma+\delta\varepsilon$, we can assume, without loss of generality, that $v$ does not have $\delta\varepsilon$ as a factor. Then one of the following cases holds:
\begin{itemize}
\item $v$ is contained in $Q(\sigma_{n-1})$;
\item $\varepsilon$ is a factor of $v$;
\item $\delta$ is a factor of $v$.
\end{itemize}
In the first case, we have $v\in J(S(\tau_n))$. In the second case, we must have $v=\varepsilon v'$ or $v=\omega\varepsilon v'$ for some path $v'$ in $Q(\tau_n)$ because we are assuming that $v$ does not contain $\delta\varepsilon$ as factor. In any of these two situations, the path $v'$ is either contained in $Q(\sigma_{n-1})$ (and has length greater than $N'+4$) or can be written as $v'=v''\delta$ or $v'=v''\delta\nu$ for some path $v''$ (of length greater than $N'+2$) contained in $Q(\sigma_{n-1})$. This yields $v\in J(\tau_n)$.

Similarly, in the third case, we must have $v=v'\delta$ or $v=v'\delta\nu$ for some path $v'$ in $Q(\tau_n)$, and in any of these situations, the path is either contained in $Q(\sigma_{n-1})$ (and has length greater than $N'+4$) or can be written as $v'=\varepsilon v''$ or $v'=\omega\varepsilon v''$ for some path $v''$ (of length greater than $N'+2$) contained in $Q(\sigma_{n-1})$. This yields $v\in J(S(\tau_n))$.

Therefore, the Jacobian algebra $\mathcal{P}(Q(\tau_n),S(\tau_n))$ has finite dimension. The theorem follows from Proposition \ref{seqofflips}, Theorem \ref{flip<->mutation} and Proposition \ref{findimisinvariant}.
\end{proof}

\section*{Acknowledgments}

I am grateful to my advisor, Professor Andrei Zelevinsky, for his support, guidance and encouragement, without any of which this work would not
have taken place. To Sachin Gautam and Professors Michael Barot, Gordana Todorov and Jerzy Weyman for helpful discussions. I also thank Manuel Vargas, Belen Fragoso, and the Representation Theory group of Universidad Nacional Aut\'onoma de M\'exico, especially Professors Michael Barot, Jos\'e Antonio de la Pe\~na and Christof Geiss, for the support I have received from them during my graduate studies. Finally, I want to thank also Pietra Delgado-Escalante and Alexei D\'iaz-Vera for their valuable help when I was typing this draft, and the anonymous referee for a number of valuable suggestions.

\end{document}